%% file: main.tex
\documentclass[11pt]{amsart}

\usepackage[T1]{fontenc}
\usepackage[utf8]{inputenc}
\usepackage{lmodern}

\usepackage{layout}
\usepackage[american]{babel}
\usepackage[babel, final]{microtype} 
\usepackage{amssymb}
\usepackage{amscd}
\usepackage{mathtools} 
\usepackage[mathscr]{eucal}
\usepackage{textcomp}
\usepackage[all,knot,arc]{xy}

\usepackage{ifpdf}
\usepackage{xifthen}
\usepackage{comment}

\usepackage{chngcntr} 
\usepackage{enumitem} 
\usepackage{multicol} 

\ifpdf
  \usepackage[pdftex,%
    dvipsnames,%
  ]{xcolor}
  \usepackage[pdftex, final]{graphicx}
  \usepackage{caption}
  \usepackage{subcaption}
  \usepackage{pinlabel}
  \usepackage[pdftex,%
    letterpaper,%
    includehead,%
    includefoot,%
    nomarginpar,%
    lmargin=1in,%
    rmargin=1in,%
    tmargin=1in,%
    bmargin=1in,%
  ]{geometry}
  \usepackage[pdftex,%
    final,%
    colorlinks=true,
    linkcolor=NavyBlue,
    citecolor=NavyBlue,
    filecolor=NavyBlue,
    menucolor=NavyBlue,
    urlcolor=NavyBlue,
    bookmarks=true,
    bookmarksdepth=3,
    bookmarksnumbered=true,
    bookmarksopen=true,
    bookmarksopenlevel=2,%
  ]{hyperref}
  \hypersetup{
    pdftitle={On the tau invariants in instanton and monopole Floer theories},
    pdfauthor={Sudipta Ghosh, Zhenkun Li, C.-M. Michael Wong},
    pdfsubject={Knot theory, Floer theory},
    pdfkeywords={Monopole Floer homology, instanton Floer homology, tau}
  }
\else
  \usepackage[dvips,%
    dvipsnames,%
  ]{xcolor}
  \usepackage[dvips, final]{graphicx}
  \usepackage{caption}
  \usepackage{subcaption}
  \usepackage{pinlabel}
  \usepackage[dvips,%
    letterpaper,%
    includehead,%
    includefoot,%
    nomarginpar,%
    lmargin=1in,%
    rmargin=1in,%
    tmargin=1in,%
    bmargin=1in,%
  ]{geometry}
  \usepackage[ps2pdf,%
    final,%
    colorlinks=true,
    linkcolor=NavyBlue,
    citecolor=NavyBlue,
    filecolor=NavyBlue,
    menucolor=NavyBlue,
    urlcolor=NavyBlue,
    bookmarks=true,
    bookmarksdepth=3,
    bookmarksnumbered=true,
    bookmarksopen=true,
    bookmarksopenlevel=2,%
  ]{hyperref}
  \hypersetup{
    pdftitle={On the tau invariants in instanton and monopole Floer theories},
    pdfauthor={Sudipta Ghosh, Zhenkun Li, C.-M. Michael Wong},
    pdfsubject={Knot theory, Floer theory},
    pdfkeywords={Monopole Floer homology, instanton Floer homology, tau}
  }
\fi

\input{macros}

\input{commands}

\begin{document}

\title[On the tau invariants in instanton and monopole Floer theories]{On the 
  tau invariants \\ in instanton and monopole Floer theories}

\author[Sudipta Ghosh]{Sudipta Ghosh}
\address{Department of Mathematics \\ Louisiana State University \\ Baton Rouge, LA 70803}
\email{\href{mailto:sghos22@lsu.edu}{sghos22@lsu.edu}}
\urladdr{\url{https://sites.google.com/view/sudipta-ghosh/}}

\author[Zhenkun Li]{Zhenkun Li}
\address{Department of Mathematics \\ Stanford University \\ Stanford, CA 
  94305}
\email{\href{mailto:zhenkun@stanford.edu}{zhenkun@stanford.edu}}
\urladdr{\url{http://sites.google.com/view/zhenkun/}}

\author[C.-M. Michael Wong]{C.-M. Michael Wong}
\address{Department of Mathematics \\ Dartmouth College \\ Hanover, NH 03755}
\email{\href{mailto:won@math.dartmouth.edu}{wong@math.dartmouth.edu}}
\urladdr{\url{http://math.dartmouth.edu/~wong/}}

\begin{abstract}

  We unify two existing approaches to the \emph{tau} invariants in instanton 
  and monopole Floer theories, by identifying $\tauG$, defined by the second 
  author via the \emph{minus} flavors $\KHIm$ and $\KHMm$ of the knot 
  homologies, with $\taushG$, defined by Baldwin and Sivek via cobordism maps 
  of the $3$-manifold homologies induced by knot surgeries. We exhibit several 
  consequences, including a relationship with Heegaard Floer theory, and use 
  our result to compute $\KHIm$ and $\KHMm$ for twist knots.
  
\end{abstract}

\maketitle

\input{sec-intro}

\input{sec-prelim}

\input{sec-tau-conc}

\input{sec-tau-conn-sum}

\input{sec-identify}

\input{sec-computation}


\bibliographystyle{mwamsalphack}
\bibliography{Index}

\end{document}

%% file: commands.tex
\newcommand{\red}[1]{\textcolor{red}{#1}}

\newcommand{\paren}[1]{\mathchoice%
  {\left( #1 \right)}%
  {( #1 )}%
  {( #1 )}%
  {( #1 )}%
}

\newcommand{\abs}[1]{\mathchoice%
  {\left\lvert #1 \right\rvert}%
  {\lvert #1 \rvert}%
  {\lvert #1 \rvert}%
  {\lvert #1 \rvert}%
}
\newcommand{\set}[1]{\mathchoice%
  {\left\{ #1 \right\}}%
  {\{ #1 \}}%
  {\{ #1 \}}%
  {\{ #1 \}}%
}
\newcommand{\setc}[2]{\mathchoice%
  {\left\lbrace #1 \, \middle\vert \, #2 \right\rbrace}%
  {\lbrace #1 \, \vert \, #2 \rbrace}%
  {\lbrace #1 \, \vert \, #2 \rbrace}%
  {\lbrace #1 \, \vert \, #2 \rbrace}%
}

\newcommand*{\R}{\mathbb{R}}
\newcommand*{\Z}{\mathbb{Z}}
\newcommand*{\C}{\mathbb{C}}

\DeclareMathOperator{\Id}{Id}

\newcommand*{\ring}{\mathcal{R}}
\DeclareMathOperator{\rk}{rk}

\newcommand*{\Ftwo}{\mathbb{F}_2}
\newcommand*{\F}{\Ftwo}

\newcommand*{\isom}{\cong}
\newcommand*{\tensor}{\otimes}                
\newcommand*{\dirsum}{\oplus}                 
\newcommand*{\bigdirsum}{\bigoplus}           

\newcommand*{\connsum}{\mathbin{\sharp}} 
\newcommand*{\bigconnsum}{\mathop{\sharp}}

\newcommand*{\concgrp}{\mathcal{C}}

\newcommand*{\bdy}{\partial}
\newcommand*{\cross}{\times}
\newcommand*{\union}{\cup}
\newcommand*{\intersect}{\cap}

\DeclareMathOperator{\Int}{Int}
\renewcommand*{\Im}{\mathrm{Im}}

\newcommand*{\mir}[1]{\overline{#1}}

\DeclareMathOperator{\gr}{gr}

\DeclareMathOperator{\HFK}{HFK}
\DeclareMathOperator{\SFH}{SFH}
\DeclareMathOperator{\KHMrm}{KHM}
\DeclareMathOperator{\KHIrm}{KHI}
\DeclareMathOperator{\KHGrm}{KHG}
\DeclareMathOperator{\SHMrm}{SHM}
\DeclareMathOperator{\SHIrm}{SHI}
\DeclareMathOperator{\SHGrm}{SHG}

\newcommand*{\HFKm}{\HFK^-}
\newcommand*{\HFKh}{\widehat{\HFK}}

\newcommand*{\KHM}{\underline{\KHMrm}}
\newcommand*{\KHI}{\underline{\KHIrm}}
\newcommand*{\KHG}{\underline{\KHGrm}}

\newcommand*{\SHM}{\underline{\SHMrm}}
\newcommand*{\SHI}{\underline{\SHIrm}}
\newcommand*{\SHG}{\underline{\SHGrm}}

\newcommand*{\bSHG}{\underline{\mathbf{SHG}}}
\newcommand*{\bKHG}{\underline{\mathbf{KHG}}}

\newcommand*{\KHMm}{\KHM^-}
\newcommand*{\KHIm}{\KHI^-}
\newcommand*{\KHGm}{\KHG^-}

\DeclareMathOperator{\HF}{HF}
\newcommand*{\HFh}{\widehat{\HF}}

\DeclareMathOperator{\HM}{HM}
\newcommand*{\HMt}{\widetilde{\HM}}

\newcommand*{\tauH}{\tau_{\mathrm{H}}}
\newcommand*{\tauM}{\tau_{\mathrm{M}}}
\newcommand*{\tauI}{\tau_{\mathrm{I}}}
\newcommand*{\tauG}{\tau_{\mathrm{G}}}
\newcommand*{\taush}{\tau^\sharp}
\newcommand*{\taushH}{\taush_{\mathrm{H}}}
\newcommand*{\taushM}{\taush_{\mathrm{M}}}
\newcommand*{\taushI}{\taush_{\mathrm{I}}}
\newcommand*{\taushG}{\taush_{\mathrm{G}}}

\newcommand*{\nush}{\nu^\sharp}

\newcommand*{\nushI}{\nush_{\mathrm{I}}}
\newcommand*{\nushG}{\nush_{\mathrm{G}}}

\newcommand*{\Ish}{\mathrm{I}^\sharp}

\newcommand*{\xistd}{\xi_{\mathrm{std}}}
\newcommand*{\tb}{\mathit{tb}}
\newcommand*{\maxtb}{\overline{\tb}}
\newcommand*{\rot}{\mathit{r}}
\newcommand*{\slink}{\mathit{sl}}
\newcommand*{\leg}{\Lambda}
\newcommand*{\trans}{\Theta}

\newcommand*{\loss}{\mathcal{L}}
\newcommand*{\lossH}{\loss_{\mathrm{H}}}
\newcommand*{\lossM}{\loss_{\mathrm{M}}}

\newcommand*{\lossh}{\widehat{\loss}}
\newcommand*{\losshH}{\lossh_{\mathrm{H}}}
\newcommand*{\losshM}{\lossh_{\mathrm{M}}}

\newcommand{\al}{\alpha}
\newcommand{\vphi}{\varphi}
\newcommand{\be}{\beta}
\newcommand{\ga}{\gamma}
\newcommand{\Ga}{\Gamma}
\newcommand{\shm}{\SHM}

\newcommand{\ra}{\rightarrow}

\newcommand{\xra}{\xrightarrow}

\DeclareMathOperator{\Eig}{Eig}

\newcommand{\intg}{\mathbb{Z}}

\newcommand*{\khm}{\KHMm}
\newcommand*{\khi}{\KHIm}
\newcommand*{\shi}{\SHI}
\newcommand*{\shg}{\SHG}
\newcommand*{\khg}{\KHGm}

\newcommand*{\swb}[1]{\widetilde{#1}}

\DeclareMathOperator{\Spin}{Spin}
\newcommand*{\Spinc}{\Spin^c}
\DeclareMathOperator{\Tors}{Tors}

%% file: sec-intro.tex
\section{Introduction}
\label{sec:intro}

Among the Floer invariants of $3$-manifolds, it is now known that various 
flavors of Heegaard Floer homology, monopole Floer homology, and embedded 
contact homology are isomorphic, while their relationship with instanton Floer 
homology remains a major open question.

The relationships between Floer invariants of knots in $3$-manifolds are even 
less understood:
Knot instanton Floer homology is not known to be isomorphic to the other knot 
homologies, and while it is known that the usual knot monopole Floer homology 
is isomorphic to the \emph{hat} flavor of knot Heegaard Floer homology
(tensored with the mod-$2$ Novikov field $\ring$)
as graded modules over $\ring$ \cite{taubes2010ech5, colin2012hf, 
  kutluhan2012hf, lekili2013heegaard}:
\begin{equation}
  \label{eqn:khm-hfk}
  \KHM (Y, K; \ring) \cong \HFKh (Y, K; \F) \otimes \ring,
\end{equation}
no analogous statement is known for the more powerful \emph{minus} flavor of 
knot Heegaard Floer homology $\HFKm (Y, K; \F)$, which is a graded module over 
$\F [U]$ rather than $\F$.  In fact,
the \emph{minus} flavors $\KHMm$ and $\KHIm$ of knot monopole and instanton 
Floer homologies have been defined only recently by the second author 
\cite{li2019direct}
using contact handle attachment maps of sutured manifolds, based on work of 
Baldwin and Sivek \cite{baldwin2016contact} and inspired by work of Etnyre, 
Vela-Vick, and Zarev \cite{etnyre2017sutured}. As such, many basic structural 
properties of $\KHMm$ and $\KHIm$
are yet unknown.

For example, a key property of $\HFKm$ for knots $K \subset S^3$ is its unique 
$\F [U]$-summand, the negative of whose maximal Alexander $\Z$-grading is a 
concordance invariant $\tauH (K)$.%
\footnote{Technically, $\tauH$---usually simply denoted $\tau$---may depend on 
  the coefficient ring. In this article, we always take $\tauH (K)$ to mean 
  $\tauH (K; \F)$.}
In fact, $\tauH$ defines a homomorphism $\tauH \colon \concgrp \to \Z$ from the 
smooth concordance group $\concgrp$.  Moreover, $\abs{\tauH (K)}$ also gives a 
lower bound on the smooth $4$-genus $g_4 (K)$.
More generally, $\tauH$ can be defined for nullhomologous knots $K$ in a 
connected, oriented, closed $3$-manifold $Y$, with a choice of a Seifert 
surface $S$.
Inspired by this,
the second author \cite{li2019direct} similarly defines $\tauM (Y, K, S)$ and 
$\tauI (Y, K, S)$ to be the negative of the maximal Alexander $\Z$-grading of 
the non--$U$-torsion elements of $\KHMm$ and $\KHIm$.%
\footnote{When $Y = S^3$, we abbreviate these by $\tauM (K)$ and $\tauI (K)$.}
However, for knots $K \subset S^3$, these have not been shown to be concordance 
invariants or to give $4$-genus bounds.

In a different approach to the \emph{tau} invariants, Baldwin and Sivek 
\cite{baldwin2020framed} define a concordance invariant $\nushI$ using 
cobordism maps between the framed instanton Floer homology $\Ish$ of $S^3$ and 
of the integer surgeries $S^3_n (K)$ along $K$, and homogenize $\nushI$ to 
obtain a concordance invariant $\taushI$.%
\footnote{In \cite{baldwin2020framed}, $\taushI$ is simply denoted $\taush$; we 
  add the subscript $\mathrm{I}$ to separate it from the monopole version 
  $\taushM$.}
They show that $\abs{\taushI (K)} \leq g_4 (K)$, and that $2 \taushI$ gives a 
homomorphism $\taushI \colon \concgrp \to \R$ that is in fact a 
\emph{slice-torus invariant}, as defined by Lewark \cite{lewark2014slice} 
following Livingston \cite{livingston2004tau}; however, defined via a 
homogenization process, $\taushI$ is not known to be an integer (or even a 
rational number).  Nonetheless, these properties of $\taushI$ are sufficient 
for Baldwin and Sivek to use to determine $\Ish$ of all nonzero rational 
surgeries on $20$ of the $35$ nontrivial prime knots in $S^3$ through $8$ 
crossings, and establish several other results.  While it is not explicitly 
stated, a concordance invariant $\taushM$ can be similarly defined in the 
monopole Floer theory, via the \emph{tilde} flavor $\HMt (S^3_n (K); \ring)$.  
By construction, $\taushI$ and $\taushM$ are defined only for knots $K \subset 
S^3$.

This article represents the natural first step in understanding the structures 
of $\KHMm$ and $\KHIm$ and their comparisons with $\HFKm$. In the following, we 
shall replace the subscripts $\mathrm{M}$ and $\mathrm{I}$ (for ``monopole'' 
and ``instanton'') in $\tauM$ and $\tauI$ by the subscript $\mathrm{G}$ (for 
``gauge-theoretic'') in $\tauG$, when the statement applies to both theories.  
To begin, our main theorem identifies the \emph{tau} invariants, answering the 
question posed in (a previous version of) \cite{baldwin2020framed}:

\begin{theorem}
  \label{thm:main}
  For all knots $K \subset S^3$, we have $\tauG (K) = \taushG (K)$.
\end{theorem}

We immediately have the following corollaries in the instanton setting:

\begin{corollary}
  \label{cor:integers}
  For all knots $K \subset S^3$, the invariant $\taushI (K)$ is an integer. In 
  other words, $\taushI$ defines a homomorphism $\taushI \colon \concgrp \to 
  \Z$.  \qed
\end{corollary}

\begin{corollary}[cf.~{\cite[Proposition~5.4]{baldwin2020framed}}]
  \label{cor:g4}
  For all knots $K \subset S^3$, we have $\abs{\tauI (K)} \leq g_4 (K)$.  \qed
\end{corollary}

As mentioned above, Baldwin and Sivek \cite[Theorem~1.6]{baldwin2020framed} 
show that $2 \taushI$ is a slice-torus invariant, and use this to show that 
$\taushI (K)$ agrees with $g_4 (K)$ when $K$ is quasipositive.  Moreover, as 
Lewark \cite{lewark2014slice} proves that slice-torus invariants agree with the 
negative of the signature for alternating knots, they obtain $\taushI (K) = - 
\sigma (K) / 2$ for such knots.\footnote{We follow the convention where the 
  right-handed trefoil has signature $-2$.} Lewark also proves that the values 
of all slice-torus invariants agree on homogeneous knots, which gives $\taushI 
(K) = \tauH (K)$ for
such knots.

In the monopole setting, the statements in the preceding paragraph can be 
readily proved for $\taushM$ also.  Thus, \fullref{thm:main} immediately 
implies the following for knots in $S^3$:

\begin{corollary}[cf.~{\cite[Theorem~1.6]{baldwin2020framed}}]
  \label{cor:slice-torus}
  The invariant $2 \tauG$ is a slice-torus invariant. If $K$ is a quasipositive 
  knot, then $\tauG (K) = g_4 (K)$. If $K$ is an alternating knot, then $\tauG 
  (K) = - \sigma (K) / 2$. If $K$ is a homogeneous knot, then $\tauG (K) = 
  \tauH (K)$.\footnote{These facts combined show that $\tauI = \tauH$ for all 
    prime knots through $9$ crossings, except possibly $9_{42}$, $9_{44}$, and 
    $9_{48}$.} \qed
\end{corollary}

In fact, in the monopole setting, we can strengthen this last statement to hold 
for all knots:

\begin{theorem}
  \label{thm:tau-hf}
  For all knots $K \subset S^3$, we have $\tauM (K) = \tauH (K)$.
\end{theorem}

\begin{proof}
  Baldwin and Sivek \cite[Section~10]{baldwin2020framed} detail how the 
  Heegaard Floer $\tauH$ invariant can also be expressed as the homogenization 
  of a concordance invariant coming from surgeries, as explained to them by 
  Jennifer Hom. (One may reasonably denote such an invariant by $\taushH$.) 
  They then use this to show that if
  \[
    \dim_{\C} \Ish (Y; \C) = \dim_{\F} \HFh (Y; \F)
  \]
  holds for all $Y$ obtained via integer surgery along a knot in $S^3$, then 
  $\taushI (K) = \taushH (K) = \tauH (K)$ for all $K \subset S^3$ 
  \cite[Proposition~1.24]{baldwin2020framed}. The exact same proof can be 
  adapted to show that if
  \begin{equation}
    \label{eqn:hm-hf}
    \rk_{\ring} \HMt (Y; \ring) = \dim_{\F} \HFh (Y; \F)
  \end{equation}
  holds for all $Y$ obtained via integer surgery, then $\taushM (K) = \taushH 
  (K) = \tauH (K)$. But \eqref{eqn:hm-hf} is simply the isomorphism between 
  monopole and Heegaard Floer homologies for $3$-manifolds 
  \cite{taubes2010ech5, colin2012hf, kutluhan2012hf}. Thus, our claim follows 
  from \fullref{thm:main}.
\end{proof}

The significance of \fullref{thm:tau-hf} is that it represents the first step 
towards proving the generalization of the isomorphism between $\KHM$ and 
$\HFKh$ in \eqref{eqn:khm-hfk} to the \emph{minus} flavor:

\begin{conjecture}
  \label{conj:khm-hfk-minus}
  Let $Y$ be a connected, oriented, closed $3$-manifold, and let $K \subset Y$ 
  be an oriented, nullhomologous knot. Then there is an isomorphism of graded 
  modules over $\ring [U]$:
  \[
    \KHMm (Y, K; \ring) \cong \HFKm (Y, K; \F) \otimes \ring [U].
  \]
\end{conjecture}

\fullref{cor:slice-torus} has another implication, as pointed out to the 
authors by Steven Sivek:

\begin{corollary}
  \label{cor:tb-r}
  Suppose that $\leg \subset (S^3, \xistd)$ is a Legendrian knot of smooth knot 
  type $K$; then
  \[
    \tb (\leg) + \abs{\rot (\leg)} \leq 2 \tauG (K) - 1.
  \]
\end{corollary}

\begin{proof}
  Consider the positive and negative transverse pushoffs $\trans_\pm (\leg)$, 
  which have self-linking numbers $\slink (\trans_\pm (\leg)) = \tb (\leg) \mp 
  \rot (\leg)$ respectively. By \cite[Theorem~6.1]{baldwin2020framed}, $\slink 
  (\trans) \leq 2 \taushG (K) - 1$ for all transverse representatives $\trans$ 
  of $K$.  (\cite[Theorem~6.1]{baldwin2020framed} is a statement for $\taushI$, 
  but the same argument works for $\taushM$.) Thus, the result follows from 
  \fullref{thm:main}.
  (Note that \cite[Theorem~6.1]{baldwin2020framed} is in fact the key 
  ingredient in proving that $\taushI (K) = g_4 (K)$ for quasipositive knots 
  $K$.)
\end{proof}

\begin{remark}
  \label{rmk:plamenevskaya}
  The analogous statement that
  \begin{equation}
    \label{eq:plamenevskaya}
    \tb (\leg) + \abs{\rot (\leg)} \leq 2 \tauH (K) - 1,
  \end{equation}
  first proved by Plamenevskaya \cite{Olga2004tb}, implies \fullref{cor:tb-r} 
  for $\tauM$ via \fullref{thm:tau-hf}. Alternatively, one could also prove it 
  using \cite[Lemma~3.10]{li2019decomposition}.
\end{remark}

Below, we describe the strategy to prove \fullref{thm:main}.
To simplify our notation, we first set up some conventions for the rest of the 
article.

\subsection*{Conventions}
  The coefficient ring for monopole Floer homologies is always taken to be the 
  mod-$2$ Novikov field, and that for instanton Floer homologies is always 
  taken to be the field $\C$ of complex numbers.  In both cases, we shall 
  denote the coefficient ring by $\ring$.  Similar to $\tauG$, we shall denote 
  both $\shm$ and $\shi$ by $\shg$ when a statement applies to both sutured 
  monopole and sutured instanton Floer homologies, and likewise denote by 
  $\KHG$ (resp.\ $\KHGm$) the knot monopole and instanton Floer homologies 
  $\KHM$ (resp.\ $\KHMm$) and $\KHI$ (resp.\ $\KHIm$).

\subsection{Strategy}
\label{ssec:strategy}

The astute reader may have noticed that we did not state the concordance 
invariance of $\tauG$, or its additivity under connected sum, as a corollary of 
\fullref{thm:main}.  The reason is that, in order to prove \fullref{thm:main}, 
we shall in fact \emph{first} prove the concordance invariance of $\tauG$:

\begin{proposition}
  \label{prop:tau-conc}
  For all knots $K \subset S^3$, the integer $\tauG (K)$ is a concordance 
  invariant.
\end{proposition}

To establish \fullref{prop:tau-conc}, we shall also prove the key property that 
$\KHGm$ has a unique $\ring [U]$-summand (also known as an infinite $U$-tower) 
for knots $K \subset S^3$, analogous to $\HFKm$:

\begin{proposition}
  \label{prop:KHGm-tower}
  For all knots $K \subset S^3$, $\KHGm (S^3, K)$ has a unique $\ring 
  [U]$-summand.
\end{proposition}

After establishing \fullref{prop:tau-conc}, we shall turn to the additivity of 
$\tauG$ under connected sum:

\begin{proposition}
  \label{prop:tau-conn-sum}
  For all pairs of knots $K_1, K_2 \subset S^3$, we have $\tauG(K_1 \connsum 
  K_2) = \tauG(K_1) + \tauG(K_2)$.
\end{proposition}

The rest of the proof of \fullref{thm:main} can be described roughly as 
follows. Recall that $\KHGm$ is defined in terms of a directed system of $\SHG$ 
of the knot complement $S^3 (K)$ with sutures $\Gamma_n$, over different values 
of $n$, where $\Gamma_n$ denotes a pair of parallel sutures on the boundary 
torus with $n$ full twists.  First, using bypass and surgery exact triangles 
involving $-S^3 (K)$, we reformulate $\tauG$ in terms of the twisting 
coefficient $n_0$ for which $\SHG$ of $-S^3 (K)$ with $-\Gamma_{n_0}$ sutures 
uniquely attains minimum rank. (This is conceptually similar to Baldwin and 
Sivek's notion of \emph{V-shaped} knots.) Next, noting that whether the 
inequality $\rk_{\ring} \SHG (-S^3 (K), -\Gamma_{n+1}) > \rk_{\ring} \SHG (-S^3 
(K), -\Gamma_n)$ holds is equivalent to the (non-)vanishing of certain surgery 
cobordism maps involving $-S^3 (K)$, further analysis using surgery exact 
triangles allows us to relate $\tauG$ to the (non-)vanishing of surgery 
cobordism maps involving $\Ish$ or $\HMt$ of $-S^3_{-n} (K)$, and thence to 
$\nushG$, giving the inequality
\[
  2 \tauG (K) - 1 \leq \nushG (K) \leq 2 \tauG (K) + 1
\]
whenever $\nushG (K) \neq 0$. A homogenization argument, using the fact that 
$\tauG$ is a concordance homomorphism, completes the proof.

\subsection{Examples}
\label{ssec:examples}

\begin{figure}[htbp]
  \captionsetup{aboveskip={\dimexpr10pt+2pt+\sactualfontsize\relax}}
  \labellist
  \small\hair 2pt
  \pinlabel {$K_m$} [t] at 64 0
  \pinlabel {$\mir{K}_m$} [t] at 289 0
  \pinlabel {\red{$\alpha$}} [bl] at 140 235
  \pinlabel {$m$ negative full twists $\left\{ \rule{0in}{0.6in} \right.$} [r] 
  at -5 170
  \pinlabel {$\left. \rule{0in}{0.6in} \right\} \, m$ positive full twists} [l] 
  at 353 170
  \endlabellist
  \includegraphics[width=0.35\textwidth]{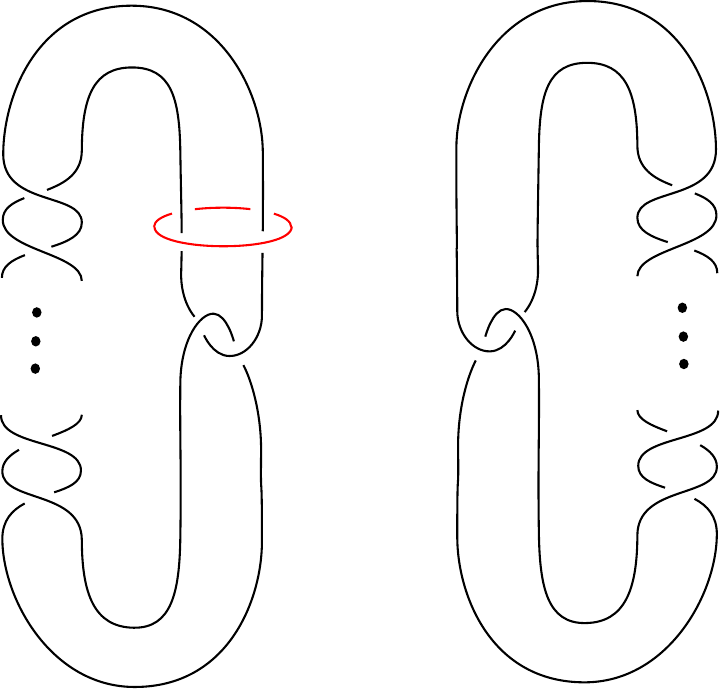}
  \caption{The twist knots $K_m$ and $\mir{K}_m$. In \fullref{sec:twist-knots}, 
    we shall perform surgery along the curve $\alpha$ in the proof of 
    \fullref{thm:twist-knots}.}
  \label{fig:twist-knots}
\end{figure}

Let $K_m \subset S^3$ be the twist knot with a positive clasp and $m$ negative 
full twists (or $-m$ positive full twists if $m < 0$), and let $\mir{K}_m$ 
denote its mirror image; see \fullref{fig:twist-knots}. (In the notation of 
Baldwin and Sivek \cite{baldwin2020framed}, their $K_n$ corresponds to our 
$K_{-n/2}$ when $n$ is even, and to our $\mir{K}_{(n+1)/2}$ when $n$ is odd.) 
Baldwin and Sivek \cite{baldwin2020framed} compute\footnote{They only compute 
  half of these, but the antisymmetry of $\nushI$ under mirroring gives the 
  other half.}
\begin{equation}
  \label{eq:nushI-bs}
  \nushI (K_m) =
  \begin{cases}
    0 & \text{for } m \leq 0,\\
    1 & \text{for } m > 0,
  \end{cases}
  \qquad
  \nushI (\mir{K}_m) =
  \begin{cases}
    0 & \text{for } m \leq 0,\\
    -1 & \text{for } m > 0,
  \end{cases}
\end{equation}
and use it to fully determine $\dim_{\C} \Ish (S^3_{p/q} (K_m))$.
One can also compute $\taushI (K_m)$ (and hence $\taushI (\mir{K}_m)$) as 
follows: Since twist knots are alternating, 
\cite[Corollary~1.10]{baldwin2020framed} says that $\taushI (K_m) = - \sigma 
(K_m) / 2$, and the signature $\sigma (K_m)$ can be directly computed from the 
$2 \times 2$ Seifert matrix. This gives
\[
  \taushI (K_m) =
  \begin{cases}
    0 & \text{for } m \leq 0,\\
    1 & \text{for } m > 0,
  \end{cases}
  \qquad
  \taushI (\mir{K}_m) =
  \begin{cases}
    0 & \text{for } m \leq 0,\\
    -1 & \text{for } m > 0.
  \end{cases}
\]
(One can also use \eqref{eq:nushI-bs} to compute $\taushI (K_m)$ without 
computing $\sigma (K_m)$, using \cite[Theorem~3.7, Proposition~5.4, and 
Corollary~1.10]{baldwin2020framed}.)
With this in hand, to illustrate \fullref{thm:main}, we provide a direct and 
complete computation of $\KHGm (-S^3, K_m)$ and $\tauG$ for this infinite 
family.

\begin{theorem}
  \label{thm:twist-knots}
  We abbreviate by $\ring$ the $\ring [U]$-module $\ring [U] / U$, and denote 
  Alexander gradings by subscripts, and direct sums by superscripts.
  \begin{enumerate}
    \item For $m \leq 0$, we have
      \begin{gather*}
        \KHGm (-S^3, K_m) \isom \KHGm (-S^3, \mir{K}_m) \isom \ring [U]_0 
        \dirsum \ring_1^{-m} \dirsum \ring_0^{-m},\\
        \tauG (K_m) = \tauG (\mir{K}_m) = 0.
      \end{gather*}
    \item For $m > 0$, we have
      \begin{align*}
        \KHGm (-S^3, K_m) &\isom \ring[U]_1 \dirsum \ring_1^{m-1} \dirsum 
        \ring_0^{m}, &
        \KHGm (-S^3, \mir{K}_m) &\isom \ring [U]_{-1} \dirsum \ring_1^{m} 
        \dirsum \ring_0^{m-1},\\
        \tauG (K_m) &= 1, & \tauG (\mir{K}_m) &= -1.
      \end{align*}
  \end{enumerate}
\end{theorem}

\subsection{Future work}
\label{ssec:future-work}

As this article represents the first step in our major goal to understand the 
structures of $\KHGm$, we present here some open questions that arise naturally 
from our discussion.

First, while \fullref{cor:slice-torus} gives $\tauG (K)$ for all alternating 
knots $K \subset S^3$, one may reasonably hope to fully determine $\KHGm (S^3, 
K)$ for such $K$. Indeed, we do so in \fullref{thm:twist-knots} for twist 
knots, which are alternating.  In knot Heegaard Floer homology, if $K$ is 
alternating with signature $\sigma = \sigma (K)$, and symmetrized Alexander 
polynomial $\Delta_K (t) = \sum_i a_i \cdot t^i$, then
\[
  \HFKm (S^3, K; \F) \isom \F [U]_{\frac{\sigma}{2}} \dirsum 
  \paren{\bigdirsum_{i \leq \frac{\sigma}{2}} (\F [U] / U)_i^{\abs{b_i}-1}} 
  \dirsum
  \paren{\bigdirsum_{i < \frac{\sigma}{2}} (\F [U] / U)_i^{\abs{b_i}}},
\]
where $b_i = \sum_{j \geq 0} a_{i+j}$; see \cite{ozsvath2003alternating}, and 
also \cite[Corollary~10.3.2]{ozsvath2015gridbook}.

\begin{question}
  \label{qn:alternating}
  Is there an analogous formula for $\KHGm (S^3, K)$ for alternating knots $K$, 
  or at least two-bridge knots $K$?
\end{question}

Another question is the mirroring of knots. It follows from our work that if 
$\mir{K}$ is the mirror of $K$, then $\tauG (\mir{K}) = -\tauG (K)$. In knot 
Heegaard Floer homology, one has the following more precise formula: If
\[
  \HFKm (S^3, K; \F) \isom \F [U]_{-\tau} \dirsum \paren{\bigdirsum_{i=1}^k (\F 
    [U] / U)_{s_i}^{n_i}},
\]
then
\[
  \HFKm (S^3, \mir{K}; \F) \isom \F [U]_{\tau} \dirsum 
  \paren{\bigdirsum_{i=1}^k (\F [U] / U)_{n_i - s_i}^{n_i}},
\]
where $\tau = \tauH (K)$; see \cite[Section~3.5]{ozsvath2004holomorphicknot}, 
and also \cite[Proposition~7.4.3]{ozsvath2015gridbook}.

\begin{question}
  \label{qn:mirror}
  Is there an analogous formula for $\KHGm (S^3, \mir{K})$ in terms of $\KHGm 
  (S^3, K)$?
\end{question}

Aside from symmetry, there are questions concerning the behavior of $\KHGm$ 
under crossing changes and with respect to skein relations. In particular, if 
$K_-$ is the result of changing a positive crossing in $K_+$ to a negative 
crossing, then there exists graded $\F [U]$-module maps
\[
  C_- \colon \HFKm (S^3, K_+) \to \HFKm (S^3, K_-), \qquad C_+ \colon \HFKm 
  (S^3, K_-) \to \HFKm (S^3, K_+),
\]
such that $C_- \circ C_+$ and $C_+ \circ C_-$ is each equal to multiplication 
by $U$.  Exploiting this, Alishahi and Eftekhary \cite{alishahi2020torsion} 
define a $U$-torsion order invariant $\mathfrak{l} (K)$ that gives a lower 
bound on the unknotting number $u (K)$. A generalized version is used by 
Juh\'asz, Miller, and Zemke \cite{juhasz2020torsion} to obtain an obstruction 
to connected knot cobordisms with a given number of local maxima.

\begin{question}
  \label{qn:torsion-order}
  Are there analogous maps for $\KHGm (S^3, K_+)$ and $\KHGm (S^3, K_-)$, and 
  consequently a $U$-torsion order invariant in knot instanton and monopole 
  Floer theory?
\end{question}

Note that a positive answer to \fullref{qn:torsion-order} would imply that $0 
\leq \tauG (K_+) - \tauG (K_-) \leq 1$, a fact that can be deduced from 
\fullref{cor:slice-torus}; see \cite[Corollary~3]{livingston2004tau}.

In regard to the oriented skein relation, Kronheimer and Mrowka 
\cite[Theorem~3.1]{kronheimer2010instanton} prove that $\KHI$ satisfies an 
exact triangle relating $K_+$, $K_-$, and their oriented resolution $K_0$. The 
analogous relation in knot Heegaard Floer homology is satisfied by both $\HFKh$ 
and $\HFKm$ \cite{ozsvath2007skein}, and like its instanton counterpart, has 
proved to be a very useful tool. One may thus ask:

\begin{question}
  \label{qn:oriented-skein}
  Does $\KHIm$ satisfy an oriented skein relation?
\end{question}

In order to answer \fullref{qn:oriented-skein}, one must necessarily generalize 
the definition of $\KHIm$ to links with multiple components. This has been 
carried out by the first and second authors 
\cite[Section~6.2]{li2019decomposition}, who define $\KHGm (Y, L)$ for 
nullhomologous links $L \subset Y$.

\begin{question}
  \label{qn:links}
  Can $\tauG$ be generalized to a multiset of values for links $L \subset Y$, 
  and if so, what properties does it satisfy?
\end{question}

The analogous notion in knot Heegaard Floer theory is that of the 
\emph{$\tauH$-set} of a link \cite{ozsvath2015gridbook, cavallo2018tau, 
  hedden2020tau}.  Notably, Hedden and Raoux \cite[Theorem~2]{hedden2020tau} 
prove that the $\tauH$-set of a link $L \subset Y$ satisfy many interesting 
properties previously known for $\tauH (S^3, K)$, including concordance 
invariance in $Y$, crossing-change inequalities in $Y$, $4$-genus bounds in 
$Y$, and, in the case $L \subset S^3 = \bdy W$ where $W$ is a definite 
$4$-manifold, an inequality for surfaces in $W$ bound by $L$.

Finally, we turn to Legendrian knot invariants. For a Legendrian knot $\leg 
\subset (Y, \xi)$ of smooth knot type $K$, Baldwin and Sivek 
\cite{baldwin2014invariants, baldwin2016equivalence} define a class $\losshM 
(\leg) \in \KHM (-Y, K)$, and show it to be equivalent to the \emph{LOSS 
  invariant} $\losshH (\leg) \in \HFKh (-Y, K)$.\footnote{For $\leg \subset 
  (S^3, \xistd)$, the Alexander grading of $\losshH (\pm \leg)$ is $(\tb (\leg) 
  \mp \rot (\leg) + 1) / 2$, which offers another proof of 
  \eqref{eq:plamenevskaya}.}  Notably, their work implies that $\losshH$ gives 
an obstruction to the existence of exact Lagrangian cobordisms between 
Legendrian knots---without adjectives such as \emph{decomposable} or 
\emph{regular}---for which there is currently no proof purely in Heegaard Floer 
theory. On the other hand, $\losshH$ has a generalization $\lossH (\leg) \in 
\HFKm (-Y, K)$, which is a non--$U$-torsion class that is mapped to $\losshH 
(\leg)$ under the natural map $\HFKm (-Y, K) \to \HFKh (-Y, K)$. Etnyre, 
Vela-Vick, and Zarev \cite{etnyre2017sutured} place $\lossH$ in the context of 
$\HFKm$ as the limit of a directed system of $\SFH$; following this strategy, 
we may also define a Legendrian invariant $\lossM (\leg) \in \KHMm (-Y, K)$, 
which is mapped to $\losshM (\leg)$ under the natural map $\KHMm (-Y, K) \to 
\KHM (-Y, K)$. By the naturality in monopole Floer theory 
\cite{baldwin2015naturality}, $\lossM$ is a well-defined class---and not only a 
class defined up to isomorphism---in $\KHMm (-Y, K)$.

\begin{question}
  \label{qn:loss}
  Is the Legendrian invariant $\lossM (\leg) \in \KHMm (-Y, K)$ effective in 
  distinguishing Legendrian knots, or in obstructing exact Lagrangian 
  cobordisms?
\end{question}

\subsection{Organization}
\label{ssec:organization}

We review the definitions of $\KHGm$, $\tauG$, and $\taushG$ in 
\fullref{sec:prelim}. In \fullref{sec:tau-conc}, we prove 
\fullref{prop:KHGm-tower} and \fullref{prop:tau-conc}, establishing that 
$\tauG$ is a concordance invariant; in \fullref{sec:tau-conn-sum}, we prove 
\fullref{prop:tau-conn-sum}, the additivity of $\tauG$. We then carry out the 
argument described in \fullref{ssec:strategy} to prove \fullref{thm:main}, 
identifying the \emph{tau} invariants in \fullref{sec:identify}. Finally, we 
compute $\KHGm$ and $\tauG$ for twist knots in \fullref{sec:twist-knots}, 
proving \fullref{thm:twist-knots}.

\subsection*{Acknowledgements}

Many of the ideas in this article were developed while ZL was visiting SG and 
CMMW at Louisiana State University, and the authors thank LSU for their 
hospitality. Part of the research was conducted while ZL was at the 
Massachusetts Institute of Technology. The authors are grateful to Steven Sivek 
for pointing out \fullref{cor:tb-r}. The authors are deeply indebted to John 
Baldwin, Tom Mrowka, Ina Petkova, Steven Sivek, and Shea Vela-Vick for the many 
helpful discussions, and more importantly, for their unwavering support.

SG was partially supported by Shea Vela-Vick's NSF Grant DMS-1907654. While at 
MIT, ZL was partially supported by Tom Mrowka's NSF Grant 1808794. CMMW was 
partially supported by NSF Grant DMS-2010863 and an AMS--Simons Travel Grant.

%% file: sec-prelim.tex
\section{Preliminaries}
\label{sec:prelim}

\subsection{\texorpdfstring{$\KHG$}{Knot instanton and monopole Floer 
    homologies} and naturality}
\label{ssec:khg-naturality}

In this article, we shall focus on oriented, based knots $(K, p) \subset S^3$ 
and $(K, p) \subset -S^3$.  As in \cite[Section~8]{baldwin2015naturality}, by 
the knot complement $S^3 (K)$ and the meridional sutures $\Gamma_{\mu}$, we 
mean the following: Let $D^2$ be the unit disk in the complex plane with 
boundary $S^1 = \bdy D^2$, and let $\vphi \colon S^1 \cross D^2 \to S^3$ be an 
embedding such that $\vphi (S^1 \cross \set{0}) = K$ and $\vphi (\set{1} \cross 
\set{0}) = p$; then
\[
  (S^3 (K), \Gamma_\mu) = (S^3 \setminus \Int (\Im (\vphi)), \mu_\vphi^+ \union 
  - \mu_\vphi^-),
\]
where $\mu_{\vphi}^\pm$ is the oriented meridian $\vphi (\set{\pm1} \cross \bdy 
D^2)$ on $\bdy S^3 (K)$.
Of course, this definition does not quite make sense yet, as it depends on the 
choice of $\vphi$.
In work of Kronheimer and Mrowka \cite{kronheimer2010knots}, the balanced 
sutured manifold $(S^3 (K), \Gamma_{\mu})$ is used to construct the knot 
instanton and monopole Floer homologies:
\[
  \KHG (S^3, K, p) = \SHG (S^3 (K), \Gamma_\mu).\footnote{The basepoint $p$ is 
    omitted in the Kronheimer--Mrowka definition.}
\]

The sutured instanton and monopole Floer homologies $\SHG (M, \gamma)$, defined 
in general for balanced sutured manifolds $(M, \gamma)$, themselves depend on 
the choice of a \emph{closure}, a closed $3$-manifold $Y$ obtained by gluing an 
auxiliary piece to $(M, \gamma)$ and then identifying the remaining boundary 
components, together with a distinguished surface $R \subset Y$.  Kronheimer 
and Mrowka assign modules to each such closure, and show that these modules are 
all isomorphic.  By refining the notion of closures, Baldwin and Sivek 
\cite{baldwin2015naturality} prove that there are in fact \emph{canonical 
  isomorphisms} relating these modules---well defined up to multiplication by a 
unit in $\ring$---and use them to build a \emph{projectively transitive system} 
$\bSHG$ for balanced sutured manifolds.  By abuse of notation, whenever we 
write $\SHG$ in the sequel, we shall mean the canonical module associated to 
$\bSHG$.\footnote{Technically, the canonical module is only defined for an 
  honestly (i.e.\ not projectively) transitive system; for projectively 
  transitive systems, one would only obtain a module modulo multiplication by a 
  unit in $\ring$, which is
  only a set. Instead, we choose to interpret $\SHG$ as an actual 
  $\ring$-module, whose \emph{elements} are well defined only up to 
  multiplication by a unit in $\ring$.}

Coming back to $\KHG (S^3, K, p)$, while the definition of $(S^3 (K), 
\Gamma_\mu)$ above depends on $\vphi$, Baldwin and Sivek 
\cite[Proposition~8.2]{baldwin2015naturality} further prove that there are 
canonical isomorphisms relating $\bSHG$ of the sutured manifolds $(S^3 (K), 
\Gamma_\mu)$ constructed using different embeddings $\vphi$ and $\vphi'$.  This 
proof hinges on the fact that the basepoint $p$ is fixed, and explains the 
notation $\KHG (S^3, K, p)$. Once again, this leads to a projectively 
transitive system $\bKHG (S^3, K, p)$, and we shall take $\KHG (S^3, K, p)$ to 
mean the associated canonical module.

\subsection{\texorpdfstring{$\KHGm$}{The \emph{minus} flavor} and 
  \texorpdfstring{$\tauG$}{\emph{tau}}}
\label{ssec:khgm-tau}

In this subsection, we recall the construction of $\KHGm$ and $\tauG$ by the 
second author \cite{li2019direct}.

Let $S \subset S^3$ be an oriented, minimal-genus Seifert surface of $K$.  The 
surface $S$ induces a framing on the boundary of the knot complement $S^3(K)$ 
and hence longitude $\lambda$ (whose orientation agrees with that of $K$).  Let 
$\vphi \colon S^1 \cross D^2 \to S^3$ be as before, with $\vphi (S^1 \cross 
\set{0}) = K$ and $\vphi (\set{1} \cross \set{0}) = p$; then define the 
balanced sutured manifold
\[
  (S^3 (K), \Gamma_n) = (S^3 \setminus \Int (\Im (\vphi)), \lambda_{\vphi,n}^+ 
  \union - \lambda_{\vphi,n}^-),
\]
where $\lambda_{\vphi,n}^+$ is the oriented longitude $\vphi (\set{e^{it} 
  \cross e^{i(-nt)}}_{t \in [0, 2\pi)})$, and $\lambda_{\vphi,n}^-$ is the 
oriented longitude $\vphi (\set{e^{it} \cross e^{i(-nt+\pi)}}_{t \in [0, 
  2\pi)})$, on $\bdy S^3 (K)$. Note that
$\Ga_n \subset \bdy S^3 (K)$ is the union of two disjoint, parallel, oppositely 
oriented simple closed curves of slope $-n$ (or, equivalently, of class 
$\pm([\lambda]-n[\mu])\in H_1(\partial{S^3(K)})$). Like $(S^3 (K), 
\Gamma_\mu)$, the sutured manifold $(S^3 (K), \Ga_n)$ depends on the choice of 
$\vphi$. By an argument similar to that of 
\cite[Proposition~8.2]{baldwin2015naturality}, there are canonical isomorphisms 
relating $\bSHG$ of $(S^3 (K), \Gamma_n)$ constructed using different 
embeddings $\vphi$ and $\vphi'$.

Note that, while our exposition so far focuses on $(S^3 (K), \Gamma_\mu)$ and 
$(S^3 (K), \Gamma_n)$, a similar construction gives the balanced sutured 
manifolds $(-S^3 (K), -\Gamma_\mu)$ and $(-S^3 (K), -\Gamma_n)$, which we shall 
use extensively.

To define $\KHGm$, maps $\psi^n_{\pm,n+1}$ are defined in 
\cite{li2019direct},\footnote{Technically, these maps are well defined only up 
  to multiplication by a unit in $\ring$, for the same reason as before. Here 
  and in the rest of the article, we say that $f = g$ if the maps $f$ and $g$ 
  on $\SHG$ (or consequently $\KHGm$) agree up to multiplication by a unit in 
  $\ring$. In particular, \eqref{eq:comm-diag} commutes up to multiplication by 
  a unit in $\ring$.} which fit into a commutative diagram
\begin{equation}
  \label{eq:comm-diag}
  \vcenter{
    \xymatrix{
      \dotsb\ar[r]&\shg(-S^3(K),-\Ga_{n})\ar[r]^{\,\,\psi^n_{-,n+1}\quad}\ar[d]^{\psi^n_{+,n+1}}&\shg(-S^3(K),-\Ga_{n+1})\ar[r]\ar[d]^{\psi^{n+1}_{+,n+2}}&\dotsb\phantom{.}\\
      \dotsb\ar[r]&\shg(-S^3(K),-\Ga_{n+1})\ar[r]^{\psi^{n+1}_{-,n+2}}&\shg(-S^3(K),-\Ga_{n+2})\ar[r]&\dotsb.
    }
  }
\end{equation}
Each horizontal row forms a \emph{directed system} of 
$\ring$-modules.\footnote{As before, this is really a directed system of 
  ``$\ring$-modules whose elements are well defined only up to multiplication 
  by a unit in $\ring$''. One could take the alternative viewpoint of choosing 
  an honest $\ring$-module representative for each $\SHG (-S^3 (K), -\Gamma_n)$ 
  by specifying an embedding $\vphi$ and a closure $(Y, R)$; however, more work 
  would be necessary to take care of the fact that $\psi^{n}_{\pm,n+1}$ is only 
  well defined (up to multiplication by a unit) for ``compatible'' closures of 
  $(-S^3 (K), \Gamma_n)$ and $(-S^3 (K), \Gamma_{n+1})$, which necessarily have 
  auxiliary pieces of the same genus.}
(Note that we choose to work primarily with $(-S^3 (K), -\Gamma_n)$ instead of 
$(S^3 (K), \Gamma_n)$, because the definition of $\psi^n_{\pm,n+1}$ makes use 
of a \emph{contact element} $\phi_\xi \in \SHG (-M, -\gamma)$ 
\cite{baldwin2016contact, baldwin2016instanton}, defined for a contact 
structure $\xi$ on $(M, \gamma)$.)

\begin{definition}[{\cite[Definition~5.4]{li2019direct}}]
  \label{defn: minus version}
  The \emph{minus knot monopole or instanton Floer homology} $\khg(-S^3,K,p)$ 
  is the direct limit of the directed system in \eqref{eq:comm-diag}, which is 
  an $\ring$-module whose elements are well defined up to multiplication by a 
  unit in $\ring$.
  The collection of maps $\set{\psi^n_{+,n+1}}_{n\in\Z_+}$ defines a map on the 
  direct limit
  \[
    U \colon \khg(-S^3,K,p)\ra \khg(-S^3,K,p),
  \]
  which gives $\KHGm (-S^3, K, p)$ an $\ring [U]$-module structure.
\end{definition}

In the following paragraphs, we describe the grading on $\KHGm (-S^3, K, p)$.   
Our description shall be brief; for more details, see \cite[Section~3 and 
4]{li2019direct}.

Fix the balanced sutured manifold $(- S^3 (K), - \Gamma)$, where $\Gamma$ is 
either $\Gamma_\mu$ or $\Gamma_n$ for some $n$. Now realize $S$ as a properly 
embedded surface $(\swb{S}, \bdy \swb{S}) \subset (-S^3 (K), -\bdy (S^3 (K))$; 
then $\bdy \swb{S} \intersect \Gamma$ must consist of exactly $2k$ points for 
some $k$.  The realization of $S$ as $\swb{S}$ in fact involves a choice, 
corresponding to the value of $k$.  By isotoping $\swb{S}$ near its boundary, 
one could create a new pair of intersection points with $\Gamma$; this is 
called the \emph{positive} or \emph{negative stabilization} of $\swb{S}$ 
depending on the isotopy. We denote by $\swb{S}^{q}$ (resp.\ $\swb{S}^{-q}$) 
the result of performing $q$ positive (resp.\ negative) stabilizations on 
$\swb{S}$. (When $q = 1$, we also denote these by $\swb{S}^\pm$.) It is proved 
\cite[Theorem~3.4]{li2019direct} that $(\swb{S}, \bdy \swb{S}) \subset (- S^3 
(K), -\bdy S^3 (K))$ induces a $\Z$-grading on $\SHG (-S^3 (K), -\Gamma)$ 
whenever $\bdy \swb{S}$ intersects $\Gamma$ at $2k$ points, where $k$ is odd, 
and a formula \cite[Proposition~4.3]{li2019direct} is given that relates the 
$\Z$-gradings associated to Seifert surfaces related by stabilizations: For all 
$r \in \Z$, we have
\begin{equation}
  \label{eq:gr-shift}
  \SHG (-S^3 (K), -\Gamma, \swb{S}^{q+2r}, i) \isom \SHG (-S^3 (K), -\Gamma, 
  \swb{S}^q, i+r),
\end{equation}
where $\SHG (-S^3 (K), -\Gamma, \swb{S}, i)$ denotes the summand in grading $i 
\in \Z$.

Now fix $(- S^3 (K), -\Gamma_n)$ for some $n$. Since the longitude $\lambda$ is 
the boundary of the Seifert surface $S$, and $\Gamma_n$ is of class $\pm 
([\lambda] - n [\mu])$, it follows that $S$ has a realization $(S_n, \bdy S_n) 
\subset (-S^3 (K), -\bdy S^3 (K))$ such that $\bdy S_n \intersect \Gamma_n$ 
consists of exactly $2n$ points. Then, for $n$ odd (resp.\ even), we obtain 
$\Z$-gradings induced by the surfaces $S_n$ (resp.\ $S_n^-$); for brevity, we 
write $S_n^{\tau (n)}$ for $S_n$ when $n$ is odd, and $S_n^-$ when $n$ is even; 
i.e.\ $\tau (n) = 0$ or $-1$.  For $(- S^3 (K), -\Gamma_\mu)$, we obtain a 
$\Z$-grading induced by the surface $S_\mu$ that intersects $\Gamma_\mu$ at 
exactly $2$ points.

It is then proved \cite[Propositions~5.5 and 5.6]{li2019direct} that, after an 
appropriate grading shift
\[
  \sigma (n) = \frac{n - 1 - \tau (n)}{2},\footnote{The difference between this 
    and the formula for $\sigma (n)$ on \cite[p.~42]{li2019direct} is due to a 
    typographical error in the cited article.}
\]
the maps $\psi^n_{-,n+1}$ in the directed system in \eqref{eq:comm-diag} become 
grading-preserving maps, i.e.\ they each decompose into maps
\[
  \psi^n_{-,n+1} \colon \SHG (-S^3 (K), -\Gamma_n, S_n^{\tau (n)}, i) [\sigma 
  (n)] \to \SHG (-S^3 (K), -\Gamma_{n+1}, S_{n+1}^{\tau (n+1)}, i) [\sigma 
  (n+1)]
\]
for $i \in \Z$. Thus, the Seifert surface $S$ induces a $\Z$-grading on $\KHGm 
(-S^3, K, p)$, known as the \emph{Alexander grading}. The maps 
$\psi^n_{+,n+1}$, and hence the action of $U$ on $\KHGm$, is then of degree 
$-1$. For knots inside $S^3$ or $-S^3$, the Alexander grading is independent of 
the choice of the Seifert surface $S$; we shall therefore suppress $S$ from the 
notation.  Thus, we obtain a decomposition
\[
  \KHGm (-S^3, K, p) = \bigdirsum_{i \in \Z} \KHGm (-S^3, K, p, i),
\]
where $\KHGm (-S^3, K, p, i)$ denotes the summand in Alexander grading $i \in 
\Z$.

Inspired by the \emph{tau} invariant in knot Heegaard Floer homology defined by 
Ozsv\'ath and Szab\'o \cite{ozsvath2003tau}, we have the following definition:

\begin{definition}[{\cite[Definition~5.7]{li2019direct}}]
  \label{def:tau}
  For a knot $K \subset S^3$, the \emph{instanton or monopole tau invariant} is 
  defined as
  \[
    \tauG (K) = \max \setc{i \in \Z}{\text{there is a homogeneous, 
        non--}U\text{-torsion element } x \in \KHGm (-S^3, K, p, i)}.
  \]
  (Here, a non--$U$-torsion element $x$ is one such that $U^j x \neq 0$ for all 
  $j \geq 0$.)
\end{definition}

In the sequel, we shall often compute the rank of $\khg (-S^3, K, p, i)$ in a 
specific Alexander grading $i \in \Z$ as an $\ring$-module. We claim that this 
completely determines the $\ring$-module isomorphism type of $\KHGm (-S^3, K, 
p, i)$: For $\khi$, this is clear, since the module is a vector space over 
$\mathbb{C}$.
For $\khm$, our claim is a consequence of the following lemma.

\begin{proposition}
  \label{lem: only rank is important}
  For any based knot $(K,p)$ in $S^3$ and any $i\in\intg$, the $\ring$-module 
  $\khm(-S^3,K,p,i)$ is free and of finite rank.
\end{proposition}

\begin{proof}
  From \cite[Proposition~5.10]{li2019direct}, we know that there exists a 
  sufficiently large $n\in \intg$, such that
  \[
    \khm(-S^3,K,p,i)\cong\shm(-S^3(K),-\Ga_n,S_n^{\tau(n)},j)
  \]
  for some $j \in \Z$, which shows the $\ring$-module is of finite rank.
  (Here, $i$ might not be equal to $j$ because of the grading shift in the 
  definition of $\khm$ that we mentioned earlier.)

  To prove that it is free, recall that by work of Kronheimer and Mrowka 
  \cite[Lemma~4.9]{kronheimer2010knots}, for a balanced sutured manifold 
  $(M,\ga)$ and a coefficient ring $\ring$ of characteristic $0$, we have an 
  isomorphism of $\ring$-modules
  \begin{equation}
    \label{eq:local-coeff}
    \SHM (M, \gamma; \ring) \cong
    \SHM(M,\ga;\Gamma_{\eta})\cong \SHM(M,\ga;\Z)\otimes_{\Z}\mathcal{R},
  \end{equation}
  which respects the grading. (Here, $\Gamma_{\eta}$ denotes a local system 
  whose fiber at every point is $\ring$, and is unrelated to the sutures 
  $\Gamma_\mu$ and $\Gamma_n$.) Sivek \cite[Section~2.2]{sivek2012monopole} 
  extends $\SHM$ to mod-$2$ coefficients, which gives the isomorphism analogous 
  to \eqref{eq:local-coeff} for the Novikov ring $\ring$ of characteristic $2$.
\end{proof}

To simplify notation, we shall omit the basepoint $p$ from the notation 
involving $\KHGm$ in the sequel; however, we emphasize again that the basepoint 
is a necessary input for naturality results that allow $\KHGm$ to be well 
defined.

\subsection{\texorpdfstring{$\taushG$}{Tau-sharp}}
\label{ssec:tau-sharp}

We now recall the definition of $\taushG$ by Baldwin and Sivek 
\cite{baldwin2020framed}. For simplicity of notation, we first focus on 
$\taushI$.  First, for a knot $K \subset S^3$, let $N (K)$ be the smallest 
integer $n \geq 0$ for which the cobordism map
\[
  \Ish (X_n, \nu_n) \colon \Ish (S^3) \to \Ish (S^3_n (K))
\]
vanishes, where $X_n$ is the trace of $n$-surgery along $K$, and $\nu_n$ 
(unrelated to $\nushI$) is some properly embedded surface in $X_n$.

\begin{definition}[{\cite[Definition~3.5]{baldwin2020framed}}]
  \label{def:nu-sharp}
  For a knot $K \subset S^3$, define $\nushI (K) \in \Z$ by the equation
  \[
    \nushI (K) = N (K) - N (\mir{K}).
  \]
\end{definition}

It is proved \cite[Theorem~3.7]{baldwin2020framed} that $\nushI (K)$ depends 
only on the smooth concordance class of $K$, and satisfies the smooth $4$-genus 
bound
\[
  \abs{\nushI (K)} \leq \max (2 g_4 (K) - 1, 0).
\]
It is then shown \cite[Theorem~5.1]{baldwin2020framed} that $\nushI$ defines a 
quasi-morphism from the smooth concordance group $\mathcal{C}$ to $\Z$, i.e.\ 
that $\nushI$ satisfies
\[
  \abs{\nushI (K_1 \connsum K_2) - \nushI (K_1) - \nushI (K_2)} \leq 1
\]
for all knots $K_1, K_2 \subset S^3$, and subsequently make the following 
definition:

\begin{definition}
  \label{def:tau-sharp}
  For a knot $K \subset S^3$, define $\taushI (K) \in \R$ as the homogenization
  \[
    \taushI (K) = \frac{1}{2} \lim_{n \to \infty} \frac{\nushI (\bigconnsum n 
      K)}{n}.
  \]
\end{definition}

One then has \cite[Proposition~5.4]{baldwin2020framed} that this concordance 
invariant defines a group homomorphism $\taush \colon \mathcal{C} \to \R$ and 
satisfies the smooth $4$-genus bound
\[
  \abs{\taushI (K)} \leq g_4 (K).
\]

Finally, $\taushM (K)$ can be defined completely analogously, where $N (K)$ 
would instead be the smallest integer $n \geq 0$ such that the cobordism map
\[
  \HMt (X_n) \colon \HMt (S^3) \to \HMt (S^3_n (K))
\]
vanishes.

%% file: sec-tau-conc.tex
\section{Concordance invariance of \texorpdfstring{$\tauG$}{tau}}
\label{sec:tau-conc}

In this section, we prove the corcordance invariance of $\tauG$, establishing 
\fullref{prop:tau-conc}.
Throughout the section, we have a knot $K\subset S^3$ and the sutures $\Ga_{n}$ 
and $\Ga_{\mu}$ on $\partial{S^3}(K)$, as described in the previous section.

Fix $n\in\intg_+$; on $\partial{S^3(K)}$, we pick a meridional curve $\al$ such 
that $\al$ intersects the sutures $\Ga_n$ twice. (Note that this is unrelated 
to the curve $\al$ in \fullref{fig:twist-knots}.) Let 
$[-1,0]\times\partial{S^3(K)}\subset S^3(K)$ be a collar of $\partial{S^3(K)}$ 
inside the knot complement $S^3(K)$, and endow a $[-1,0]$-invariant tight 
contact structure on $[-1,0]\times\partial{S^3(K)}$, so that each slice 
$\{t\}\times\partial{S^3(K)}$ for $t\in[-1,0]$ is convex and the dividing set 
is (isotopic to) $\Ga_n$. By the Legendrian Realization Principle, we can push 
$\al$ into the interior of the collar $[-1,0]\times\partial{S^3(K)}$ and get a 
Legendrian curve $\be$. With respect to the surface framing, the curve $\be$ 
has $\tb=-1$. (When talking about framings of $\be$, we will always refer to 
the surface framing with respect to $\partial{S^3(K)}$.)

Following Baldwin and Sivek \cite{baldwin2016contact}, since $\al$ intersects 
the sutures $\Ga_n$ twice, after making $\al$ Legendrian, we can glue a contact 
$2$-handle to $(S^3(K),\Ga_n)$ along $\al$, and get a new balanced sutured 
manifold $(M,\ga)$. Suppose that $(Y,R)$ is a closure of $(S^3(K),\Ga_n)$ in 
the sense of Kronheimer and Mrowka \cite{kronheimer2010knots} such that $g(R)$ 
is sufficiently large; then, by work of Baldwin and Sivek 
\cite{baldwin2014invariants}, we know that a closure $(Y_0,R)$ of $(M,\ga)$ can 
be obtained from $(Y,R)$ by performing $0$-surgery along the curve $\be$. Note 
that, inside $Y$, $\be$ is disjoint from $R$, and so the surgery can be made 
disjoint from $R$; this means that the surface $R$ survives in $Y_0$.  Now let 
$(M_{-1},\Ga_n)$ be the balanced sutured manifold obtained from 
$(S^3(K),\Ga_n)$ by performing a $(-1)$-surgery along $\be$.  Note that $\be$ 
is contained in the interior of $S^3(K)$, and so the surgery does not affect 
the boundary or as the sutures. 

Clearly, if we perform $(-1)$-surgery along $\be$ on $Y$, we will get a closure 
$(Y_{-1},R)$ of the balanced sutured manifold $(M_{-1},\Ga_n)$.  The surgery 
exact triangle proved by Kronheimer, Mrowka, Ozsv\'{a}th, and Szab\'{o} 
\cite[Theorem~2.4]{kronheimer2007monopolesandlens}, generalized to the sutured 
setting, gives the exact triangle
\begin{equation*}
  \xymatrix{
    \shg(-M_{-1},-\Ga_n)\ar[rr]&&\shg(-S^3(K),-\Ga_n)\ar[ld]^{C_{h,n}}\\
    &\shg(-M,-\ga).\ar[lu]&
  }
\end{equation*}

\begin{remark}
  \label{rmk:surgery-tri}
  Compared to the one in \cite{kronheimer2007monopolesandlens}, the surgery 
  exact triangle here seems to go in the reverse direction; this is because the 
  orientations on the sutured manifolds have been reversed.
\end{remark}

We now determine that $(M,\ga)$ and $(M_{-1},\Ga_n)$ are familiar balanced 
sutured manifolds.  First, $(M,\ga)$ is obtained from $(S^3(K),\Ga_n)$ by 
attaching a contact $2$-handle along a meridional curve $\al$, and so it is 
nothing but $(S^3(1),\delta)$, where $S^3(1)$ is obtained from $S^3$ by 
removing a $3$-ball, and $\delta$ is a connected simple closed curve on the 
spherical boundary of $S^3 (1)$. For $(M_{-1},\Ga_n)$, note that $\be$ and $K$ 
are inside the $3$-sphere $S^3$, and $\be$ is a meridian around $K$.  Thus, 
$(-1)$-surgery along $\be$ on $S^3(K)$ will result in the same $3$-manifold 
$S^3(K)$, while the framing on its boundary will increase by $1$.
In other words, we have $(M_{-1},\Ga_n)\cong(S^3(K),\Ga_{n-1})$. (Recall that 
the slope of $\Ga_{n}$ is $-n$).  Thus, the above exact triangle becomes
\begin{equation}
  \label{eq: exact triangle to 3-sphere}
  \vcenter{
    \xymatrix{
      \shg(-S^3(K),-\Ga_{n-1})\ar[rr]&&\shg(-S^3(K),-\Ga_n).\ar[ld]^{C_{h,n}}\\
      &\shg(-S^3(1),-\delta)\ar[lu]&
    }
  }
\end{equation}

\begin{lemma}
  \label{lem: the growth rate}
  Denote by $\maxtb(K)$ the maximal Thurston--Bennequin number among all 
  Legendrian representatives $\leg \subset (S^3, \xistd)$ of the smooth knot 
  type $K$.
  If $n\geq-\maxtb(K)$, then the map $C_{h,n}$ is surjective, and hence
  \begin{equation}
    \label{eq: the growth rate}
    \rk_{\ring} \shg(-S^3(K),-\Ga_n) = \rk_\ring \shg(-S^3(K),-\Ga_{n-1}+1.
  \end{equation}
\end{lemma}

\begin{proof}
  Since $n\geq-\maxtb(K)$, we can isotope $K$ to a Legendrian $\leg \subset 
  (S^3, \xistd)$ with $\tb(\leg)=-n$.  We can remove a standard Legendrian 
  neighborhood of $\leg$; then the dividing set on the boundary of the 
  complement is the sutures $\Ga_n$.  Hence, when we glue back a contact 
  $2$-handle, we get $(S^3 (1),\delta)$ with the standard tight contact 
  structure. By work of Baldwin and Sivek 
  \cite{baldwin2016contact,baldwin2016instanton}, we know that the 
  corresponding contact element is a generator of
  \[
    \shg(-S^3(1),-\delta)\cong\mathcal{R}.
  \]
  Since the contact $2$-handle attaching map $C_{h,n}$ preserves the contact 
  element, we see that $C_{h,n}$ is surjective.
\end{proof}

We now digress to prove that there is a unique $\ring [U]$-summand in $\KHGm 
(S^3, K)$.

\begin{proof}[Proof of \fullref{prop:KHGm-tower}]
  Suppose that $S$ is a minimal-genus Seifert surface of $K$, and let $g = 
  g(S)$. The main portion of this proof will be to show that a pattern emerges 
  for $\SHG (-S^3 (K), -\Gamma_n)$ for sufficiently large $n$, with gradings 
  taken into account.  More precisely, we shall use bypass exact triangles to 
  show that the rank of $\SHG (-S^3 (K), -\Gamma_{2g+k})$ increases by a fixed 
  positive integer $r$ whenever the nonnegative integer $k$ increases by $1$, 
  as expressed in the following: For $k$ odd, we have
  \begin{equation}
    \label{eq: the homology group for large n, 1}
    \begin{aligned}
      & \quad \shg(-S^3(K),-\Ga_{2 g +k},S_{2 g +k},i)\\
      & \isom \begin{cases}
        0 & \text{for } i > 2 g  + (k-1)/2,\\
        \shg(-S^3,-\Ga_{2 g },S_{2 g }^-,i - (k-1)/2) & \text{for } (k+1)/2 
        \leq i \leq 2 g + (k-1)/2,\\
        \ring^r & \text{for } -(k-1)/2 \leq i \leq (k-1)/2,\\
        \shg(-S^3,-\Ga_{2 g },S_{2 g }^-,i + (k+1)/2) & \text{for } - 2g - 
        (k-1)/2 \leq i\leq - (k+1)/2,\\
        0 & \text{for }i < -2 g  - (k-1)/2;
      \end{cases}
    \end{aligned}
  \end{equation}
  while for $k$ even, we have
  \begin{equation}\label{eq: the homology group for large n, 2}
    \begin{aligned}
      & \quad \shg(-S^3(K),-\Ga_{2 g +k},S_{2 g +k}^-,i)\\
      & \isom \begin{cases}
        0 & \text{for } i > 2 g  + k / 2,\\
        \shg(-S^3,-\Ga_{2 g },S_{2 g }^-,i - k/2) & \text{for } k/2 + 1 \leq i 
        \leq 2 g + k/2,\\
        \ring^r & \text{for } -k/2 + 1 \leq i \leq k/2,\\
        \shg(-S^3,-\Ga_{2 g },S_{2 g }^-,i + k/2) & \text{for } - 2g - k/2 + 1 
        \leq i\leq - k/2,\\
        0 & \text{for }i < -2 g  - k/2 + 1.
      \end{cases}
    \end{aligned}
  \end{equation}

  To begin, as described in \cite[p.~13]{li2019direct}, the maps 
  $\psi^n_{\pm,n+1}$ fit into bypass exact triangles proved by Baldwin and 
  Sivek \cite[Theorem~1.21]{baldwin2018khovanov}:
  \begin{equation}
    \label{eq:bypass}
    \vcenter{
      \xymatrix{
        \shg(-S^3(K),-\Ga_{n-1})\ar[rr]^{\psi^{n-1}_{\pm,n}}&&\shg(-S^3(K),-\Ga_n)\ar[ld]^{\psi^n_{\pm,\mu}}\\
        &\shg(-S^3(K),-\Ga_{\mu})\ar[lu]^{\psi^\mu_{\pm,n-1}}.
      }
    }
  \end{equation}
  (Note that \eqref{eq:bypass} is in fact two different bypass exact triangles, 
  one for positive bypasses and one for negative bypasses, written together.  
  The same is true for \eqref{eq:bypass-gr-odd} and \eqref{eq:bypass-gr-even} 
  below.)
  Examining the proof of \cite[Proposition~5.5]{li2019direct}, one obtains the 
  graded versions of the exact triangles above: Let $S_n$ and $S_\mu$, as well 
  as their positive and negative stabilizations, be as in 
  \fullref{ssec:khgm-tau}; then, for $n$ odd, we have
  \begin{equation}
    \label{eq:bypass-gr-odd}
    \vcenter{
      \xymatrix{
        \shg(-S^3(K),-\Ga_{n-1},S^{\pm}_{n-1},i)\ar[rrrr]^{\psi^{n-1}_{\pm,n}\quad} 
        &&&& \shg(-S^3(K),-\Ga_n,S_n,i)\\
        {}\save[]+<1.7in,0in>*+{\shg(-S^3(K),-\Ga_{\mu},S_{\mu}^{\mp n \pm 
            1},i);}
        \ar[u]^>>>>>>>{\psi^{\mu}_{\pm,{n-1}}}
        \ar@{<-}[urrrr]_>>>>>>>{\psi^{n}_{\pm,\mu}}
        \restore
      }
    }
  \end{equation}
  while for $n$ even, we have
  \begin{equation}
    \label{eq:bypass-gr-even}
    \vcenter{
      \xymatrix{
        \shg(-S^3(K),-\Ga_{n-1},S^{\pm2}_{n-1},i)\ar[rrrr]^{\psi^{n-1}_{\pm,n}\quad} 
        &&&& \shg(-S^3(K),-\Ga_n,S_n^{\pm},i)\\
        {}\save[]+<1.7in,0in>*+{\shg(-S^3(K),-\Ga_{\mu},S_{\mu}^{\mp n \pm 
            2},i).}
        \ar[u]^>>>>>>>{\psi^{\mu}_{\pm,{n-1}}}
        \ar@{<-}[urrrr]_>>>>>>>{\psi^{n}_{\pm,\mu}}
        \restore
      }
    }
  \end{equation}

  We shall in general be applying \eqref{eq:bypass-gr-odd} and 
  \eqref{eq:bypass-gr-even} with $n = 2g + k$, where $k > 0$. The key 
  observation is that the homology group in the bottom rows of 
  \eqref{eq:bypass-gr-odd} and \eqref{eq:bypass-gr-even} is zero for many 
  gradings $i$, which give us an isomorphism in the top row. Precisely, it is 
  well known (for example, see \cite{kronheimer2010knots}) that for $\abs{i} > 
  g$,
  \[
    \shg(-S^3(K),-\Ga_{\mu},S_{\mu},i)=0,
  \]
  and so by the grading shift in \eqref{eq:gr-shift}, for $k$ odd, we have
  \begin{align*}
    \shg(-S^3(K),-\Ga_{\mu},S^{-2g-k+1}_{\mu},i)& =0 &\text{for } i < -g + (g + 
    (k-1)/2) &= (k-1)/2,\\
    \shg(-S^3(K),-\Ga_{\mu},S^{2g+k-1}_{\mu},i)& =0 &\text{for } i > g + (-g - 
    (k-1)/2) &= -(k-1)/2;
  \end{align*}
  thus, the positive (resp.\ negative) bypass exact triangle in 
  \eqref{eq:bypass-gr-odd} splits for $i < (k-1)/2$ (resp.\ $i > -(k-1)/2$), 
  and we obtain, for $i < (k-1)/2$,
  \begin{equation}
    \label{eq:shg-isom-odd-1}
    \begin{aligned}
      \SHG (-S^3 (K), -\Gamma_{2g+k}, S_{2g+k}, i) & \isom \SHG (-S^3 (K), 
      -\Gamma_{2g+k-1}, S_{2g+k-1}^+, i)\\
      & \isom \SHG (-S^3 (K), -\Gamma_{2g+k-1}, S_{2g+k-1}^-, i + 1)
    \end{aligned}
  \end{equation}
  where the last isomorphism follows also from \eqref{eq:gr-shift}, and for $i 
  > -(k-1)/2$,
  \begin{equation}
    \label{eq:shg-isom-odd-2}
    \SHG (-S^3 (K), -\Gamma_{2g+k}, S_{2g+k}, i) \isom \SHG (-S^3 (K), 
    -\Gamma_{2g+k-1}, S_{2g+k-1}^-, i).
  \end{equation}
  Similarly, for $k$ even, the positive and negative bypass exact triangles in 
  \eqref{eq:bypass-gr-even} respectively give, for $i - 1 < (k-2)/2$ (i.e.\ for 
  $i < k/2$),
  \begin{equation}
    \label{eq:shg-isom-even-1}
    \begin{aligned}
      \SHG (-S^3 (K), -\Gamma_{2g+k}, S_{2g+k}^-, i) & \isom \SHG (-S^3 (K), 
      -\Gamma_{2g+k}, S_{2g+k}^+, i - 1)\\
      & \isom \SHG (-S^3 (K), -\Gamma_{2g+k-1}, S_{2g+k-1}^{2}, i - 1)\\
      & \isom \SHG (-S^3 (K), -\Gamma_{2g+k-1}, S_{2g+k-1}, i),
    \end{aligned}
  \end{equation}
  and for $i > -(k-2)/2 = -k/2 + 1$,
  \begin{equation}
    \label{eq:shg-isom-even-2}
    \begin{aligned}
      \SHG (-S^3 (K), -\Gamma_{2g+k}, S_{2g+k}^-, i) &\isom \SHG (-S^3 (K), 
      -\Gamma_{2g+k-1}, S_{2g+k-1}^{-2}, i)\\
      & \isom \SHG (-S^3 (K), -\Gamma_{2g+k-1}, S_{2g+k-1}, i-1).
    \end{aligned}
  \end{equation}

  Now for $k$ odd, by setting $n = 2g+k$ in \cite[Theorem~2.20~(1)]{li-ye}, we 
  see that the $\ring$-module $\SHG (-S^3 (K), -\Gamma_{2g+k}, S_{2g+k}, i)$ is 
  supported only in gradings $-2g - (k-1)/2 \leq i \leq 2g + (k-1)/2$.  This, 
  together with \eqref{eq:shg-isom-even-1} and \eqref{eq:shg-isom-even-2}, 
  implies that for $k$ even, $\SHG (-S^3 (K), -\Gamma_{2g+k}, S_{2g+k}^-, i)$ 
  is supported only in gradings $-2g - k/2 + 1 \leq i \leq 2g + k/2$. We call 
  these the \emph{possible gradings}.

  Therefore, in essence, what \eqref{eq:shg-isom-odd-1} and 
  \eqref{eq:shg-isom-even-1} say is that the summands of $\SHG (-S^3 (K), 
  -\Gamma_{2g+k})$ in the bottom $2g + k - 1$ possible gradings are 
  respectively isomorphic to the summands of $\SHG (-S^3 (K), 
  -\Gamma_{2g+k-1})$ in the bottom $2g + k - 1$ possible gradings (possibly 
  with a grading shift), and \eqref{eq:shg-isom-odd-2} and 
  \eqref{eq:shg-isom-even-2} give the analogous statement for the top $2g + k - 
  1$ possible gradings.  (Since $\SHG (-S^3 (K), -\Gamma_{2g+k-1})$ has $4g + k 
  - 1$ possible gradings, this means that summands in $2 (2g + k - 1) - (4g + k 
  -1) = k-1$ ``middle'' gradings are ``sampled'' twice.) As $\SHG (-S^3 (K), 
  -\Gamma_{2g+k})$ has $4g + k$ possible gradings, these isomorphisms 
  completely determine $\SHG (-S^3 (K), -\Gamma_{2g+k})$ in terms of $\SHG 
  (-S^3 (K), -\Gamma_{2g+k-1})$, except in the case $k = 1$, wherein the middle 
  grading $\SHG (-S^3 (K), -\Gamma_{2g+1}, S_{2g+1}, 0)$ is not determined.  
  Simply letting
  \[
    r = \rk_{\ring} \SHG (-S^3 (K), -\Gamma_{2g+1}, S_{2g+1}, 0),
  \]
  we establish \eqref{eq: the homology group for large n, 1} and \eqref{eq: the 
    homology group for large n, 2} by inducting on $k$.

  Now, taking $k$ to be sufficiently large, \eqref{eq: the growth rate} implies 
  that the integer $r$ must in fact be $1$. Thus, for a fixed grading $i = -g 
  -m \leq -g$, we have that
  \begin{align*}
    & \quad \SHG (-S^3 (K), -\Gamma_{2g+k}, S_{2g+k}^{\tau(2g+k)}, i) 
    [\sigma]\\
    & \isom
    \begin{cases}
      \SHG (-S^3 (K), -\Gamma_{2g+k}, S_{2g+k}^{\tau(2g+k)}, -m + (k-1)/2) & 
      \text{if } k \text{ is odd},\\
      \SHG (-S^3 (K), -\Gamma_{2g+k}, S_{2g+k}^{\tau(2g+k)}, -m + k/2) & 
      \text{if } k \text{ is even},
    \end{cases}\\
    & \isom \ring
  \end{align*}
  whenever $k \geq m + 1$, where the last isomorphism follows from \eqref{eq: 
    the homology group for large n, 1} and \eqref{eq: the homology group for 
    large n, 1}. By definition, this means that
  \[
    \KHGm (-S^3, K, i) \isom \ring
  \]
  for all $i \leq -g$. Finally, this together with 
  \cite[Corollary~5.11]{li2019direct} imply that there is a submodule (and at 
  most one such submodule) in $\KHGm (-S^3, K)$ isomorphic to $\ring [U]$.  
  Since $\ring$ is a field in our context, we conclude that it is in fact a 
  unique $\ring [U]$-summand.
\end{proof}

In the following, we will continue to denote by $g$ the genus $g (K)$ of a knot 
$K$.

Strictly speaking, we did not have to prove \fullref{prop:KHGm-tower} for the 
arguments of this section. Its significance, however, is that it explains the 
definition of $\tauG$ as a natural, unique definition for knots in $S^3$.

Having proved that $\KHGm (-S^3, K)$ has a unique $\ring [U]$-summand, we now 
return to the main set-up of the section to prove the concordance invariance of 
$\tauG$. Our next major goal is to re-characterize $\tauG$ in terms of the 
(non-)vanishing of a map on $\KHGm (-S^3, K)$ induced by the maps $C_{h,n}$. We 
begin with the following lemma:

\begin{lemma}
  \label{lem: C_h is a map on KHG minus}
  The maps
  \[
    C_{h,n} \colon \shg(-S^3(K),-\Ga_n)\ra\shg(-S^3(1),-\delta),
  \]
  which appear in the exact triangle \eqref{eq: exact triangle to 3-sphere}, 
  induce a surjective map
  \[
    C_h \colon \khg(-S^3,K)\ra\shg(-S^3(1),-\delta).
  \]
  Furthermore, $C_h$ commutes with the action of $U$on $\KHGm (-S^3, K)$.
\end{lemma}

\begin{proof}
  The lemma follows from \fullref{lem: the growth rate} and the following two 
  commutative diagrams, one for positive bypasses and one for negative 
  bypasses:
  \begin{equation*}
    \label{eq:C_h,n-comm}
    \vcenter{
      \xymatrix{
        \shg(-S^3,-\Ga_n)\ar[rr]^{\psi^n_{\pm,n+1}}\ar[dr]_{C_{h,n}}&&\shg(-S^3,-\Ga_{n+1})\ar[dl]^{C_{h,n+1}}\\
        &\shg(-S^3(1),-\delta).
      }
    }
  \end{equation*}
  To prove these commutative diagrams, recall that the maps $\psi^n_{\pm,n+1}$ 
  are constructed via bypass attachments, which can be interpreted as contact 
  handle attachments (see \cite[Section~3]{ozbagci} and  
  \cite[Section~5]{baldwin2016contact}), and so is $C_{h,n}$.  The 
  commutativity for the diagrams follows from the observation that the contact 
  handle attaching regions for $\psi^n_{\pm,n+1}$ and $C_{n,h}$ are disjoint 
  from each other.
\end{proof}

As a quick aside, we exhibit an immediate consequence of \fullref{lem: C_h is a 
  map on KHG minus}:

\begin{corollary}
  \label{cor:U-1}
  There is an exact triangle
  \begin{equation*}
    \xymatrix{
      \khg(-S^3,K)\ar[rr]^{\Psi}&&\khg(-S^3,K)\ar[dl]^{C_h}\\
      &\shg(-S^3(1),-\delta)\ar[lu]&
    }
  \end{equation*}
\end{corollary}

\begin{proof}
  The maps $C_{h,n}$ in the exact triangle \eqref{eq: exact triangle to 
    3-sphere} commute with the maps $\psi^n_{-,n+1}$ in the directed system, 
  and so we can pass to the direct limit and still have an exact triangle.
\end{proof}

The significance of \fullref{cor:U-1} is the following. There is an exact 
triangle in Heegaard Floer theory that involves the modules that are analogous 
to those appearing in \fullref{cor:U-1}:
\[
  \xymatrix{
    \HFKm(-S^3,K)\ar[rr]^{U-1}&&\HFKm(-S^3,K)\ar[dl]\\
    &\HFh(-S^3)\ar[lu]&
  }
\]
One key difference is that, in this context, the map in the top row is defined 
algebraically. Thus, we are led to ask the following natural question:

\begin{question}
  Does $\Psi$ admit an interpretation as $U-1$, where $U$ denotes the action of 
  $U$?
\end{question}

We believe establishing a positive answer to this question would have 
topological applications.

In any case, we are now ready to re-characterize $\tauG$:

\begin{proposition}
  \label{prop: another definition of tau}
  The invariant $\tauG (K)$ admits an alternative definition:
  \[
    \tauG(K) = \max \setc{i \in \Z}{\text{the restriction of } C_h \text{ to } 
      \khg(-S^3,K,i) \text{ is non-trivial} \,}
  \]
\end{proposition}

\begin{proof}
  We claim that an element $[x] \in \KHGm (-S^3, K)$ is non--$U$-torsion if and 
  only if $C_h ([x]) \neq 0$.
  Let $[x] \in \KHGm (-S^3, K, i)$ be a $U$-torsion element with $i \leq g$ 
  (which is the top possible grading); then it is represented by an element
  \[
    x\in \shg(-S^3(K),-\Ga_{2g+k},S_{2g+k},i)[\sigma] \isom
    \SHG \paren{-S^3 (K), -\Gamma_{2g+k}, S_{2g+k}, i + g + \frac{k-1}{2}}
  \]
  for some odd $k > \max (0, -g -i)$,
  such that there exists an even $l\in\intg_+$ with
  \[
    \psi_{+,2g+k+l}^{2g+k+l-1}\circ \dotsb \circ\psi_{+,2g+k+1}^{2g+k}(x)=0.
  \]
  (Note that by \cite[Proposition 5.10]{li2019direct}, the direct system 
  defining $\KHGm (-S^3, K, i)$ has stabilized at $\shg(-S^3(K),-\Ga_{2g+k})$  
  for $k > -g -i$.
  We choose $k$ to be odd to avoid having to work with stabilizations of the 
  grading surface.
  Also, if the image of $x$ vanishes after an odd number of applications of 
  $\psi_+$, we may simply apply $\psi_+$ once more; we choose $l$ to be even 
  for ease of presentation in the next paragraph.)
  From the commutative diagram in the proof of \fullref{lem: C_h is a map on 
    KHG minus}, we know that
  \[
    C_{h,2g+k}(x)=C_{h,2g+k+l}\circ\psi_{+,2g+k+l}^{2g+k+l-1}\circ\dotsb\circ\psi_{+,2g+k+1}^{2g+k}(x)=0,
  \]
  and so $C_h ([x]) = 0$.

  Conversely, let $[x] \in \KHGm (-S^3, K)$ be a non--$U$-torsion element; then 
  it is represented by an element $x\in \shg(-S^3(K),-\Ga_{2g+k},S_{2g+k},i + g 
  + (k-1)/2)$ such that
  \[
    \psi_{+,2g+k+l}^{2g+k+l-1}\circ \dotsb \circ\psi_{+,2g+k+1}^{2g+k}(x)\neq0
  \]
  for all even $l \in \Z_+$. (Note that this implies a statement for odd $l$ as 
  well.) Taking into account the gradings as in the first rows of 
  \eqref{eq:bypass-gr-odd} and \eqref{eq:bypass-gr-even}, for a given, even $l 
  \in \Z_+$, this is an element of
  \[
    \shg \paren{-S_3(K),-\Gamma_{2g+k+l},S_{2g+k+l},i + g + \frac{k-1}{2} -l}.
  \]
  Now the idea is that, for large $n$, the map $C_{h,n}$ is an isomorphism when 
  restricted to the ``middle'' possible gradings; and we can ensure that our 
  element lies in those ``middle'' gradings by taking $l$ to be sufficiently 
  large.  Precisely, from the proof of \cite[Proposition~4.28]{li-ye}, 
  $C_{h,2g+k+l}$ is an isomorphism when restricted to 
  $\shg(-S^3(K),-\Ga_{2g+k+l},S_{2g+k+l},j)$ for
  \[
    -\frac{k-1}{2}-l\leq j\leq \frac{k-1}{2}+l.
  \]
  Since we chose $k \geq -g -i + 1$, we have that
  \[
    -\frac{k-1}{2}-l \leq i + g + \frac{k-1}{2} - l
  \]
  for all $l$; and if we take $l \geq 2g$, then we will have
  \[
    i + g + \frac{k-1}{2} - l \leq 2g + \frac{k-1}{2} \leq \frac{k-1}{2} + l.
  \]
  Then for these choices, we see that
  \[
    C_{h,2g+k}(x)=C_{h,2g+k+l}\circ\psi_{+,2g+k+l}^{2g+k+l-1}\circ\dotsb\circ\psi_{+,2g+k+1}^{2g+k}(x)\neq0,
  \]
  which implies that $C_h ([x]) \neq 0$. The proposition follows immediately.
\end{proof}

\begin{remark}
  \label{rmk:non-triv-is-surj}
  By the same argument as in the proof of \fullref{lem: only rank is 
    important}, we can show that in \fullref{prop: another definition of tau}, 
  the map $C_h$ being non-trivial is equivalent to it being surjective.
\end{remark}

With the alternative definition of $\tauG$, we can now prove that it is a 
concordance invariant:

\begin{proof}[Proof of \fullref{prop:tau-conc}]
  Suppose $K_0$ and $K_1$ are concordant; then there exists a properly embedded 
  annulus $A\subset [0,1]\times S^3$ such that
  \[
    (\{0\}\times S^3,A\cap \{0\}\times S^3)\cong(S^3,K_0),~(\{1\}\times 
    S^3,A\cap \{1\}\times S^3)\cong(S^3,K_1).
  \]
  The idea of the proof is that $A$ induces a grading-preserving cobordism map 
  $F_A \colon \KHGm (-S^3, K_0) \to \KHGm (-S^3, K_1)$ that commutes with 
  $C_h$, which will imply the result for $\tauG$ via \fullref{prop: another 
    definition of tau}.

  The first step is to analyze the cobordism map induced by $A$ on $\SHG (-S^3 
  (K_0), \Gamma_n)$.
  For each $n$, the pair $([0,1]\times S^3,A)$ induces a cobordism $W_n$ from 
  $Y_{0,n}$ to $Y_{1,n}$, where $Y_{i,n}$ is a closure of $(-S^3(K_i),-\Ga_n)$, 
  and $W_n$ induces a map
  \[
    F_{A,n} \colon \shg(-S^3(K_0),-\Ga_n)\ra \shg(-S^3(K_1),-\Ga_{n})
  \]
  as follows. There are two ways to describe $W_n$, which are both useful; 
  below, we briefly recall both of these descriptions from \cite{li2018gluing}.

  First, take a parametrization of $A \isom [0,1]\times S^1$. Then, a tubular 
  neighborhood of $A\subset [0,1]\times S^3$ can be identified with $A\times 
  D^2 \isom [0,1]\times S^1\times D^2$, with
  \[
    (A\times D^2) \cap (\set{0,1} \times S^3) \cong \{0,1\}\times S^1\times 
    D^2.
  \]
  Thus, we know that
  \begin{equation}
    \label{eq:A-D^2-comp}
    \partial \paren{([0,1]\times S^3)\setminus (A\times D^2)} \isom
    -S^3(K_0)\cup ([0,1]\times S^1\times\partial D^2)\cup S^3(K_1).
  \end{equation}
  Choosing a closure $Y_{0,n}$ of $(-S^3 (K_0),-\Ga_n)$,
  we can write
  \begin{equation*}
    \label{eq:W_n}
    W_n \isom -\paren{([0,1]\times S^3) \setminus (A\times D^2)}\cup 
    \paren{[0,1] \times (Y_{0,n}\setminus S^3(K_0))},
  \end{equation*}
  via a natural identification
  \[
    [0,1]\times S^1\times\partial D^2 \isom [0,1]\times\partial S^3(K_0).
  \]

  A second description of $W_n$ is as follows. As
  \[
    \bdy S^3 (K_0) \isom \bdy S^3 (K_1) \isom S^1 \times D^2,
  \]
  from \eqref{eq:A-D^2-comp}, $([0,1]\times S^3)\setminus (A\times D^2)$ can be 
  obtained from $([0,1]\times S^3(K_1))$ by attaching a set of $4$-dimensional 
  handles $\mathcal{H}$ to the interior of $\set{1} \times S^3(K_0)$,
  as in \cite[Lemma~3.3]{li2018gluing}. Thus, as above,
  choosing
  a closure $Y_{0,n}$ of $(-S^3(K_0),-\Ga_n)$, we can attach the same set of 
  handles $\mathcal{H}$ to $\{1\}\times Y_{0,n}\subset [0,1]\times Y_{0,n}$, 
  and the result is
  again $W_n$.

  We break down the rest of the proof into four claims, as detailed below.

  \emph{Claim 1.} The maps $F_{A,n}$ give rise to a map
  \[
    F_A \colon \khg(-S^3,K_0)\ra \khg(-S^3,K_1).\footnote{The basepoints $p_i$ 
      for $K_i$ are specified by $p_i = \set{i} \times p$ in the 
      parametrization $A \isom [0, 1] \times S^1$.}
  \]
  To prove the claim, it suffices to show that we have a commutative diagram
  \begin{equation*}
    \xymatrix{
      \shg(-S^3(K_0),-\Ga_n)\ar[r]^{F_{A,n}}\ar[d]^{\psi^n_{-,n+1}}&\shg(-S^3(K_1),-\Ga_n)\ar[d]^{\psi^n_{-,n+1}}\\
      \shg(-S^3(K_0),-\Ga_{n+1})\ar[r]^{F_{A,n+1}}&\shg(-S^3(K_1),-\Ga_{n+2}).
    }
  \end{equation*}
  The commutativity of this diagram follows from the fact that the attaching 
  regions for the handles associated to $F_{A,n}$ and to $\psi^n_{-,n+1}$ are 
  disjoint: When constructing $F_{A,n}$, we attached handles to $[0,1]\times 
  Y_{0,n}$ along the region $\{1\}\times\Int S^3(K_0)$, while when constructing 
  the map $\psi^n_{-,n+1}$, we attached handles to $[0,1]\times Y_{i,n}$ along 
  the region $\set{1} \times [0, 1] \times \partial S^3(K_i)$ see 
  \cite[Section~3]{li2018gluing}).

  \emph{Claim 2.} $F_A$ commutes with the $U$ map on $\khg$.
  The proof of this claim is completely analogous to one for Claim~1, with 
  $\psi_+$ instead of $\psi_-$.

  The two claims above show that $F_A$ is a homomorphism of $\ring 
  [U]$-modules.

  \emph{Claim 3.} There is a commutative diagram
  \begin{equation*}
    \xymatrix{
      \khg(-S^3,K_0)\ar[rr]^{F_A}\ar[dr]_{C_{h}}&&\khg(-S^3,K_1)\ar[dl]^{C_{h}}\\
      &\shg(-S^3(1),-\delta),&
    }
  \end{equation*}
  where $C_h$ is defined as in \fullref{lem: C_h is a map on KHG minus}.

  To prove the claim, it suffices to prove that the following digram commutes 
  for all $n$:
  \begin{equation}
    \label{eq: commutative diagram for F_A, C_h}
    \vcenter{
      \xymatrix{
        \shg(-S^3(K_0),-\Ga_n)\ar[r]^{F_{A,n}}\ar[d]^{C_{h,n}}&\shg(-S^3(K_1),-\Ga_n)\ar[d]^{C_{h,n}}\\
        \shg(-S^3(1),-\delta)\ar[r]^{\Id}&\shg(-S^3(1),-\delta).
      }
    }
  \end{equation}
  As above, suppose we have a closure $Y_{0,n}$ for $(-S^3(K_0),-\Ga_{n})$. Let 
  $Y_{1,n}$ be the corresponding closure for $(-S^3(K_1),-\Ga_{n})$ as in the 
  construction of $W_n$ above. Recall from the construction of $C_{h,n}$ that 
  it is the map associated to a $2$-handle attached along a meridian curve 
  $\al\subset\partial S^3(K_0)$; we can push $\alpha$ slightly into the 
  interior and get a curve $\be$. Then we get a closure $Y'_{0}$ for 
  $(-S^3(1),-\delta)$ by performing $0$-surgery on $Y_{0,n}$ along $\be$.  Note 
  the difference between $S^3(K_0)$ and $S^3(K_1)$ is contained in the 
  interior, and so we also have the curve $\be\subset S^3(K_1)\subset Y_{1,n}$.  
  Thus, we can obtain another closure $Y_1'$ for $(-S^3(1),-\delta)$.  We can 
  form a cobordism $W_n'$ from $Y_0'$ to $Y_1'$ by attaching the set of 
  $4$-dimensional handles $\mathcal{H}$ as in the proof of Claim~1 to 
  $Y_0'\times\{1\}\subset Y_0'\times[0,1]$, and the attaching region is 
  contained in $\Int(S^3(K_0))\subset Y_0'$. Hence, there is a commutative 
  diagram just as in the proof of Claim~1:
  \begin{equation*}
    \xymatrix{
      \shg(-S^3(K_0),-\Ga_n)\ar[r]^{F_{A,n}}\ar[d]^{C_{h,n}}&\shg(-S^3(K_1),-\Ga_n)\ar[d]^{C_{h,n}}\\
      \shg(-S^3(1),-\delta)\ar[r]^{F_A'}&\shg(-S^3(1),-\delta),
    }
  \end{equation*}
  where $F_A'$ is the map induced by the cobordism $W_n'$.
  
  So to prove \eqref{eq: commutative diagram for F_A, C_h}, it suffices to show 
  that $W_n'$ is actually a product $[0,1]\times Y_0'$, which will imply that 
  $F_A'=\Id$.  To do this, recall that $W_n'$ is obtained from $[0,1] \times 
  Y_0'$ by attaching a set of handles $\mathcal{H}$, while the attachment 
  regions are contained in $\Int S^3(K_0)\subset \Int S^3(1) \subset 
  \{1\}\times Y_0'$.  This means that we can split $W_n'$ into two parts
  \[
    W_n' \isom W_n''\cup \paren{[0,1]\times(Y_0'\setminus S^3(1))},
  \]
  where $W_n''$ is obtained from $[0,1]\times S^3(1)$ by attaching the set of 
  handles $\mathcal{H}$. Recall that $(S^3(1),\delta)$ is obtained from 
  $(S^3(K_0),\Ga_n)$ by attaching the contact $2$-handle $h$, and so 
  topologically,
  \[
    S^3(1) \isom S^3(K_0)\cup B^3.
  \]
  Note the $3$-ball $B^3$ is attached to $S^3(K_0)$ along part of the boundary, 
  and the set of handles $\mathcal{H}$ is attached to $[0,1]\times S^3(1)$ 
  within the region $\Int (S^3(K_0))\subset \{1\}\times S^3(1)$, and so the two 
  attaching regions are disjoint. Thus, we have
  \begin{align*}
    W_n''&\isom[0,1]\times S^3(1)\cup\mathcal{H}\\
    &\isom \paren{[0,1]\times(S^3(K_0)\cup B^3)}\cup\mathcal{H}\\
    &\isom([0,1]\times S^3(K_0))\cup\mathcal{H}\cup ([0,1]\times B^3)\\
    &\isom \paren{([0,1]\times S^3)\setminus (A\times D^2)} \cup ([0,1]\times 
    B^3).
  \end{align*}
  Here, $[0,1]\times B^3$ is glued to $([0,1]\times S^3)\setminus (A\times 
  D^2)$ along a thickened annulus.  From here, it is straightforward to check 
  that the resulting manifold $W_n''$ is diffeomorphic to $[0,1]\times S^3(1)$.

  \emph{Claim 4.} The map
  \[
    F_A \colon \khg(-S^3,K_0)\to\khg(-S^3,K_1)
  \]
  preserves the grading.

  By definition, we know that for any fixed $j\in \intg$, we can pick a large 
  enough odd $n$ so that, for $i=0,1$,
  \[
    \khg(-S^3,K_i,j) \isom \shg 
    \paren{-S^3(K_i),-\Ga_{i,n},S_{i,n},j+\frac{n-1}{2}}.
  \]
  (Here, $\Ga_{i,n}$ is a set of sutures on $-S^3 (K_i)$ of slope $-n$, and 
  $S_{i,n}$ is a minimal-genus Seifert surface of $K_i$ that intersects 
  $\Ga_{i,n}$ at exactly $2n$ points.)
  Hence, to show that $F_A$ preserves the grading, we need only to show that 
  $F_{A,n}$ preserves the grading.
  Note we can identify the boundaries
  \[
    \partial S^3(K_0) \isom \partial S^3(K_1)
  \]
  via the parametrization $A \isom [0,1]\times S^1$, and we can assume that 
  under the above identification,
  \[
    S_{0,n}\cap \partial S^3(K_0) \isom S_{1,n}\cap \partial S^3(K_1).
  \]

  Now let $Y_{0,n}$ be a closure of $(-S^3(K_0),-\Ga_n)$, and let 
  $\overline{S}_{0,n}$ be the closure of $S_{0,n}$ in $Y_{0,n}$, as in the 
  construction of gradings; see \cite[Section~3]{li2019direct}.  Then we have a 
  corresponding closure $Y_{1,n}$ for $(-S^3(K_1),-\Ga_n)$, inside which there 
  is the closure $\overline{S}_{1,n}$ of $S_{1,n}$.  To describe this surface, 
  recall that
  \[
    Y_{1,n} \isom -S^3(K_1) \cup_{\partial S^3(K_0) \isom \partial S^3(K_1)} 
    (Y_{0,n}\setminus S^3(K_0))
  \]
  as in the construction of $W_n$ at the beginning of the proof; then 
  concretely, $\overline{S}_{1,n}$ is defined to be
  \[
    \overline{S}_{1,n} = S_{1,n}\cup(\overline{S}_{0,n}\setminus S^3(K_0)).
  \]
  Using the Mayer--Vietoris sequence, we see that
  \[
    H_2(([0,1]\times S^3)\setminus(A\times D^2))=0.
  \]
  Therefore, the closed surface $-S_{0,n} \cup A\cup S_{1,n} \subset 
  ([0,1]\times S^3)\setminus(A\times D^2)$ bounds a $3$-chain $c 
  \subset([0,1]\times S^3)\setminus(A\times D^2)$.
  Now inside $W_n$, let
  \[
    d = c \cup \paren{[0,1]\times(\overline{S}_{0,n}\setminus S^3(K_0))},
  \]
  where the two pieces are glued along
  \[
    A \isom [0,1]\times S^1 \isom [0,1]\times{\bdy (\overline{S}_{0,n}\setminus 
      S^3(K_0))}.
  \]
  It is straightforward to check that
  \[
    \partial d \isom -\overline{S}_{0,n}\cup \overline{S}_{1,n}.
  \]
  Hence we conclude that
  \[
    [\overline{S}_{0,n}]=[\overline{S}_{1,n}]\in H_2(W_n),
  \]
  whence it follows that $F_{A,n}$ preserves the grading.

  The four claims above together prove the existence of a grading-preserving 
  homomorphism $F_A \colon \KHGm (-S^3, K_0) \to \KHGm (-S^3, K_1)$ of $\ring 
  [U]$-modules that commutes with the map $C_h$. By \fullref{prop: another 
    definition of tau}, $\tauG$ is the maximum grading for which $C_h$ is 
  non-trivial, and thus our proof is complete.
\end{proof}

Having achieved our main goal of the section, we end it with an application to 
\emph{ribbon concordance}, which is a knot concordance that admits a handle 
decomposition with only $0$-, $1$-, but not $2$-handles. In recent work of 
Daemi, Lidman, Vela-Vick, and the third author \cite{daemi2019obstructions}, it 
is proved that the map on $\KHI$ associated to a ribbon concordance is 
injective.  We may quickly extend this result to $\KHIm$:

\begin{corollary}
  Suppose $A$ is a ribbon concordance from $K_1$ to $K_2$ in $[0,1] \times 
  S^3$; then the map
  \[
    F_A \colon \khg(-S^3,K_0)\to \khg(-S^3,K_1)
  \]
  defined in the proof of \fullref{prop:tau-conc} is injective.
\end{corollary}

\begin{proof}
  By \cite[Theorem~4.4]{daemi2019obstructions}, the map
  \[
    F_{A,n} \colon \shg(-S^3(K_0),-\Ga_n)\ra\shg(-S^3(K_1),-\Ga_n)
  \]
  is injective for all $n \in \Z$. Passing to the direct limit, we see that 
  $F_A$ is also injective.
\end{proof}

These results may be compared to that of Zemke \cite{zemke2019ribbon}, who 
first proves the analogous statement for both $\HFKh$ and $\HFKm$.

%% file: sec-tau-conn-sum.tex
\section{Additivity of \texorpdfstring{$\tau$}{tau} under connected sum}
\label{sec:tau-conn-sum}

In this subsection, we prove the additivity of the $\tauG$ under connected sum, 
establishing \fullref{prop:tau-conn-sum}.
To begin, we establish the superadditivity of $\tauG$.

\begin{proposition}
  \label{prop: sup additivity}
  Suppose $K_0$ and $K_1$ are two knots in $S^3$; then
  \[
    \tau_G(K_0\connsum K_1) \geq \tau_G(K_0)+\tau_G(K_1).
  \]
\end{proposition}

\begin{proof}
  Suppose $K_0$ and $K_1$ are two knots in $S^3$, and suppose $m$ and $n$ are 
  two sufficiently large, odd integers. Suppose further that $S_0$ and $S_1$ 
  are minimal-genus Seifert surfaces of $K_0$ and $K_1$ respectively. We can 
  attach a $1$-handle $h^1$ to connect the two balanced sutured manifolds 
  $(S^3(K_0),\Ga_m)$ and $(S^3(K_1),\Ga_n)$. Let $(M_0,\ga_0)$ be the resulting 
  balanced sutured manifold; then we have
  \begin{equation}
    \label{eq: C h 1}
    C_{h^1} \colon 
    \shg(-S^3(K_0),-\Ga_{m})\otimes\shg(-S^3(K_1),-\Ga_n)\xra{\cong}\shg(-M_0,-\ga_0).
  \end{equation}

  On $(M_0,\ga_0)$, we can attach a contact $2$-handle $h^2_2$ along the curve 
  $\al$, as depicted in \fullref{fig: contact handles}, and the resulting 
  balanced sutured manifold is $(S^3(K_0\connsum K_1),\Ga_{m+n})$. (This $\al$ 
  is not the same as the one in \fullref{fig:twist-knots}.) Thus, there is a 
  map
  \[
    C_{h^2_2} \colon \shg(-M_0,-\ga_0)\ra \shg(-S^3(K_0\connsum 
    K_1),-\Ga_{m+n}).
  \]

  \begin{figure}[htbp]
    \captionsetup{aboveskip={\dimexpr10pt+2pt+\sactualfontsize\relax}}
    \labellist
    \small\hair 2pt
    \pinlabel {$S^3 (K_0)$} [t] at 59 0
    \pinlabel {$S^3 (K_1)$} [t] at 366 0
    \pinlabel {\textcolor{red}{$\ga_0$}} [l] at 95 300
    \pinlabel {\textcolor{red}{$\ga_1$}} [l] at 388 300
    \pinlabel {\textcolor{blue}{$\al$}} [t] at 366 167
    \pinlabel {\textcolor{blue}{$\mu_0$}} [t] at 59 49
    \pinlabel {\textcolor{blue}{$\mu_1$}} [t] at 366 49
    \pinlabel {$h^1$} [t] at 210 136
    \endlabellist
    \includegraphics[width=0.5\textwidth]{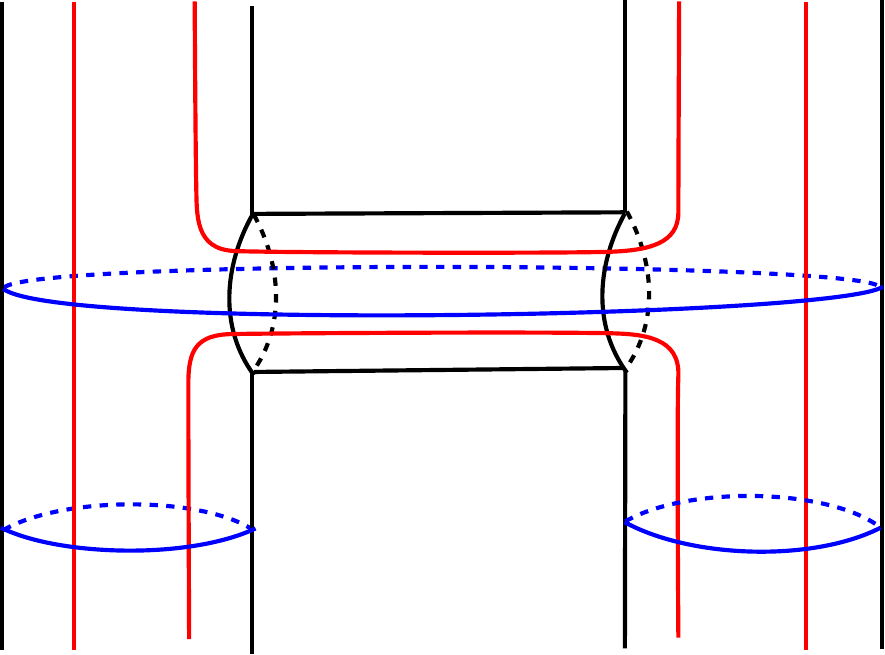}
    \caption{The contact handles: $h^2_0$ is attached along $\mu_0$, $h^2_1$ is 
      attached along $\mu_1$, and $h^2_2$ is attached along $\al$.}
    \label{fig: contact handles}
  \end{figure}

  Inside $(M_0,\ga_0)$, there is a surface $S_0 \sqcup S_1$, whose associated 
  grading is the one we are interested in.  However, the surface $S_0 \sqcup 
  S_1$ intersects the curve $\al$, along which we attach the $2$-handle 
  $h^2_2$, and so it does not survive in $(S^3(K_0\connsum K_1),\Ga_{m+n})$ as 
  a properly embedded surface.  To circumvent this problem, we add to it a 
  strip $P$, as described in the next paragraph.
  
  \begin{figure}[htbp]
    \captionsetup{aboveskip={\dimexpr10pt+2pt+\sactualfontsize\relax}}
    \labellist
    \small\hair 2pt
    \pinlabel {$\bdy S_0$} [t] at 64 0
    \pinlabel {$\bdy S_1$} [t] at 367 0
    \pinlabel {\textcolor{blue}{$\al$}} [ ] at 393 157
    \pinlabel {\rule{0.5pt}{0.13in}} [t] at 334 165
    \pinlabel {$P$} [ ] at 334 131
    \pinlabel {\rule{0.5pt}{0.13in}} [t] at 882 165
    \pinlabel {$S$} [ ] at 882 131
    \endlabellist
    \includegraphics[width=\textwidth]{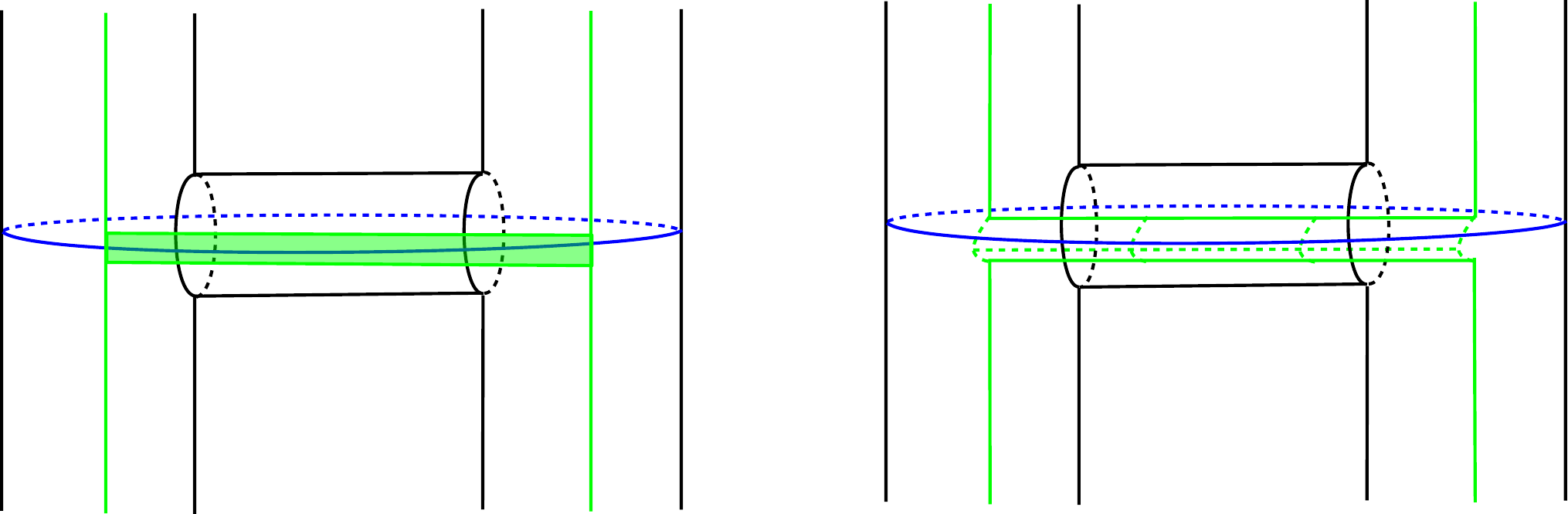}
    \caption{Left: The strip $P$ as the shaded region. Right: Pushing off the 
      interior of $P$ into the interior of $M_0$ in the construction of $S$.}
    \label{fig: adding a strip}
  \end{figure}

  See \fullref{fig: adding a strip}. Pick a strip $P\subset \partial M_0$, 
  which serves as a $2$-dimensional $1$-handle attached to the surfaces $S_0$ 
  and $S_1$. Let $S$ be the union $S_0 \cup S_1\cup P$, with the interior of 
  $P$ being pushed off into the interior of $M_0$; then $S$ is a properly 
  embedded surface inside $(M_0,\ga_0)$ and is disjoint from $\al$.  Thus, 
  after attaching the contact $2$-handle along $\al$, $S$ survives in 
  $S^3(K_0\connsum K_1,\Ga_{m+n})$, and it is obvious that $S$ is a Seifert 
  surface of $K_0\connsum K_1$. Since $\al\cap S=\emptyset$, the map 
  $C_{h^2_2}$ preserves the gradings induced by $S$ and its stabilizations.  To 
  compare the gradings induced by $S_0 \sqcup S_1$ and $S$, note that their 
  difference, the $2$-dimensional $1$-handle $P$, is chosen to be on 
  $\partial{M}_1$.  Hence, we know that
  \[
    [S_0,\partial S_0]+[S_1,\partial{S}_2]=[S,\partial S]\in H_2(M_0,\partial 
    M_0).
  \]
  In \cite{li2019decomposition}, the first and second authors prove that the 
  gradings induced by $S_0 \sqcup S_1$ and $S$ differ by an overall grading 
  shift. To pin down the exact grading shift, observe that the decomposition of 
  $(M,\ga)$ along $S_0 \sqcup S_1$ and $S$ are both taut; this fact allows us 
  to identify the maximal non-vanishing gradings.  Thus, combining with the 
  fact that $C_{h^2_2}$ preserves the grading, we have the following lemma.

  \begin{lemma}
    \label{lem: grading preserving}
    Suppose $m$ and $n$ are sufficiently large, odd integers. Then, for all
    $i,j\in\intg$, the map $C_{h^2_2} \circ C_{h^1}$ shifts the grading as 
    follows:
    \begin{align*}
      C_{h^2_2}\circ C_{h^1} & \colon 
      \shg(-S^3(K_0),-\Ga_m,S_0,i)\otimes\shg(-S^3(K_1),-\Ga_n,S_1,j)\\
      & \quad \ra \shg(-S^3(K_0\connsum K_1),-\Ga_{m+n},S^-,i+j+1). \qedhere
    \end{align*}
  \end{lemma}

  Let $\mu_0\subset\partial{S^3(K_0)}$ and $\mu_1\subset\partial{S^3(K_1)}$ be 
  meridians of $K_0$ and $K_1$ respectively. See \fullref{fig: contact 
    handles}.  We can attach contact $2$-handles $h^2_0$ and $h^2_1$ along 
  $\mu_0$ and $\mu_1$ respectively; the resulting balanced sutured manifolds 
  are both $(S^3(1),\delta)$.
  Thus, we have maps
  \begin{align*}
    C_{h^2_0} & \colon \shg(-S^3(K_0),-\Ga_m)\ra\shg(-S^3(1),-\delta),\\
    C_{h^2_1} & \colon \shg(-S^3(K_1),-\Ga_n)\ra\shg(-S^3(1),-\delta).
  \end{align*}
  The curves $\mu_0$ and $\mu_1$ are disjoint from the contact handles $h^1$ 
  and $h^2_2$, and so they survive in $(S^3(K_0\connsum K_1),\Ga_{m+n})$. Both 
  $\mu_0$ and $\mu_1$ become meridians of $K_0\connsum K_1$, and so the contact 
  $2$-handle attaching maps associated to them (viewed as attachment maps from 
  $(-S^3 (K_0 \connsum K_1), -\Gamma_{m+n})$) are the same:
  \[
    C^{\sharp}=C_{h_0^2}=C_{h_1^2} \colon \shg(-S^3(K_0\connsum 
    K_1),-\Ga_{m+n})\ra\shg(-S^3(1),-\delta).
  \]
  The commutativity of contact handle attachments then gives us the following 
  commutative diagram:
  \begin{equation}
    \label{eq: commutative diagram for contact handles}
    \vcenter{
      \xymatrix{
        (K_0,\Ga_m)\otimes(S^3(1),\delta)\ar[dd]^{C_{h^1}}_\isom&(K_0,\Ga_m)\otimes(K_1,\Ga_{n})\ar[l]_{\Id\tensor 
          C_{h^2_1}}\ar[r]^{C_{h^2_0}\tensor\Id}\ar[d]^{C_{h^1}}_\isom&(K_1,\Ga_n)\otimes(S^3(1),\delta)\ar[dd]^{C_{h^1}}_\isom\\
        &(M_0,\ga_0)\ar[ld]_{C_{h^2_1}}\ar[rd]^{C_{h^2_0}}\ar[d]^{C_{h^2_2}}&\\
        (K_0,\Ga_m)\ar[rd]^{C_{h^2_0}}&(K_0\connsum 
        K_1,\Ga_{m+n})\ar[d]^{C^{\connsum}}&(K_1,\Ga_n)\ar[ld]_{C_{h^2_1}}\\
        &(S^3(1),\delta).
      }
    }
  \end{equation}
  Here and below, for the sake of space, we often denote $\SHG (-S^3(1), 
  -\delta)$ by $(S^3(1), \delta)$, denote $\SHG (-S^3(K), -\Gamma)$ by $(K, 
  \Gamma)$, and denote $\SHG (-M_0, -\ga_0)$ by $(M_0, \ga_0)$:

  Since $m$ and $n$ are chosen to be odd and sufficiently large, by 
  \cite[Proposition~5.8]{li2019direct},
  elements in $\khg(-S^3,K_0)$ and $\khg(-S^3,K_1)$
  of sufficiently large gradings
  can be
  found in
  $\shg(-S^3(K_0),-\Ga_m)$ and $\shg(-S^3(K_1),-\Ga_n)$
  respectively,
  as in the previous section.
  In particular, let $x_0\in \shg(-S^3(K_0),-\Ga_m)$ be an element representing 
  a non--$U$-torsion element in $\KHGm (-S^3, K_0)$ of maximal grading; then, 
  by \fullref{prop: another definition of tau},
  \[
    \gr_{S_0}(x_0)=\tau_G(K_0)+\frac{m-1}{2}, \qquad C_{h^2_0}(x_0)\neq0,
  \]
  where $\gr_{S_0}$ means the grading with respect to $S_0$, and the term 
  $(m-1)/2$
  represents the grading shift in the definition of $\KHGm$.
  Similarly, we can pick $y_0\in \shg(-S^3(K_1),-\Ga_n)$ to represent a 
  non--$U$-torsion element in $\KHGm (-S^3, K_1)$ of maximal grading; then
  \[
    \gr_{S_1}(y_0)=\tau_G(K_1)+\frac{n-1}{2}, \qquad C_{h^2_1}(y_0)\neq0.
  \]
  Let
  \[
    z_0 = C_{h^2_2}\circ C_{h^1}(x_0\otimes y_0) \in \SHG (-S^3 (K_0 \connsum 
    K_1), -\Ga_{m+n});
  \]
  then we know from \fullref{lem: grading preserving} that
  \[
    \gr_{S^-}(z_0)=\tau_G(K_0)+\tau_G(K_1)+\frac{m+n}{2}.
  \]
  From the commutative diagram \eqref{eq: commutative diagram for contact 
    handles}, we know that
  \begin{equation}
    \label{eq: z_0 is in the infinite U tower}
    C^{\sharp}(z_0) = C^{\sharp} \circ C_{h^2_2} \circ C_{h^1} (x_0 \tensor 
    y_0) = C_{h^2_0} \circ C_{h^1} \circ (\Id \tensor C_{h^2_1}) (x_0 \tensor 
    y_0) = C_{h^2_0} (x_0) \neq 0,
  \end{equation}
  where the third equality uses the fact that $C_{h^2_1 (y_0)} \neq 0$. Hence, 
  by \fullref{prop: another definition of tau},
  we have
  \[
    \tau_G(K_0\connsum K_1)+\frac{m+n}{2}\geq 
    gr_{S^-}(z_0)=\tau_G(K_0)+\tau_G(K_1)+\frac{m+n}{2},
  \]
  from which the proposition follows.
\end{proof}

We now upgrade the inequality in \fullref{prop: sup additivity} to an equality:

\begin{proof}[Proof of {\fullref{prop:tau-conn-sum}}]
  We keep all notation from the proof of \fullref{prop: sup additivity}. In 
  particular, we have an element
  \[
    z_0 = C_{h_2^2} \circ C_{h^1} (x_0 \tensor y_0) \in \SHG (-S^3 (K_0 \connsum 
    K_1), -\Ga_{m+n}),
  \]
  where $x_0 \in \SHG (-S^3 (K_0), -\Ga_m)$ and $y_0 \in \SHG (-S^3 (K_1), 
  -\Ga_n)$ represent non--$U$-torsion elements in $\KHGm (-S^3, K_0$ and $\KHGm 
  (-S^3, K_1)$ of maximal gradings respectively.

  By \eqref{eq: z_0 is in the infinite U tower}, we see that $z_0$ in fact 
  corresponds to a non--$U$-torsion element
  \[
    z_0^- \in \KHGm (-S^3, K_0 \connsum K_1).
  \]
  If we assume the contrary of the proposition, i.e.\
  \begin{equation*}
    \tau_G(K_0\connsum K_1)>\tau_G(K_0)+\tau_G(K_1),
  \end{equation*}
  then we are assuming that
  $z_0^-$ is not the starting point of the unique infinite $U$-tower;
  in other words, it
  has a pre-image under $U$. Translating back to
  $\SHG (-S^3(K_0\connsum K_1), -\Ga_{m+n})$, this means that there is an 
  element
  \[
    z_{1}\in\shg(-S^3(K_0\connsum K_1),-\Ga_{m+n-1})
  \]
  such that
  \begin{equation*}
    \psi^{m+n-1}_{+,m+n}(z_{1})=z_0.
  \end{equation*}
  By the positive bypass exact triangle in \eqref{eq:bypass}, we see that 
  \[
    \psi^{m+n}_{+,\infty} (z_0) = 0.
  \]
  We claim that this will lead to a contradiction.

  Indeed,
  consider the maps
  \begin{align*}
    \psi^m_{+,\mu} \colon \SHG (-S^3 (K_0), -\Ga_m) \to \SHG (-S^3 (K_0), 
    -\Ga_\mu),\\
    \psi^n_{+,\mu} \colon \SHG (-S^3 (K_1), -\Ga_n) \to \SHG (-S^3 (K_1), 
    -\Ga_\mu),
  \end{align*}
  which each fit into the positive bypass exact triangle in \eqref{eq:bypass}.
  Let $\be_0 \in \bdy S^3 (K_0)$ and $\be_1 \in \bdy S^3 (K_1)$ be the arcs 
  along which bypasses corresponding to these maps are attached; we may view 
  $\be_0$ and $\be_1$ in $(S^3 (K_0 \connsum K_1), \Ga_{m+n})$, as in 
  \fullref{fig: the by-passes}.
  (Since $\be_0$ and $\be_1$ are both disjoint from the $1$-handle $h^1$ and 
  the $2$-handle $h^2_2$, they survive in $(S^3(K_0\connsum K_1),\Ga_{m+n})$.)

  \begin{figure}[htbp]
    \captionsetup{aboveskip={\dimexpr10pt+2pt+\sactualfontsize\relax}}
    \labellist
    \small\hair 2pt
    \pinlabel {$S^3 (K_0)$} [t] at 59 0
    \pinlabel {$S^3 (K_1)$} [t] at 362 0
    \pinlabel {$\beta_0$} [b] at 66 268
    \pinlabel {$\beta_1$} [t] at 362 266
    \pinlabel {$\alpha$} [t] at 362 167
    \pinlabel {$h^1$} [t] at 212 136
    \endlabellist
    \includegraphics[width=0.5\textwidth]{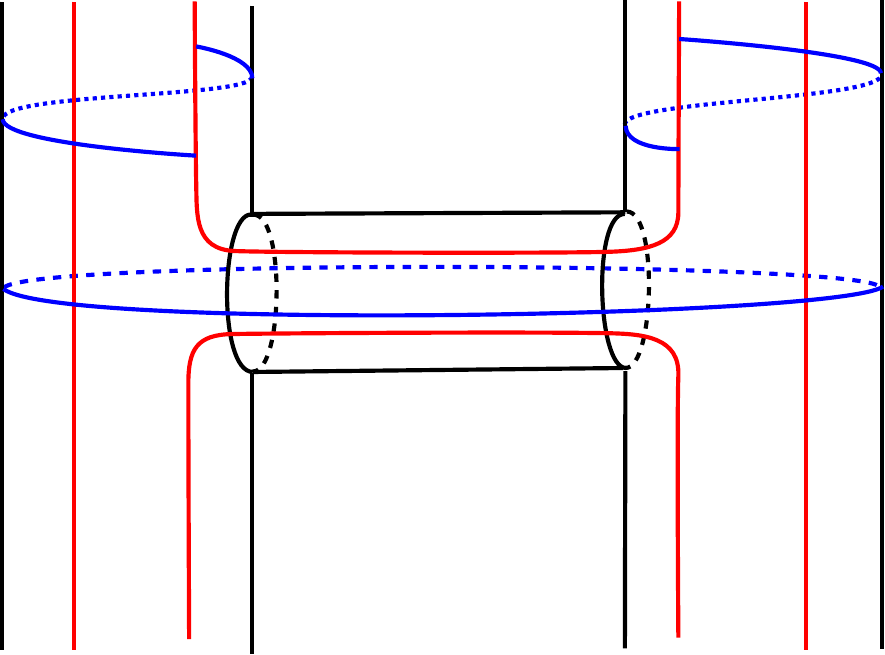}
    \caption{The arcs $\be_0$ and $\be_1$, along which bypasses are attached, 
      viewed in $(S^3 (K_0 \connsum K_1), \Ga_{m+n})$.}
    \label{fig: the by-passes}
  \end{figure}

  Inside $(S^3(K_0\connsum K_1),\Ga_{m+n})$, the arcs $\be_1$ and $\be_2$ are 
  isotopic; thus, they both correspond to the bypass map
  \[
    \psi^{m+n}_{+,\mu} \colon \SHG (-S^3 (K_0 \connsum K_1), -\Ga_{m+n}) \to 
    \SHG (-S^3 (K_0 \connsum K_1), -\Ga_\mu).
  \]
  For concreteness,
  suppose
  this bypass map
  is constructed via a bypass attached along $\be_2$.  Since $\be_2$ is 
  disjoint from $h^1$ and $h^2_2$, there is a commutative diagram as follows:
  \begin{equation*}
    \xymatrix{
      (K_0,\Ga_m)\otimes(K_1,\Ga_n)\ar[rr]^{\Id\otimes\psi^{n}_{+,\mu}}\ar[d]^{C_{h^2_2}\circ 
        C_{h^1}}&&(K_0,\Ga_m)\otimes(K_1,\Ga_{\mu})\ar[d]^{C_{h^2_2}\circ 
        C_{h^1}}\\
      (K_0\connsum K_1,\Ga_{m+n})\ar[rr]^{\psi^{m+n}_{+,\mu}}&&(K_0\connsum 
      K_1,\Ga_{\mu})
    }
  \end{equation*}
  (Here, we are using the simplified notation as in \eqref{eq: commutative 
    diagram for contact handles}.)

  From the commutativity, we know that
  \begin{equation*}
    \label{eq: a second form of psi(z_0)}
    \begin{aligned}
      \psi^{m+n}_{+,\mu}(z_0) & = \psi^{m+n}_{+,\mu}\circ C_{h^2_2}\circ 
      C_{h^1}(x_0\otimes y_0)\\
      & = C_{h^2_2}\circ C_{h^1}\circ(\Id\otimes\psi^{n}_{+,\mu})(x_0\otimes 
      y_0)\\
      & = C_{h^2_2}\circ C_{h^1}(x_0\otimes y_{\mu}),
    \end{aligned}
  \end{equation*}
  where $y_\mu = \psi^n_{+,\mu} (y_0)$.
  Since $y_0$ corresponds to a non--$U$-torsion element in $\KHGm (-S^3, K_1)$ 
  of maximal grading,
  we know that $y_0\notin \Im \psi^{n-1}_{+,n}$,
  and so by the exactness of \eqref{eq:bypass}, we know that $y_{\mu}\neq0$.  

  Now we claim that the following diagram commutes:
  \begin{equation}
    \label{eq: commutative diagram, be_1, h^1 and h^2_2}
    \vcenter{
      \xymatrix{
        (K_0,\Ga_m)\otimes(K_1,\Ga_{\mu})\ar[rr]^{\psi^m_{+,\mu}\otimes 
          \Id}\ar[d]^{C_{h^2_2}\circ 
          C_{h^1}}&&(K_0,\Ga_{\mu})\otimes(K_1,\Ga_{\mu})\ar[d]^{\isom}\\
        (K_0\connsum K_1,\Ga_{\mu})\ar[rr]^{=}&&(K_0\connsum K_1,\Ga_{\mu}).
      }
    }
  \end{equation}
  (The isomorphism in the right column arises from the fact that the two 
  sutured manifolds have the same closure; the same is true for $C_{h^1}$ on 
  the left column, but we display it explicitly so that $C_{h^2_2}$ makes 
  sense.) Since $x_\mu = \psi^m_{+,\mu} (x_0) \neq 0$ (as $x_0$ represents a 
  non--$U$-torsion element in $\KHGm (-S^3, K_0)$ of maximal grading), this 
  will show that
  \[
    \psi^{m+n}_{+,\mu} (z_0) = C_{h^2_2} \circ C_{h^1} (x_0 \otimes y_\mu) = 
    \psi^m_{+,\mu} (x_0) \tensor y_\mu = x_\mu \tensor y_\mu \neq 0,
  \]
  giving us the desired contradiction.

  The rest of the proof is devoted to proving the commutativity of \eqref{eq: 
    commutative diagram, be_1, h^1 and h^2_2}.
  Let $(M_1,\ga_1)$ be the result of attaching the handle $h^1$ to 
  $(S^3(K_0),\Ga_{m})\sqcup (S^3(K_1),\Ga_{\mu})$.
  (This gives us a map
  \[
    C_{h^1} \colon \SHG (-S^3 (K_0), -\Ga_{m}) \tensor \SHG (-S^3 (K_1), 
    -\Ga_\mu) \to \SHG (-M_1, -\ga_1),
  \]
  similar to \eqref{eq: C h 1}, but with sutures $-\Ga_\mu$ instead of $-\Ga_n$ 
  on $-S^3 (K_1)$.)
  Our strategy is to analyze the contact $2$-handle attachment along $\alpha$, 
  corresponding to $C_{h^2_2}$ and viewed in $(-M_1, -\ga_1)$, and compare it 
  to the bypass attachment along $\beta_0$, corresponding to $\psi^m_{+,\mu}$.  
  A bypass attachment is in fact the composition of a contact $1$-handle and a 
  contact $2$-handle (see, for example, \cite[Section~5]{baldwin2016contact}); 
  in our context, we shall work with the preclosures of the sutured manifolds 
  (see \cite[Section~4.2]{baldwin2016contact} for details of the relevant 
  constructions), where the contact $1$-handle associated to $\psi^m_{+,\mu}$ 
  will be identified with a part of the auxiliary surface associated to $(M_1, 
  \gamma_1)$, and the attaching curve of the contact $2$-handle associated to 
  $\psi^m_{+,\mu}$ will be identified with an isotopic copy of $\alpha$.

  \begin{figure}[htbp]
    \captionsetup{aboveskip={\dimexpr10pt+2pt+\sactualfontsize\relax}}
    \labellist
    \small\hair 2pt
    \pinlabel {$S^3 (K_0)$} [t] at 57 0
    \pinlabel {$S^3 (K_1)$} [t] at 364 0
    \pinlabel {\textcolor{red}{\rule{0.3in}{0.5pt}}} [l] at 88 40
    \pinlabel {\textcolor{red}{$\ga_1 \setminus \hat{\ga}_1$}} [l] at 132 40
    \pinlabel {\textcolor{red}{\rule{0.3in}{0.5pt}}} [l] at 89 280
    \pinlabel {\textcolor{red}{$\ga_1 \setminus \hat{\ga}_1$}} [l] at 133 280
    \pinlabel {\textcolor{red}{$\ga_1 \setminus \hat{\ga}_1$}} [t] at 362 40
    \pinlabel {\textcolor{red}{$\hat{\ga}_1$}} [t] at 362 117
    \pinlabel {\textcolor{blue}{$\al \setminus \hat{\al}$}} [t] at 62 164
    \pinlabel {\textcolor{blue}{$\hat{\al}$}} [t] at 362 164
    \endlabellist
    \includegraphics[width=0.5\textwidth]{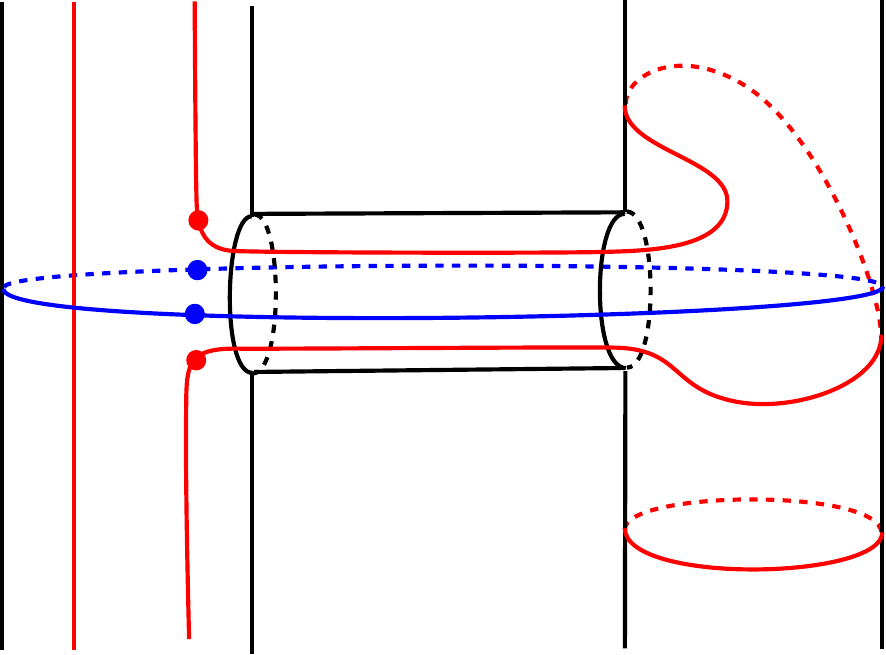}
    \caption{The arcs $\hat{\ga}_1$ and $\hat{\al}$, which we think of as 
      isotopic relative to their endpoints. Their endpoints are denoted by the 
      red and blue dots respectively. The short arcs $\zeta$ are omitted.}
    \label{fig: meridional suture}
  \end{figure}

  See \fullref{fig: meridional suture}.
  Because we have the sutures $\Ga_\mu$ on $\bdy S^3 (K_1)$, we see that after 
  the $1$-handle $h^1$ is added, one component of $\ga_1$ is simply a meridian 
  on the $\bdy S^3 (K_1)$ part of the boundary (of $S^3 (K_0 \connsum K_1)$), 
  while the other component, which intersects the $\bdy S^3 (K_0)$ part of the 
  boundary, also wraps around the $\bdy S^3 (K_1)$ part of the boundary like a 
  meridian. We may thus view a part of this latter component, an arc 
  $\hat{\ga}_1$, as isotopic to a part of $\al$, which we call $\hat{\al}$, 
  relative to their endpoints.
  More precisely, while the arcs $\hat{\ga}_1$ and $\hat{\al}$ do not have the 
  same endpoints; however, from \fullref{fig: meridional suture}, one can pair 
  up the endpoints obviously by short arcs $\zeta$.

  Suppose $T_1$ is a connected auxiliary surface of $(M_1,\ga_1)$; then we can 
  form the preclosure
  \[
    \widetilde{M}=M_1\cup [-1,1]\times T_1.
  \]
  From \cite[Section~4.2.2]{baldwin2016contact}, there is an auxiliary surface 
  $T$ for $(S^3(K_0),\Ga_{m})\sqcup (S^3(K_1),\Ga_{\mu})$,
  obtained from $T_1$ by attaching a $2$-dimensional $1$-handle 
  $\overline{h}^1$, which corresponds to the $3$-dimensional $1$-handle $h^1$, 
  as in \fullref{fig:aux-surf}, so that we also have
  \[
    \widetilde{M}=\paren{S^3(K_0)\sqcup S^3(K_1))} \cup \paren{[-1,1]\times T}.
  \]
  In this description, we can think of $h^1$ as a thickening of $\overline{h}^1 
  \subset T$.

  \begin{figure}[htbp]
    \captionsetup{aboveskip={\dimexpr10pt+2pt+\sactualfontsize\relax}}
    \labellist
    \small\hair 2pt
    \pinlabel \textcolor{red}{{$\hat{\ga}_1$}} [t] at 103 209
    \pinlabel \textcolor{red}{{$\hat{\ga}_1$}} [t] at 347 209
    \pinlabel \textcolor{blue}{{\rule{0.3in}{0.5pt}}} [r] at 322 295
    \pinlabel \textcolor{blue}{{$\hat{\al} \cap \overline{h}^1$}} [r] at 273 295
    \pinlabel \textcolor{blue}{{\rule{0.3in}{0.5pt}}} [l] at 359 295
    \pinlabel \textcolor{blue}{{$\hat{\al} \cap \overline{h}^1$}} [l] at 406 295
    \pinlabel {$\overline{h}^1$} [ ] at 302 251
    \pinlabel \textcolor{blue}{{$\hat{\al}_T$}} [t] at 574 150
    \pinlabel {$T_1$} [t] at 104 0
    \pinlabel {$T$} [t] at 347 0
    \pinlabel {$T$} [t] at 572 0
    \pinlabel {$T_2$} [t] at 819 0
    \endlabellist
    \includegraphics[width=\textwidth]{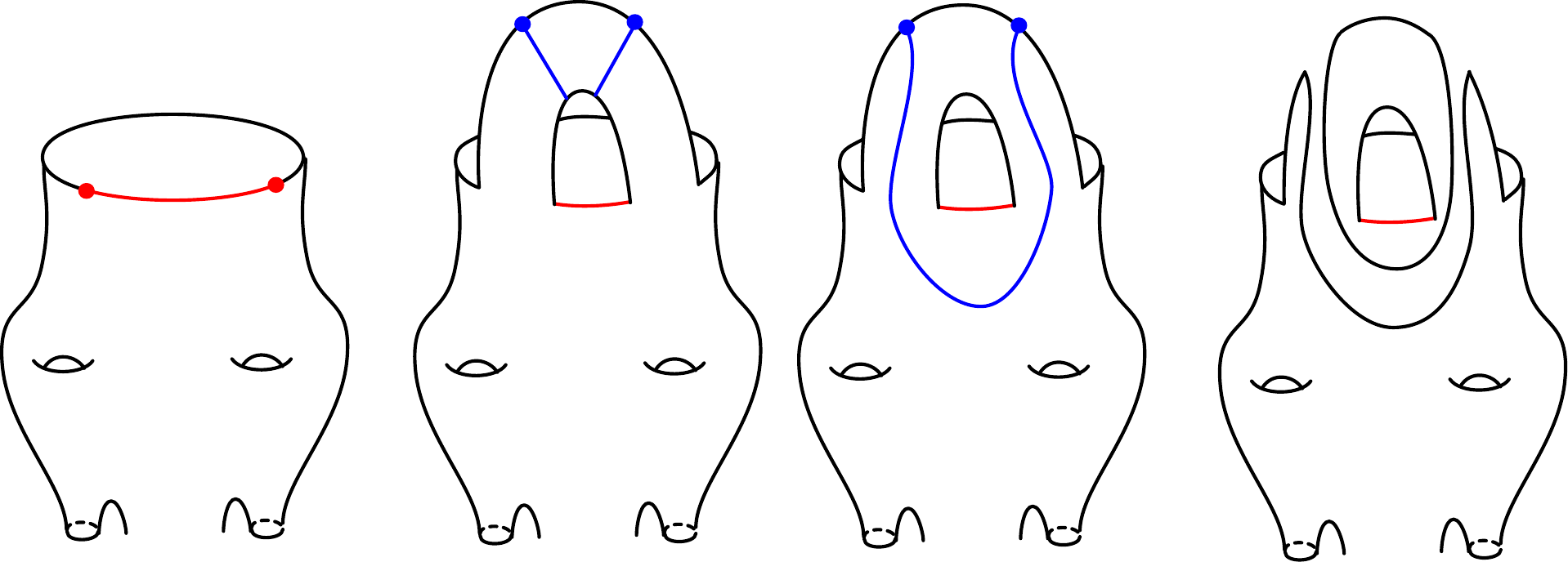}
    \caption{Constructing an auxiliary surface $T_2$ for $(M_2, \ga_2)$, from 
      an auxiliary surface $T_1$ for $(M_1, \ga_1)$. In the second diagram, 
      only a part of $\hat{\al}$ is on $T$; we isotope all of $\hat{\al}$ onto 
      $T$ in the third diagram. As shown, the auxiliary surfaces may have 
      non-zero genus; their irrelevant boundary componentsare omitted.}
    \label{fig:aux-surf}
  \end{figure}

  From \cite[Section~4.2.3]{baldwin2016contact}, attaching the contact 
  $2$-handle $h^2_2$ along $\al$ corresponds to performing a $0$-surgery along 
  a push off of $\al$ on the level of preclosures. Since $\hat{\al}$ is 
  isotopic to $\hat{\ga}_1$ relative to their endpoints, we can isotope 
  $\hat{\al}$ onto $\hat{\ga}_1$ (using $\zeta$) and hence onto $T$,
  to give a properly embedded arc $\hat{\al}_T \subset T$,
  as depicted in \fullref{fig:aux-surf}. The product neighborhood of 
  $\hat{\al}_T$
  corresponds to a contact $1$-handle $h^1_0$ attached to 
  $(S^3(K_0),\Ga_{m})\sqcup (S^3(K_1),\Ga_{\mu})$; this is the contact 
  $1$-handle associated to $\psi^m_{+,\mu}$.
  Let $(M_2,\ga_2)$ be the balanced sutured manifold obtained by attaching 
  $h^1_0$, and let $T_2 = T \setminus \hat{\al}$; then $T_2$ is an auxiliary 
  surface for $(M_2,\ga_2)$, and thus
  \[
    \widetilde{M}=(M_2,\ga_2)\cup[-1,1]\times T_2.
  \]
  Because $\hat{\al}_T$ is isotopic to $\hat{\al}$, we can think of $h^1_0$ as 
  attached to
  $(S^3(K_0),\Ga_{m})\sqcup (S^3(K_1),\Ga_{\mu})$
  along the two end points of $\al\setminus \hat{\al}$ on $\Ga_m \subset \bdy 
  S^3 (K_0)$.  We further isotope $\hat{\al}_T$ to an arc $\hat{\al}_{h^1_0}$, 
  on the boundary of $h^1_0$, that intersects $\ga_2$ exactly once, and let
  \[
    \al'=(\al\setminus\hat{\al})\cup\hat{\al}_{h^1_0} \subset \partial{M}_2.
  \]
  Then the $0$-surgery (with respect to the surface framing) along a push off 
  of $\al$ corresponds to a $0$-surgery along a push off of $\al'$, and hence 
  to a contact $2$-handle attachment along $\al'$. The $1$-handle $h^1_0$ and 
  the $2$-handle attached along $\al'$ together correspond to a bypass attached 
  along $\al\setminus \hat{\al}$.

  Now under the same identification of the endpoints as before---by the short 
  arcs $\zeta$---%
  we see that $\al \setminus \hat{\al}$
  is isotopic to the arc $\be_0$ relative to their endpoints (if we allow the 
  endpoints to move along $\Gamma_m$), viewed on $(S^3 (K_0), \Ga_m)$; compare 
  \fullref{fig: the by-passes} and \fullref{fig: meridional suture}.  (They are 
  not isotopic when viewed on $(S^3 (K_0 \connsum K_1), \Ga_{m+n})$.)
  Thus, we see that the map $C_{h^2_2}$ corresponds to the map associated to a 
  bypass attached along $\beta_0$, which is $\psi^m_{+,\mu}$, and
  the proposition follows.
\end{proof}

Having achieved our goal of the section, we end it by spelling out an immediate 
corollary:

\begin{corollary}
  \label{cor: tau_G of the mirror image}
  For all knots $K\subset S^3$, we have $\tau_G(\mir{K})=-\tau_G(K)$,
  where $\mir{K}$ is the mirror image of $K$.
\end{corollary}

\begin{proof}
  This is a direct consequence of \fullref{prop:tau-conc} and 
  \fullref{prop:tau-conn-sum}.
\end{proof}

%% file: sec-identify.tex
\section{Identifying the tau invariants}
\label{sec:identify}

In this section, we identify the invariants $\tauG$ and $\taushG$, proving 
\fullref{thm:main}. While the instanton and monopole Floer theories are 
formally similar, there are some differences in their definitions. For example, 
the definition of $\SHM (M, \gamma)$ involves a decomposition into $\Spinc$ 
structures of $Y$ (where $(Y, R)$ is a closure of $(M, \gamma)$), which are in 
bijection with $H^2 (Y)$ (see \cite{kronheimer2007monopoles}); the definition 
of $\SHI (M, \gamma)$ involves a generalized eigenspace decomposition by 
actions of surfaces, corresponding to $H^2 (Y) / \Tors$ (see, for example, 
\cite[Corollary~7.6]{kronheimer2010knots}).  To identify $\tauG$ with 
$\taushG$, we have to work directly with these objects above. As $\taushG$ is 
defined only in the instanton setting in \cite{baldwin2020framed}, we focus on 
$\tauI$ and $\taushI$ throughout the section, and discuss the changes necessary 
for the monopole setting.

\subsection{\texorpdfstring{$\tau_I$}{The instanton tau invariant} revisited}
Suppose $K\subset S^3$ is an oriented knot. Let $\lambda$ and $\mu$ be the 
longitude and meridian on $\partial S^3(K)$, respectively.

Let $\Ga_{n}$ be the suture on $\partial S^3(K)$ consisting of two curves of 
slope $-n$, or equivalently, of class $\pm(\lambda-n\mu)$. Let $S$ be a minimal 
genus Seifert surface of $K$. For any $n\in\intg$, take $S_n$ to be an isotopy of $S$ that intersects the suture $\Gamma_n$ exactly $2n$ times. Let $S_n^{\tau(n)}$ be either a negative stabilization of $S_n$ when $n$ is even, or the surface $S_n$ itself when $n$ is odd. The surface $S_n^{\tau(n)}$ induces a $\intg$-grading on 
$\shi(-S^3(K),-\Ga_n)$ for any $n\in\intg$. We write this grading as
\[
  \shi(-S^3(K),-\Ga_n,S_n^{\tau(n)},i).
\]

Recall $\mu$ is a meridian of $K$ on $\partial S^3(K)$. Let $\al\subset{\rm 
  int}(S^3(K))$ be a push off of $\mu$ into the interior of $(S^3(K),\Ga_n)$.  
We can perform $0$-surgery along $\al$ and obtain a map
\[
  F_{\al,n}\colon \shi(-S^3(K),-\Ga_n)\ra\shi(-M,-\Ga_n),
\]
where $(M,\ga)$ is obtained from $(S^3(K),\Ga_n)$ by performing a $0$-surgery 
along $\al$. (The map $F_{\al,n}$ here is in fact the same as the map $C_{h,n}$ 
in \fullref{sec:tau-conc}: While the former is associated to a $0$-surgery and 
the latter a contact $2$-handle, they are in fact the same by construction; see 
\cite[Section~4.2.3]{baldwin2016contact}.) Note $\al$ is a meridian of $K$, and 
so it is an unknot inside $S^3$.  Hence the result of performing $0$-surgery 
along $\al\subset S^3$ is $S^1\times S^2$.  The $2$-sphere $S^2$ comes from 
capping off the Seifert disk bounded by $\alpha$ inside $S^3$. Since $K$ 
intersects the Seifert disk and hence $S^2\subset S^1\times S^2$ at one point, 
by the 3-dimensional light bulb theorem, $K$ is isotopic to $S^1\subset 
S^1\times S^2$. Hence the complement $M=S^1\times S^2(K)$ is simply a solid 
torus. The meridian disk $D$ of this solid torus is $S^2\cap S^1\times S^2(K)$.  
Also we have $\partial D=\mu$ on $\partial (S^1\times S^2(K))x=\partial 
S^3(K)$, so it intersects the suture $\Ga_{n}$ twice.  Hence we know that
\[
  \shi(-M,-\Ga_n)\cong\mathbb{C}.
\]

Let $-K$ be the same knot as $K$ but with reversed orientation. The longitude 
and meridian for $-K$ is minus those of $K$. Hence the sutured manifold 
associated $-K$ is $(-S^3(K),\Ga_n)$, as opposed to $(-S^3 (K), -\Ga_n)$, which 
is associated to $K$. Also, $-S_n^{\tau(n)}$ is a Seifert surface of $-K$.  
Hence $-S_n^{\tau(n)}$ induces a grading
\[
  \shi(-S^3(K),\Ga_n,-S_n^{\tau(n)},i).
\]

\begin{lemma}
  \label{lem: iso between K and -K}
For any $i\in\intg$, we have an isomorphism
\[
  \shi(-S^3(K),\Ga_n,-S_n^{\tau(n)},i)\cong\shi(-S^3(K),-\Ga_n,S_n^{\tau(n)},i).
\]
\end{lemma}

\begin{proof}
Let us first ignore the surfaces $S_n^{\tau(n)}$ and $-S_n^{\tau(n)}$ and the gradings they induce. 
Suppose $(Y,R)$ is a closure of $(-S^3(K),-\Ga_n)$, then $(Y,-R)$ is a closure 
of $(-S^3(K),\Ga_n)$. Let
\[
  \Eig(\mu(R),i)=i\text{-generalized eigenspace of } \mu(R).
\]
Then
\[
  \shi(-S^3(K),-\Ga_n)=\Eig(\mu(R),2g(R)-2),
\]
and
\[
  \shi(-S^3(K),\Ga_n)=\Eig(-\mu(R),2g(R)-2)=\Eig(\mu(R),2-2g(R))
\]
However, there is a symmetry
\[
  \Eig(\mu(R),2g(R)-2)\cong \Eig(\mu(R),2-2g(R)),
\]
which is an analog of conjugating the $\Spinc$ structure in monopole and 
Heegaard Floer theory; see, for example, \cite[Lemma~2.3]{baldwin2019lspace}.  
This symmetry also commutes with the cobordism maps.  Hence we know that
\[
  \shi(-S^3(K),\Ga_n)\cong \shi(-S^3(K),-\Ga_n).
\]
If we take the grading into consideration, then we can pick a closure $(Y,R)$ 
of $(-S^3(K),-\Ga_n)$ so that $S_n^{\tau(n)}$ extends to a closed surface $\overline{S}_n$. Then 
$(Y,-R)$ is a closure of $(-S^3(K),\Ga_n)$ so that $-S_n^{\tau(n)}$ extends to $-\overline{S}_n$. 
Thus
\begin{align*}
\shi(-S^3(K),\Ga_n,-S_n^{\tau(n)},i)&=\Eig(\mu(R),2-2g(R))\cap\Eig(\mu(\overline{S}_n),-2i)\\
&\cong\Eig(\mu(R),2g(R)-2)\cap\Eig(\mu(\overline{S}_n),2i)\\
&=\shi(-S^3(K),-\Ga_n,S_n^{\tau(n)},i).
\end{align*}
Again this isomorphism commutes with cobordism maps.
\end{proof}

\begin{corollary}
$\tau_I(-K)=\tau_I(K)$.
\end{corollary}

\begin{proof}
Straightforward.
\end{proof}

\begin{corollary}
  \label{cor: minimal non vanishing grading}
For $n$ odd and large enough,
$$\min\setc{i}{\exists \, x\in\shi(-S^3(K),-\Ga_n,S_n,i) \text{ s.t.\ } 
  F_{\al,n}(x)\neq0}=-\tauI(K)-\frac{n-1}{2}.$$
\end{corollary}

\begin{proof}
Since $\partial S^3(K)$ is an annulus, we can isotope $\Ga_n$ to be $-\Ga_n$. 
Hence there is a diffeomorphism
$$f\colon(S^3(K),\Ga_n)\ra (S^3(K),-\Ga_n)$$
which restricts to the identity outside a collar of the boundary. Hence, under 
this diffeomorphism, the surface $S_n$ (note we have assumed that $n$ is odd, so $S_n^{\tau(n)}=S_n$) is preserved:
$$f\colon(S^3(K),\Ga_n,S_n)\ra (S^3(K),-\Ga_n,S_n).$$
Thus this diffeomorphism induces an isomorphism
$$\shi(-S^3(K),-\Ga_n,S_n,i)\cong \shi(-S^3(K),\Ga_n,S_n,i).$$
Hence combined with \fullref{lem: iso between K and -K}, we know that
\begin{align*}
\shi(-S^3(K),-\Ga_n,S_n,i)&\cong\shi(-S^3(K),\Ga_n,S_n,i)\\
&=\shi(-S^3(K),\Ga_n,-S_n,-i)\\
&\cong\shi(-S^3(K),-\Ga_n,S_n,-i).
\end{align*}
Since the diffeomorphism $f$ restricts to the identity outside a collar of 
$\partial S^3(K)$, we can take $\al$ to be inside the region where $f$ is the 
identity. Hence the above isomorphisms commute with the maps associated to the 
zero surgery along $\al$. Thus the above symmetric implies that the maximum of 
the grading restricted to where $F_{\al,n}$ is non-trivial is minus the minimum 
of such gradings.
\end{proof}

\subsection{The sutured manifold \texorpdfstring{$(S^3(K),\Ga_n)$}{associated 
    to the knot complement}}
Let
$$d_n={\rm dim}_{\mathbb{C}}\shi(-S^3,-\Ga_n).$$
We have the following.

\begin{lemma}
  \label{lem: increasing and decreasing of d_n}
If $d_n>d_{n-1}$, then $d_{n+1}>d_n$. If $d_n<d_{n-1}$, then $d_{n-1}<d_{n-2}$.
\end{lemma}

\begin{proof}
There are two bypass exact triangles:
\begin{equation*}
\xymatrix{
\shi(-S^3(K),-\Ga_n)\ar[rr]^{\psi^n_{\pm,n+1}}&&\shi(-S^3(K),-\Ga_{n+1})\ar[dl]\\
&\shi(-S^3(K),-\Ga_{\mu})\ar[ul]&
}    \end{equation*}
Here $\Ga_{\mu}$ consists of two meridians. We can also perform $0$-surgeries 
along $\al\subset{\rm int}(S^3(K))$ and obtain a commutative diagram:
\begin{equation*}
\xymatrix{
\shi(-S^3(K),-\Ga_n)\ar[rr]^{\psi_{\pm,n+1}^n}\ar[dd]^{F_{\al,n}}&&\shi(-S^3(K),-\Ga_{n+1})\ar[dl]\ar[dd]^{F_{\al,n+1}}\\
&\shi(-S^3(K),-\Ga_{\mu})\ar[ul]\ar[dd]&\\
\shi(-M,-\Ga_n)\ar[rr]&&\shi(-M,-\Ga_{n+1})\ar[dl]\\
&\shi(-M,-\Ga_{\mu})\ar[ul]&
}    \end{equation*}
Here, as above, $M=S^1\times S^2(K)$ is a solid torus with meridian disk $D$ 
and $\partial D=\mu\subset \partial S^1\times S^2(K)=S^3(K)$. Hence we know 
that $(-M,-\Ga_{\mu})$ is not taut since $D$ is a boundary compressing disk. 
Then we know that
$$\shi(-M,-\Ga_{\mu})=0,$$
and we have a commutative diagram
\begin{equation*}
\xymatrix{
\shi(-S^3(K),-\Ga_n)\ar[rr]^{\psi_{\pm,n+1}^n}\ar[dd]^{F_{\al,n}}&&\shi(-S^3(K),-\Ga_{n+1})\ar[dd]^{F_{\al,n+1}}\\
&&\\
\shi(-M,-\Ga_n)\ar[rr]^{\cong}&&\shi(-M,-\Ga_{n+1})\\
}
\end{equation*}
Thus, we know that $F_{\al,n+1}\neq 0$ if $F_{\al,n}\neq 0$.  There is a 
surgery exact triangle associated to $\al$:
\begin{equation*}
\xymatrix{
\shi(-S^3(K),-\Ga_n)\ar[rr]^{F_{\al,n}}&&\shi(-M,-\Ga_{n})\ar[dl]^{G_{\al,n-1}}\\
&\shi(-S^3(K),-\Ga_{n-1})\ar[ul]&
}
\end{equation*}
Since we have computed that
$$\shi(-M,-\Ga_n)\cong\mathbb{C},$$
we know that
$$d_{n}>d_{n-1}$$
if and only if $F_{\al_{n}}\neq0$ if and only if $G_{\al,n-1}=0$. Thus, 
$d_{n}>d_{n-1}$ implies that $F_{\al,n}\neq0$ and hence $F_{\al,n+1}\neq 0$ and 
finally $d_{n+1}>d_n$.

Similarly, we have a commutative diagram
\begin{equation*}
\xymatrix{
\shi(-S^3(K),-\Ga_n)\ar[rr]^{\psi_{\pm,n+1}^n}&&\shi(-S^3(K),-\Ga_{n+1})\\
&&\\
\shi(-M,-\Ga_{n+1})\ar[uu]^{G_{\al,n}}\ar[rr]^{\cong}&&\shi(-M,-\Ga_{n+2})\ar[uu]^{G_{\al,n+1}}\\
}
\end{equation*}
and know that $G_{\al,n}\neq0$ if $G_{\al,n+1}\neq0$. Hence $d_{n}<d_{n-1}$ 
implies that $G_{\al,n-1}\neq0$ and hence $G_{\al,n-2}\neq0$ and finally 
$d_{n-1}<d_{n-2}$.
\end{proof}

\begin{corollary}
The sequence
$$\{d_n={\rm dim}_{\mathbb{C}}\shi(-S^3(K),-\Ga_n)\}$$
has a unique minimum.
\end{corollary}

\begin{proof}
It follows directly from \fullref{lem: increasing and decreasing of d_n} and 
the fact that $d_n\geq0$.
\end{proof}

\begin{definition}
Let $n_0$ be the index of the unique minimus of the sequence $\{d_n\}$.
\end{definition}

Recall that we have a map
$$G_{\al,n}\colon \shi(-M,-\Ga_{n+1})\ra\shi(-S^3(K),-\Ga_n)$$
and
$$\shi(-M,-\Ga_{n+1})\cong \mathbb{C}.$$
We can define $$x_n=G_{\al,n}(\mathbf{1}_{n+1})$$
where $\mathbf{1}_{n+1}\in\shi(-M,-\Ga_{n+1})$ is a generator. From the 
commutative diagram in \fullref{lem: increasing and decreasing of d_n}, there 
are non-zero complex numbers $z_n$ so that
$$x_{n+1}=z_{n}\cdot\psi_{\pm,n+1}^n(x_n).$$
We choose suitable generators $\mathbf{1}_{n+1}\in\shi(-M,-\Ga_{n+1})$ so that 
all $z_n=1$ and simply having
\begin{equation}\label{eq: relate x_n with x_n+1}
x_{n+1}=\psi_{\pm,n+1}^n(x_n).
\end{equation}
Let
$$x_{n}=\sum_{i\in\intg}x_{n,i},$$
where $x_{n,i}\in \shi(-S^3(K),-\Ga_n,S^{\tau(n)}_n,i)$. 
\begin{definition}
Define
$$l_n=\max\setc{i}{x_{n,i}\neq0}-\min\setc{i}{x_{n,i}\neq0}.$$
\end{definition}

\begin{lemma}\label{lem: a bound for l_-n}
For $n$ odd and large enough, we know that
$$l_{-n}=2\tau_I(\mir{K})+n.$$
Here $\mir{K}$ is the mirror of $K$.
\end{lemma}
\begin{proof}
There is an orientation preserving diffeomorphism
$$(-S^3(K),-\Ga_{-n},S_n)\cong(S^3(\mir{K}),\Ga_{n},S_n).$$
Here $S_n$ on the right hand side is a Seifert surface of $\mir{K}$. Hence we 
have a commutative diagram
\begin{equation*}
\xymatrix{
\shi(-S^3(K),-\Ga_{-n},S_n,i)\ar[rr]^{\cong}&&\shi(S^3(\mir{K}),\Ga_{n},S_n,i)\\
&&\\
\shi(-M,-\Ga_{-n})\ar[rr]^{\cong}\ar[uu]^{G_{\al,n}}&&\shi(M,\Ga_{n})\ar[uu]^{G'_{\al,n}}\\
}
\end{equation*}
Note we know that $\mir{K}\subset S^1\times S^2$ is also isotopic to the knot 
$S^1$ so we keep writing $M$ on the bottom right of the above diagram. Note 
from \cite{li2018gluing}, we have a natural identification
$$\shi(S^3(\mir{K}),\Ga_{n},S,i)\cong 
\shi(-S^3(\mir{K}),\Ga_{n},S,i)^{\vee},{~\rm 
  and~}\shi(M,\Ga_n)\cong\shi(-M,\Ga_n)^{\vee}$$
where $V^{\vee}$ denotes the vector space dual to $V$. We have a commutative 
diagram
\begin{equation*}
\xymatrix{
\shi(S^3(\mir{K}),\Ga_{n},S,i)\ar[rr]^{\cong}&&\shi(-S^3(\mir{K}),\Ga_{n},S,i)^{\vee}\\
&&\\
\shi(M,\Ga_{n})\ar[rr]^{\cong}\ar[uu]^{G'_{\al,n}}&&\shi(-M,\Ga_{n})^{\vee}\ar[uu]^{F_{\al,n}^{\vee}}\\
}
\end{equation*}
Here, the map $F_{\al,n}^{\vee}$ is the dual of the map $$F_{\al,n}\colon 
\shi(-S^3(\mir{K}),\Ga_{n},S,i)\ra\shi(-M,\Ga_{n}).$$
The reason why $G'_{\al,n}$ is the dual of $F_{\al,n}$ is that they are 
induced by the cobordism but with reversed orientations. We 
know when $F_{\al,n}$ is zero due to \fullref{cor: minimal non vanishing 
  grading}.  Translating back to $x_{-n}$, we obtain
$$\max\setc{i}{x_{-n,i}\neq0}=\tau_I(\mir{K})+\frac{n-1}{2}$$
and
$$\min\setc{i}{x_{-n,i}\neq0}=-\tau_I(\mir{K})-\frac{n-1}{2}.$$
Hence we have
\[
  l_{-n}=2\tau_I(\mir{K})+n. \qedhere
\]
\end{proof}

\begin{lemma}\label{lem: decreasing of l_n}
For any $n\in \intg$, if $l_n>0$, then $l_{n+1}\leq l_{n}-1$.
\end{lemma}
\begin{proof}
Write $$x_{n}=\sum_{i\in\intg}x_{n,i}.$$
Suppose $x_{n,i}=0$ for $i>i_{max}$ and $i<i_{min}$, and $x_{n,i}\neq0$ for 
$i=i_{max}$ and $i=i_{min}$. Then
$$l_n=i_{max}-i_{min}+1.$$

Note according to the grading induced by $S^{\tau}$ ($S$ if $n$ is odd, and a 
negative stabilization of $S$ if $n$ is even), $\psi_{-,n+1}^n$ is grading 
preserving and $\psi_{+,n+1}^n$ shifts the grading down by one. Thus we can apply 
\eqref{eq: relate x_n with x_n+1} and for $i\geq i_{max}$,
$$x_{n+1,i}=\psi_{+,n+1}^n(x_{n,i+1})=0,$$
for $i<i_{min}$,
$$x_{n+1,i}=\psi_{-,n+1}^n(x_{n,i+1})=0.$$
Hence we conclude that
\[
  l_{n+1}\leq l_n-1. \qedhere
\]
\end{proof}

\begin{corollary}
  \label{cor: a bound for n_0}
We have the inequality
$$n_0\leq2\tauI(\mir{K}).$$
\end{corollary}

\begin{proof}
The notation $l_n$ tracks the difference between the top and bottom 
non-vanishing graded part of $x_n$. Hence as long as $x_n$ is non-zero, at 
least one of its graded parts must be non-zero and hence by definition $l_n>0$.  
Note there is an exact triangle between \begin{equation*}
\xymatrix{
 \shi(-S^3, -\Ga_n)\ar[rr]^{}&&\shi(-S^3(K),-\Gamma_{n+1})\ar[dl]\\
&\shi(-S^3(1),-\delta)\ar[ul]&
}    \end{equation*} The third term has dimension 1. So adopting our notations 
of $G_{\alpha,n}$ and $d_n$, there are only two possibilities:
\begin{enumerate}
  \item $d_{n}=d_{n+1}+1$, and it is equivalent to having $G_{\alpha,n}$ not 
    being zero; or
  \item $d_{n}=d_{n+1}-1$, and it is equivalent to having $G_{\alpha,n}$ being 
    zero.
\end{enumerate}

Hence $G_{\alpha,n}=0$ implies that $d_{n}=d_{n+1}-1$. Since $n_0$ is defined 
to the minimal integer k such that $d_{k}=d_{k+1}-1$ holds. So $n$ is at least 
$n_0$. Finally from \fullref{lem: a bound for l_-n} and \fullref{lem: 
  decreasing of l_n}, we have
\[
  l_{2\tauI(\mir{K})}=0.\qedhere
\]
\end{proof}

\subsection{Identifying \texorpdfstring{$\tau_I$ with $\taushI$}{tau with 
    tau-sharp}}

\begin{lemma}\label{lem: lower bound of nu^sharp}
For any knot $K$ so that $\nu^{\sharp}(K)\neq0$, we have
$$\nu^{\sharp}(K)>-2\tau_I(\mir{K})-2.$$
\end{lemma}

\begin{proof}
Take $n=2\tau_I(\mir{K})+2$. Then from \fullref{cor: a bound for n_0} we know 
that $d_{n-1}>d_{n-2}$ and hence the map $F_{\al,n-1}$ is surjective as in the 
proof of \fullref{lem: increasing and decreasing of d_n}. Take a curve of class 
$n\mu-\lambda$ on $\partial S^3(K)$ and $\be$ is a push off of that curve into 
the interior of $S^3(K)$, which has linking number 1 with $\al$ (so $\be$ is 
closer to $\partial S^3(K)$ than $\al$). There is a surgery exact triangle 
associated to $\be$:
\begin{equation*}
\xymatrix{
\shi(-S^3(K),-\Ga_{n-1})\ar[rr]&&\shi(-Y_{-n}(K),-\Ga_{n-1})\ar[dl]\\
&\shi(-S^3(K),-\Ga_{\mu})\ar[ul]&
}
\end{equation*}
Here the $3$-manifold $Y_{-n}$ is obtained from $S^3$ by performing $0$-surgery along $\be$, with respect to the $\partial S^3(K)$-surface-framing. Inside $S^3$, the curve $\be$ is isotopic to $K$ and the $0$-surface-framing is $(-n)$-Seifert-framing. Hence $Y_{-n}\cong S^3_{-n}(K)$. Let $(Y_{-n}-N(K))$ be the complement of $K$ inside $Y_{-n}$. Note $(Y_{-n}-N(K))$ can also be obtained from $S^3(K)$ by performing $0$-surgery along $\be\subset{\rm int}(S^3(K))$. So we can also regard the suture $\Ga_{n-1}$ as on $\partial (Y_{-n}-N(K))$. Inside $Y_{-n}$, the original knot $K$ is an unknot, since it bounds a disk $D$ obtained by capping off the annulus cobounded by $K$ and $\be$ through the $0$-surgery. This disk D intersects $\partial S^3(K)$ along a curve of class $n\mu-\lambda$. Hence it has two intersection with the suture $\Ga_{n-1}$. Hence, we conclude that
$$\shi(-(Y_{-n}-N(K)),-\Ga_{n-1})\cong I^{\sharp}(-Y_{-n})=I^{\sharp}(-S^3_{-n}(K)).$$
Next, note the curve $\al$ is inside both $S^3(K)$ and $Y_{-n}$, so we can also 
perform a $0$-surgery along $\al$ and obtain a commutative diagram:
\begin{equation*}
\xymatrix{
\shi(-S^3(K),-\Ga_{n-1})\ar[rr]\ar[dd]^{F_{\al,n-1}}&&\shi(-Y_{-n}(K),-\Ga_{n-1})\ar[dl]\ar[dd]^{V_{-n}}\\
&\shi(-S^3(K),-\Ga_{\mu})\ar[ul]\ar[dd]&\\
\shi(-M,-\Ga_{n-1})\ar[rr]&&\shi(-S^3(U),-\Ga_{n-1})\ar[dl]\\
&\shi(-M,-\Ga_{\mu})\ar[ul]&
}
\end{equation*}
Here, $M=S^1\times S^2(K)$ is a solid torus as above,
$$\shi(-M,-\Ga_{n-1})\cong\mathbb{C},$$
and
$$\shi(-M,-\Ga_{\mu})=0$$
as in the proof of \fullref{lem: increasing and decreasing of d_n}. Note $\be$ 
is closer to $\partial{S^3(K)}$ than $\al$, so inside $S^3$ $\al$ is 
also a meridian of $\be$. Hence performing a $0$-surgery along $\al$ inside 
$Y_{-n}=S^3_{-n}(K)$ makes it to become $S^3$ again. The knot $K\subset Y_{-n}$ 
bounds a disk $D$ which is disjoint from $\al$ so inside $S^3$ after performing 
the $0$-surgery, the knot still bounds the disk $D$, so in the above diagram it 
is written as $U\subset S^3$ to emphasize the fact that it is an unknot. Again 
$D$ intersects $\Ga_{n-1}$ twice so
$$\shi(-S^3(U),-\Ga_{n-1})\cong I^{\sharp}(-S^3)\cong\mathbb{C}.$$
Thus we have a simplified commutative diagram
\begin{equation*}
\xymatrix{
\shi(-S^3(K),-\Ga_{n-1})\ar[rr]\ar[dd]^{F_{\al,n-1}}&&\shi(-Y_{-n}(K),-\Ga_{n-1})=I^{\sharp}(-S^{3}_{-n}(K))\ar[dd]^{V_{-n}}\\
&&\\
\shi(-M,-\Ga_{n-1})\ar[rr]^{\cong}&&\shi(-S^3(U),-\Ga_{n-1})=I^{\sharp}(-S^3)
}
\end{equation*}
Hence we know that $V_{-n}\neq0$ since $F_{\al,n-1}\neq 0$. Next, note that 
$V_n$ is the cobordism map associated to the $0$-surgery along $\al$, so we 
know that
$$V_{-n}=W_{-n}^{\vee},$$
the dual of the map
$$W_{-n}\colon I^{\sharp}(S^3)\ra I^{\sharp}(S^3_{-n}(K)).$$
Hence $V_{-n}\neq0$ implies that $W_{-n}\neq0$. From \cite{baldwin2020framed}, we know that if $\nushI(K)\neq0$, then $K$ is 
$V$-shaped and $W_{-n}\neq0$ if and only if $-n<\nushI(K)$. Hence we conclude 
that
\[
  -2\tauI(\mir{K})-2=-n<\nushI(K). \qedhere
\]
\end{proof}

\begin{corollary}
  \label{cor: upper bound of nu^sharp}
For any knot $K$ so that $\nushI(K)\neq0$, we have
$$2\tauI(K)-2<\nushI(K)<2\tauI(K)+2.$$
\end{corollary}

\begin{proof}
Since $\tau_I$ is a concordance homomorphism, we know that
$$\tauI(K)=-\tauI(\mir{K}).$$
Hence
$$2\tauI(K)-2<\nushI(K)$$
follows direct from \fullref{lem: lower bound of nu^sharp}. For the other 
inequality, observe that
$$\nushI(K)=-\nushI(\mir{K}).$$
So we can apply the above inequality to $\mir{K}$ and have
\[
  -2\tauI(K)-2<-\nushI(K). \qedhere
\]
\end{proof}

We are now ready to identify $\tauI$ with $\taushI$:

\begin{proof}[Proof of \fullref{thm:main} in the instanton setting]
If $\taushI(K)\geq 1$, then we know that
$$\taushI(\bigconnsum nK)\geq n,$$
and hence
$$\nushI(\bigconnsum nK)>0$$
for large enough $n$. So we can apply \fullref{cor: upper bound of nu^sharp} to 
conclude that
$$\taushI(K)=\tauI(K).$$
If $\taushI(K)\leq 1$, then we can pick a knot $K_0$ with 
$\taushI(K_0)>0$. Since $\taushI$ is a concordance homomorphism, we 
know that
$$\taushI(K)=\taushI(K\connsum nK_0)-\taushI(\bigconnsum nK_0).$$
We can take $n$ large enough so that
$$\taushI(K\connsum nK_0)>0.$$
Hence from the above argument we know that
$$\taushI(K\connsum nK_0)=\tauI(K\connsum nK_0)$$
and
$$\taushI(\bigconnsum nK_0)=\tauI(\bigconnsum nK_0).$$
Since $\tau_I$ is also a concordance homomorphism, we conclude that
\[
  \taushI=\tauI(K). \qedhere
\]
\end{proof}

Finally, we state the necessary changes for the monopole setting:

\begin{proof}[Proof of \fullref{thm:main} in the monopole setting]
  The proof is similar to that in the instanton setting, with the symmetry 
  isomorphism
  \[
    \SHI (-S^3 (K), \Ga_n, -S_n, i) \isom \SHI (-S^3 (K), -\Ga_n, S_n, i)
  \]
  replaced by
  \[
    \SHM (-S^3 (K), \Ga_n, -S_n, i) \isom \SHM (-S^3 (K), -\Ga_n, S_n, i),
  \]
  for example, in \fullref{cor: minimal non vanishing grading}. While the 
  former follows from a symmetry in the generalized eigenspaces associated to 
  $\mu (R)$, the latter follows from the conjugation symmetry in the $\Spinc$ 
  decomposition in monopole Floer theory.
\end{proof}

%% file: sec-computation.tex
\section{Computation for twist knots}
\label{sec:twist-knots}

In this section, we compute $\khg$ for the family of knots $K_m$ as in 
\fullref{fig:twist-knots}, proving \fullref{thm:twist-knots}. We divide 
\fullref{thm:twist-knots} into four propositions: \fullref{prop: khg of K_m, 
  m<=0}, \fullref{prop: khg of K_m, m>0}, \fullref{prop: khg of K_m-mir, m<=0}, 
and \fullref{prop: khg of K_m-mir, m>0}

Note that, in particular, $K_1$ is the right-handed trefoil, $K_0$ is the 
unknot, and $K_{-1}$ is the figure-eight knot.

From the Seifert algorithm, we can easily construct a genus-$1$ Seifert surface 
for $K_m$, which we denote by $S_m$. Hence, $g(K_m)=1$.  Also,
it is straightforward to compute the (symmetrized) Alexander polynoial of $K_m$ 
to be
\begin{equation}\label{eq: alexander polynomial}
\Delta_{K_m}(t)=mt+(1-2m)+mt^{-1}.
\end{equation}

First, we will compute $\KHG (-S^3,K_m)$. Suppose $(S^3(K_m),\Ga_{\mu})$ is the 
balanced sutured manifold obtained by taking meridional sutures on knot 
complements. There is a curve $\al\subset \Int S^3(K_m)$ as in 
\fullref{fig:twist-knots} so that we have a surgery exact triangle:
\begin{equation*}
\xymatrix{
\shg(-S^3(K_m),-\Ga_{\mu})\ar[rr]&&\shg(-S^3(K_{m+1}),-\Ga_{\mu})\ar[dl]\\
&\shg(-M,-\Ga_{\mu})\ar[lu]&
}	
\end{equation*}

Here, $K_{m}$ is described as above, and $M$ is obtained from $S^3(K_m)$ by performing a 0-Dehn surgery along $\al$. We can use the surface $S_{m,\mu}$ which intersects the suture $\Ga_{\mu}$ twice to construct a grading on the sutured monopole and instanton Floer homologies. Let $S_{m,n}$ be an isotopy of $S_{m,\mu}$ so that $S_{m,n}$ intersects the suture $\Ga_n$ exactly $2n$ times. Since $\al$ is disjoint from $S_{m,\mu}$, all the Seifert surfaces $S_{m,\mu}$ and $S_{m,n}$ survives in $M$, which we call $S_{\mu}$ and $S_n$ respectively. Also, there is a graded version of the exact triangle (note we omit the surfaces from the following exact triangle): 

\begin{equation}\label{eq: graded surgery exact triangle along alpha}
  \vcenter{
\xymatrix{
\shg(-S^3(K_m),-\Ga_{\mu},i)\ar[rr]&&\shg(-S^3(K_{m+1}),-\Ga_{\mu},i)\ar[dl]\\
&\shg(-M,-\Ga_{\mu},i)\ar[lu]&
}
}
\end{equation}

Since $S_{m,\mu}$ has genus one and intersects the suture twice, all the graded 
sutured monopole and instanton Floer homologies in \eqref{eq: graded surgery 
  exact triangle along alpha} could only possibly be non-trivial for $-1\leq 
i\leq 1$. To understand what is $\shg(-M,-\Ga_{\mu})$, from \cite{kronheimer2010instanton} and \cite{li2018contact}, the surgery 
exact triangle \eqref{eq:bypass} is just the same as the oriented skein exact 
triangle and $\shg(-M,-\Ga)$ is isomorphic to the knot monopole or instanton 
Floer homology of the oriented smoothing of $K_m$, which is a Hopf link.  
Applying oriented Skein relation again on Hopf links, we can conclude that
\begin{equation}\label{eq: rank at most 4}
\rk_\ring(\shg(-M,-\Ga_{\mu}))\leq 4.
\end{equation}

For the monopole and instanton knot Floer homologies of $K_1$ (trefoil), we 
could look at the surgery exact triangle along the curve $\be$ in \fullref{fig: 
  trefoil} and argue in the same way as in \cite{kronheimer2010instanton} to conclude
$$\rk_\ring(\shg(-S^3(K_1),-\Ga_{\mu})\leq 3.$$
Using the Alexander polynomial in \eqref{eq: alexander polynomial} and 
\cite{kronheimer2010knots,kronheimer2010instanton}, we know that
\begin{equation}\label{eq: shg of trefoil}
\shg(-S^3(K_1),-\Ga_{\mu},S_{1,\mu},i)\cong\mathcal{R}
\end{equation}
for $i=-1,0,1$ and it vanishes in all other gradings.

\begin{figure}[htbp]
  \labellist
  \small\hair 2pt
  \pinlabel {\textcolor{red}{$\beta$}} [br] at 44 138
  \endlabellist
  \includegraphics[width=0.3\textwidth]{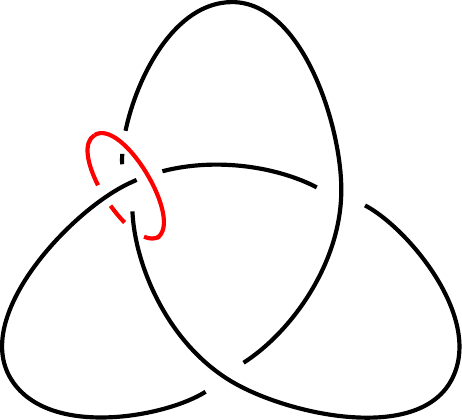}
  \caption{The trefoil and the circle $\be$.}
  \label{fig: trefoil}
\end{figure}

Now let $m=1$ in \eqref{eq: graded surgery exact triangle along alpha}. We know 
from \eqref{eq: alexander polynomial} that
$$\rk_\ring(\shg(-S^3(K_2),-\Ga_{\mu})\geq 7.$$
Then, from the exactness and inequalities \eqref{eq: rank at most 4} and 
\eqref{eq: shg of trefoil}, we know that 

\begin{equation*}
\rk_\ring(\shg(-M,-\Ga_{\mu}))=4.
\end{equation*}

After further examining each gradings, we know that
\begin{equation*}
  \label{eq: shg of (M,gamma_mu)}
\shg(-M,-\Ga_{\mu},S_{\mu},i)=
  \begin{cases}
\mathcal{R}&\text{for }i=1,-1,\\
\mathcal{R}^2&\text{for }i=0,\\
0&\text{otherwise.}
\end{cases}
\end{equation*}

Thus, by using the same argument and the induction, we can compute, for $m>0$, that
\begin{equation}\label{eq: khg tilde of K_m, m positive}
\shg(-S^3(K_m),-\Ga_{\mu},S_{m,\mu},i)=
  \begin{cases}
\mathcal{R}^m& \text{for }i=1,-1,\\
\mathcal{R}^{2m-1}& \text{for }i=0,\\
0& \text{otherwise.}
\end{cases}
\end{equation}

Since $K_0$ is the unknot, we can use the same technique to compute for, $m\leq0$, that
\begin{equation*}
  \label{eq: khg tilde of K_m, m non positive}
\shg(-S^3(K_m),-\Ga_{\mu},S_{m,\mu},i)=
\begin{cases}
\mathcal{R}^{-m}&\text{for }i=1,-1,\\
\mathcal{R}^{1-2m}&\text{for }i=0,\\
0&\text{otherwise.}
\end{cases}
\end{equation*}

Now we are ready to compute the minus version. Recall that the Seifert surface 
induces a framing on the boundary of the knot complements as well as $M$. Write 
$\Ga_n$ the suture consists of two curves of slope $-n$. We have a graded 
version of by-pass exact triangles \eqref{eq:bypass-gr-even} for even $n$ as 
well as
\eqref{eq:bypass-gr-odd}
for odd $n$.

A simple case to analyze is when $m<0$. For the knot $K_{m}$ with $m<0$, take 
$n=1$ in \eqref{eq:bypass-gr-odd}, we have \fullref{tb: psi +, 1}.
\begin{figure}[htbp]
\begin{equation*}	
\xymatrix@R-1pc{
  \gr&{}\save[]+<0in,1em>*+{}\ar@{-}[dddd]+<0in,-1em>\restore&\Ga_{\mu}\ar[r]^{\psi_{+,1}^{\mu}}&\Ga_1\ar[r]^{\psi_{+,2}^{1}}&\Ga_2\\
1&&\mathcal{R}^{-m}&&\mathcal{R}^{-m}\\
0&&\mathcal{R}^{1-2m}&\mathcal{R}^{-m}&\mathcal{R}^{b}\\
-1&&\mathcal{R}^{-m}&\mathcal{R}^{a}&\mathcal{R}^{c}\\
-2&&&\mathcal{R}^{-m}&\mathcal{R}^{-m}
}
\end{equation*}
\caption{The map $\psi_{+,2}^1$ for $K_m$. Each row is the positive bypass 
  exact triangle in a particular grading.  The leftmost column indicates the 
  gradings. We use letters like $a$, $b$, and $c$ to indicate that, \emph{a 
    priori,} we don't know what the rank is.}\label{tb: psi +, 1}
\end{figure}

Here, as in \cite[Section 4]{li2019direct}, the top and bottom non-vanishing grading of $\shg(-S^3(K_m),-\Ga_n)$ can be computed via sutured manifold decomposition and coincide with the top and bottom non-vanishing grading of $\shg(-S^3(K_m),-\Ga_{\mu}).$

From the graded exact triangles on the rows of the table and an extra exact 
triangle \eqref{eq: exact triangle to 3-sphere}, we know that
$$b\geq 1-m, \qquad c\geq a+m, \qquad b+c\leq a+1.$$
Hence, the only possibility is $b=1-m$, $c=a+m$. Now take $n=2$ in 
\eqref{eq:bypass-gr-even}, we have \fullref{tb: psi +, 2}.

\begin{figure}[htbp]
\begin{equation*}	
\xymatrix@R-1pc{
  \gr&{}\save[]+<0in,1em>*+{}\ar@{-}[ddddd]+<0in,-1em>\restore&\Ga_{\mu}\ar[r]^{\psi_{+,2}^{\mu}}&\Ga_2\ar[r]^{\psi_{+,3}^{2}}&\Ga_3\\
2&&\mathcal{R}^{-m}&&\mathcal{R}^{-m}\\
1&&\mathcal{R}^{1-2m}&\mathcal{R}^{-m}\ar[r]^{(\psi_{+,3}^2)_1}&\mathcal{R}^{1-m}\\
0&&\mathcal{R}^{-m}&\mathcal{R}^{1-m}\ar[r]^{(\psi_{+,3}^2)_0}&\mathcal{R}\\
-1&&&\mathcal{R}^{a+m}&\mathcal{R}^{a+m}\\
2&&&\mathcal{R}^{-m}&\mathcal{R}^{-m}
}
\end{equation*}
\caption{The map $\psi_{+,3}^2$ for $K_m$. We denote by $(\psi_{+,3}^2)_i$ the 
  restriction of the map $\psi_{+,3}^2$ to the grading
  $i$.}\label{tb: psi +, 2}
\end{figure}

Here, $\shg(-S^3(K_m),-\Ga_3)$ can be computed by taking $k=1$ in \eqref{eq: 
  the homology group for large n, 1}. We know from \cite[Section 5]{li2018gluing} that
$$\khg(-S^3,K_m,i)\cong \shg(-S^3(K_m),-\Ga_3,S_{m,3},i+1)$$
for $i=1,0,$ and $-1$, and the $U$ maps on $\khg(-S^3,K_m,i)$ for $i=1$ and 
$2$ coincide with the maps $\psi^i_{+,2}$ as in \fullref{tb: psi +, 2}. From 
the exactness, we know that $U$ map is actually zero at grading 1 and has a 
kernel of rank $-m$ at grading 0. Hence, we conclude:

\begin{proposition}
  \label{prop: khg of K_m, m<=0}
Suppose $m\leq 0$ and the knot $K_m$ is described as above. Then
$$\khg(-S^3,K_m)\cong \mathcal{R}[U]_{0}\oplus(\mathcal{R}_1)^{-m}\oplus(\mathcal{R}_0)^{-m},$$
and hence
\(
  \tau_G(K_m)=0.
\) \qed
\end{proposition}

To compute $\khg$ of $K_m$ for $m>0$, we first deal with the case $m=1$. Now 
$K_1$ is a right-handed trefoil,
which has
$\maxtb(K_1)=1$, and, hence, from \fullref{lem: the growth rate}, we know that
$$\rk_\ring\shg(-S^3(K_1),-\Ga_1)=\rk_\ring\shg(-S^3(K_1),-\Ga_0)+1.$$
Now let us compute $\shg(-S^3(K_1),-\Ga_0)$. Pick $S_0$ to be a genus 1 Seifert surface of $K$ so that $S_0$ is disjoint from $\Ga_0$. We can use the surface $S_0^-$, a negative stabilization of $S_0$ as in \cite[Definition 3.1]{li2019direct} to construct a grading on $\shg(-S^3(K_1),-\Ga_0)$. From the construction of grading and the adjunction inequality, there could only be three non-vanishing grading $-1,0,$ and $1$. For the grading $1$ part, we can apply \cite[Lemma 3.2 and Lemma 4.2]{li2019direct} and get
$$\shg(-S^3(K_1),-\Ga_0,S_0^-,1)\cong\shg(M',\ga'),$$
where  the balanced sutured manifold $(M',\ga')$ is obtained from $(-S^3(K),-\Ga_0)$ by a (sutured manifold) decomposion along the surface $S_0$
Since $K$ is a fibred knot, the underlining manifold $M'$ is just a product $[-1,1]\times S_0$. The suture $\ga'$ is not just $\{0\}\times\partial{S}$ but is actually three parallel copies of $\{0\}\times\partial{S}$ on $[-1,1]\times\partial{S}$. We can find an annulus $A\subset [-1,1]\times \partial{S}$ which contains the suture $\ga'$. Then, we can push the interior of $A$ into the interior of $S\times[-1,1]$ and get a properly embedded surface. If we further decompose $(M',\ga')$ along (the pushed off of) $A$, then we get a disjoint union of a product balanced sutured manifold $(S\times[-1,1],\partial{S}\times\{0\})$ with a solid torus with four longitudes as the suture. The sutured monopole and instanton Floer homologies of the first are both of rank 1 and the second of rank 2, as in \cite{kronheimer2010instanton} and \cite{li2018contact}. Hence, we conclude
$$\shg(-S^3(K_1),-\Ga_0,S_0^-,1)\cong\mathcal{R}^2.$$
For the other two gradings, note that from the grading shifting property in \cite[Proposition 4.3]{li2019direct}, we have
\begin{align*}
\shg(-S^3(K_1),-\Ga_0,S_0^-,i)&=\shg(-S^3(K_1),-\Ga_0,S_0^+,i-1)\\
&=\shg(-S^3(K_1),-\Ga_0,(-S_0)^-,1-i).
\end{align*}
The second equality follows from the basic observation that if we reverse the orientation of the surface $S_0^+$, then we get $(-S_0)^-$.  Hence,
$$\shg(-S^3(K_1),-\Ga_0,S_0^-,-1)=\shg(-S^3(K_1),-\Ga_0,(-S_0)^-,2)=0$$
by the adjunction inequality and
$$\shg(-S^3(K_1),-\Ga_0,S_0^-,0)=\shg(-S^3(K_1),-\Ga_0,(-S_0)^-,1)\cong\mathcal{R}^2.$$
by the same argument as above. Thus, as a conclusion,
$$\shg(-S^3(K_1),-\Ga_1)\cong\mathcal{R}^5.$$
Similarly, there are only three possible non-vanishing gradings 
$-1,0,1$. We have already known that the homology at top and bottom gradings 
are of rank 1 each, so the middle grading has rank $3$. Let $n=1$ in 
\eqref{eq:bypass-gr-odd}; we have \fullref{tb: psi +, -, 1}.

\begin{figure}[htbp]
\begin{equation*}	
\xymatrix@R-1pc{
  \gr&{}\save[]+<0in,1em>*+{}\ar@{-}[ddddd]+<0in,-1em>\restore&\Ga_{\mu}\ar[r]^{\psi_{+,1}^{\mu}}&\Ga_1\ar[r]^{\psi_{+,2}^{1}}&\Ga_2&{}\save[]+<0in,1em>*+{}\ar@{-}[ddddd]+<0in,-1em>\restore&\Ga_{\mu}\ar[r]^{\psi_{-,1}^{\mu}}&\Ga_{1}\ar[r]^{\psi_{-,2}^{1}}&\Ga_2\\
2&&&&&&&\mathcal{R}&\mathcal{R}\\
1&&\mathcal{R}&&1&&\mathcal{R}&\mathcal{R}^3&\mathcal{R}^b\\
0&&\mathcal{R}&\mathcal{R}&\mathcal{R}^b&&\mathcal{R}&\mathcal{R}&\mathcal{R}^c\\
-1&&\mathcal{R}&\mathcal{R}^3&\mathcal{R}^c&&\mathcal{R}&&\mathcal{R}\\
-2&&&\mathcal{R}&\mathcal{R}&&&&
}
\end{equation*}
\caption{The map $\psi^1_{+,2}$ (on the left) and $\psi^1_{-,2}$ (on the right) 
  for $K_1$.}\label{tb: psi +, -, 1}
\end{figure}
From the exactness, we know that $b=c=2$. The rest of the computation is straightforward and we conclude that
\begin{equation*}
  \label{eq: khg of right handed trefoil}
\khg(-S^3,K_1)\cong \mathcal{R}[U]_{1}\oplus\mathcal{R}_{0}.
\end{equation*}

Now we have the map
\[
C_{1,h,1}\colon \shg(-S^3(K_1),-\Ga_1)\ra\shg(-S^3(1),\delta)
\]
and by the description of $\khg(-S^3,K_1)$ above, \fullref{prop: another 
  definition of tau}, and the fact that $C_{1,h,n}$ commutes with $\psi_{-,n}$ 
(Claim~1 in the proof of \fullref{prop:tau-conc}), we know that
$$C_{1,h,1}\colon \shg(-S^3(K_1),-\Ga_1,1)\ra\shg(-S^3(1),-\delta)$$
is surjective, and, since $\shg(-S^3(K_1),-\Ga_1,1)$ has rank 1 it is actually 
an isomorphism (for the monopole case, the argument is essentially the same as 
in the proof of \fullref{lem: only rank is important}). Now we go back to the 
surgery exact triangle in \eqref{eq: graded surgery exact triangle along 
  alpha}, which corresponds to surgeries on the curve $\al\subset\Int 
S^3(K_m)$.  Since $\al$ is disjoint from the boundary, and as above, disjoint 
from all Seifert surfaces $S_{m,n}$, we have the following exact triangle for 
any $m$ and $n$
(where we again omit the surfaces):
\begin{equation}\label{eq: graded surgery exact triangle along alpha, Gamma_n}
  \vcenter{
\xymatrix{
\shg(-S^3(K_m),-\Ga_{n},i)\ar[rr]&&\shg(-S^3(K_{m+1}),-\Ga_{n},i)\ar[dl]\\
&\shg(-M,-\Ga_{n},i)\ar[lu]&
}
}
\end{equation}

There are contact $2$-handle attaching maps
$$C_{m,h,n}\colon \shg(-S^3(K_m),-\Ga_n)\ra\shg(-S^3(1),-\delta),$$
where the contact $2$-handle is attached along a meridional curve on the knot 
complements. We can attach a contact $2$-handle along the same curve on the 
boundary of M, and the handle attaching maps commute with the maps in the exact 
triangle \eqref{eq: graded surgery exact triangle along alpha, Gamma_n}. Thus, 
we have a diagram:
\begin{equation}\label{eq: graded surgery exact triangle plus commutative diagram}
  \vcenter{
\xymatrix{
\shg(-S^3(K_m),-\Ga_{n},i)\ar[rr]\ar[dd]^{C_{m,h,n}}&&\shg(-S^3(K_{m+1}),-\Ga_{n},i)\ar[dl]^{\tau_{m,n,i}}\ar[dd]^{C_{m+1,h,n}}\\
&\shg(-M,-\Ga_{n},i)\ar[lu]\ar[dd]^{\substack{\quad\quad\\\quad\quad\\C_{M,h,n}}}&\\
\shg(-S^3(1),-\delta)\ar[rr]^{\phi_{\infty}\quad\quad\quad\quad\quad}&&\shg(-S^3(1),-\delta)\ar[dl]^{\phi_{1}}\\
&\shg(-S^2\times S^1(1),-\delta)\ar[ul]^{\phi_{0}}&
}
}
\end{equation}
Here, $S^2\times S^1$ is obtained from $S^3$ by performing a $0$-surgery along the unknot $\al$. The balanced sutured manifold $(S^2\times S^1(1),\delta)$ is obtained from $S^2\times S^1$ by removing a $3$-ball and assigning a connected simple closed curve on the spherical boundary as the suture. Its sutured monopole and instanton Floer homologies are computed in \cite{baldwin2016instanton} and \cite{li2018contact} and are both of rank 2. Thus, the exactness tells us that $\phi_{\infty}=0$, $\phi_1$ is injective, and $\phi_0$ is surjective. 

Now take $m=0,n=1,$ and $i=1$, we know that
$$\shg(-M,-\Ga_1,S_1,1)\cong\shg(-S^3(K_1),-\Ga_1,S_{1,1},1)\cong\mathcal{R},$$
and $C_{M,h,n}$ is injective. Then, take $m$ to be an arbitrary non-negative 
integer and $n=1,i=1$ in \eqref{eq: graded surgery exact triangle plus 
  commutative diagram}. From \eqref{eq: khg tilde of K_m, m positive}, we know 
that
$$\shg(-S^3(K_m),-\Ga_{\mu},S_{m,\mu},1)\cong\mathcal{R}^m.$$
By performing sutured manifold decompositions along $S_{m,n}$ and applying 
\cite[Lemma 4.2]{li2019direct}, we know that
$$\shg(-S^3(K_m),-\Ga_{1},S_{m,1},1)\cong\shg(-S^3(K_m),-\Ga_{\mu},S_{m,\mu},1)\cong\mathcal{R}^m.$$
Recall from above discussions we have
$$\shg(-M,-\Ga_1,S_1,1)\cong\mathcal{R},$$
so in the exact triangle \eqref{eq: graded surgery exact triangle plus 
  commutative diagram}, we know that $\tau_{m,1,1}$ is surjective. Then, we can 
use the commutativity part of \eqref{eq: graded surgery exact triangle plus 
  commutative diagram} and conclude that
$$C_{m+1,h,n}\colon 
\shg(-S^3(K_{m+1}),-\Ga_1,S_{m,1},1)\ra\shg(-S^3(1),-\delta)$$
is surjective. From the fact that $\psi_{\pm,n+1}^n$ commutes with $C_{h,n}$ as in 
Claim 1 and 2 in the proof of \fullref{prop:tau-conc}, we know that this 
surjectivity means that the unique $U$ tower in $\khg(-S^3,K_{m},p_m)$ starts 
at grading 1:
$$\tau_G(K_m)=1$$
for $m>0$. 

Take $n=1$ in \eqref{eq:bypass-gr-odd}, then we have \fullref{tb: psi +, 1, 
  positive m}.

\begin{figure}[htbp]
\begin{equation*}	\xymatrix@R-1pc{
    \gr&{}\save[]+<0in,1em>*+{}\ar@{-}[dddd]+<0in,-1em>\restore&\Ga_{\mu}\ar[r]^{\psi_{+,1}^{\mu}}&\Ga_1\ar[r]^{\psi_{+,2}^{1}}&\Ga_2\\
1&&\mathcal{R}^m&&\mathcal{R}^m\\
0&&\mathcal{R}^{2m-1}&\mathcal{R}^m\ar[r]^{(\psi_{+,2}^1)_0}&\mathcal{R}^b\\
-1&&\mathcal{R}^m&\mathcal{R}^a&\mathcal{R}^c\\
-2&&&\mathcal{R}^m&\mathcal{R}^m
}
\end{equation*}
\caption{The map $\psi_{+,2}^1$ for $K_m$.}\label{tb: psi +, 1, positive m}
\end{figure}

The fact that $\tau_G(K_m)=1$ means that $(\psi_{+,2}^1)_0\neq0$, as $(\psi^1_{+,2})_0$ corresponds to the $U$ map at grading $1$ part of $\khg(-S^3,K_m,p_m)$. Thus, from the exactness we know that
$$b\geq m+1,~c\geq a-m.$$
From the exact triangle \eqref{eq: exact triangle to 3-sphere} we know that
$$b+c\leq a+1$$
and hence $b=m+1,c=a-m$. Thus, we conclude:

\begin{proposition}
  \label{prop: khg of K_m, m>0}
Suppose $m>0$ and $K_m$ is as above. Then
$$\khg(-S^3,K_m)\cong\mathcal{R}[U]_1\oplus(\mathcal{R}_1)^{m-1}\oplus(\mathcal{R}_0)^{m},$$
and hence $\tau_G(K_m)=1$. \qed
\end{proposition}

We could also compute the $\khg$ of the knots $\mir{K}_m$, the mirror image of 
$K_m$. For $m\leq0$, the computation is exactly the same as before, and we 
conclude:

\begin{proposition}
  \label{prop: khg of K_m-mir, m<=0}
Suppose $m\leq0$ and the knot $\mir{K}_m$ is as above. Then
$$\khg(-S^3,\mir{K}_m)\cong 
\mathcal{R}[U]_{0}\oplus(\mathcal{R}_1)^{-m}\oplus(\mathcal{R}_0)^{-m},$$
and hence
\(
  \tau_G(\mir{K}_m)=0.
\) \qed
\end{proposition}

For $m>0$, we have a diagram similar to \eqref{eq: graded surgery exact 
  triangle plus commutative diagram}, as follows.

\begin{equation}\label{eq: graded surgery exact triangle plus commutative diagram, 2}
  \vcenter{
\xymatrix{
\shg(-S^3(\mir{K}_{m+1}),-\Ga_{n},i)\ar[rr]&&\shg(-S^3(\mir{K}_{m}),-\Ga_{n},i)\ar[dl]^{\tau_{m,n,i}}\\
&\shg(-M,-\Ga_{n},i)\ar[lu]&
}}
\end{equation}

Let us first compute the case $m=1$, when $\bar{K}_m$ is the left-handed 
trefoil. In this case, take $n=1$ in \eqref{eq:bypass-gr-odd}, then we get 
\fullref{tb: psi +, 1, mirror case}.
\begin{figure}[htbp]
\begin{equation*}	
\xymatrix@R-1pc{
  \gr&{}\save[]+<0in,1em>*+{}\ar@{-}[dddd]+<0in,-1em>\restore&\Ga_{\mu}\ar[r]^{\psi_{+,1}^{\mu}}&\Ga_1\ar[r]^{\psi_{+,2}^{1}}&\Ga_2\\
1&&\mathcal{R}&&\mathcal{R}^m\\
0&&\mathcal{R}&\mathcal{R}^m\ar[r]^{(\psi_{+,2}^1)_0}&\mathcal{R}^b\\
-1&&\mathcal{R}^m&\mathcal{R}^a&\mathcal{R}^c\\
-2&&&\mathcal{R}&\mathcal{R}
}
\end{equation*}
\caption{The map $\psi_{+,2}^1$ for $\mir{K}_1$.}\label{tb: psi +, 1, mirror case}
\end{figure}

The left-handed trefoil is not right veering in the sense of 
\cite{baldwin2018khovanov}, so from their discussion we conclude that 
$(\psi_{+,2}^1)_0=0$. (This is how they prove that the second top grading of 
the instanton knot Floer homology of a non--right-veering knot is non-trivial.  
Though they only work in the instanton case, the monopole case is exactly the 
same.) Thus we conclude that $b=0$.

In \eqref{eq: graded surgery exact triangle plus commutative diagram, 2}, let 
$m=0,n=2,i=0$. Note the grading is induced by $S_{m,2}^+$, i.e., a Seifert 
surface of the knot $\mir{K}_m$ which intersects the suture $\Ga_{2}$ 
transversely at four points and with a positive stabilization.  With the 
gradings as in the first row of \eqref{eq:bypass-gr-even}, we have
$$\shg(-S^3(\mir{K}_1),-\Ga_2,S_{1,2}^{\tau(2)},0)=\mathcal{R}^b=0, \qquad 
\shg(-S^3(\mir{K}_0),-\Ga_2,S_{0,2}^{\tau(2)},0)\cong \mathcal{R}.$$
Here, $\mir{K}_0$ is the unknot and we have computed the $\shg$ of a solid 
torus with any possible sutures in \cite{li2019direct}. Thus, we conclude that
$$\shg(-M,-\Ga_2,S_{2}^{\tau(2)},0)\cong\mathcal{R}.$$

Use the exactness and the induction, then we have
$$\shg(-S^3(\mir{K}_m),\Ga_2,0)\cong\mathcal{R}^{c_m}, \qquad c_m\leq m-1.$$

For the knot $\mir{K}_m$, take $n=2$ in \eqref{eq:bypass-gr-even}, then we have 
\fullref{tb: psi +, 2, mirror case}.
\begin{figure}[htbp]
\begin{equation*}	
\xymatrix@R-1pc{
  \gr&{}\save[]+<0in,1em>*+{}\ar@{-}[ddddd]+<0in,-1em>\restore&\Ga_{\mu}\ar[r]^{\psi_{+,2}^{\mu}}&\Ga_2\ar[r]^{\psi_{+,3}^{2}}&\Ga_3\\
2&&\mathcal{R}^m&&\mathcal{R}^m\\
1&&\mathcal{R}^{2m-1}&\mathcal{R}^m\ar[r]^{(\psi_{+,3}^2)_1}&\mathcal{R}^{c_m}\\
0&&\mathcal{R}^m&\mathcal{R}^{c_m}\ar[r]^{(\psi_{+,3}^2)_0}&\mathcal{R}\\
-1&&&?&?\\
2&&&\mathcal{R}^m&\mathcal{R}^m
}
\end{equation*}
\caption{The map $\psi_{+,3}^2$ for $K_m$.}\label{tb: psi +, 2, mirror case}
\end{figure}
Thus, we conclude from the exactness that $c_m=m-1$, $(\psi^2_{+,3})_1=0$, and 
$(\psi^2_{+,3})_0=0$. As above, the two maps $(\psi^2_{+,3})_1$ and 
$(\psi^2_{+,3})_0$ correspond to the $U$ maps of $\khg(-S^3,\mir{K}_m)$ at 
grading $1$ and $0$, respectively. Hence, we conclude:

\begin{proposition}
  \label{prop: khg of K_m-mir, m>0}
Suppose $m>0$ and the knot $\mir{K}_m$ is as above. Then
$$\khg(-S^3,\mir{K}_m)\cong 
\mathcal{R}[U]_{-1}\oplus(\mathcal{R}_1)^{m}\oplus(\mathcal{R}_0)^{m-1},$$
and hence
\(
  \tau_G(\mir{K}_m)=0.
\)
\qed
\end{proposition}